\documentclass[twoside, onecolumn, 10pt]{article}
\usepackage{latexsym}
\usepackage{amssymb}
\usepackage{amsmath}
\usepackage{marvosym}

\usepackage{mathrsfs}
\newcommand{\R}{{\mathscr R}}

\newtheorem{theorem}{Theorem}[section]
\newtheorem{lemma}[theorem]{Lemma}
\newtheorem{corollary}[theorem]{Corollary}

\def\nb{\nonumber}

\begin{document}

\begin{center}
{\LARGE How to solve three fundamental linear matrix inequalities in the L\"owner partial ordering}
\end{center}

\begin{center}
{\large Yongge  Tian$^*$}
\end{center}

\begin{center}
{\footnotesize {\it  CEMA, Central University of Finance and Economics, Beijing 100081, China}}
\end{center}

\renewcommand{\thefootnote}{\fnsymbol{footnote}}
\footnotetext{$^*$E-mail: yongge.tian@gmail.com}

\noindent {\bf Abstract.}
This paper shows how to solve analytically the three fundamental linear matrix inequalities
 $$
 AXB \succcurlyeq C \, (\succ  C), \ \ AXA^*\succcurlyeq B \, (\succ  B), \ \  AX +(AX)^{*} \succcurlyeq B \, (\succ  B)
 $$
in the L\"owner partial ordering by using ranks, inertias and generalized inverses of matrices.
\\


\medskip


\section{Introduction}

\renewcommand{\theequation}{\thesection.\arabic{equation}}
\setcounter{section}{1}
\setcounter{equation}{0}

Throughout this paper,
\begin{enumerate}
\item[]${\mathbb C}^{m\times n}$ stands for the set
of all $m\times n$  complex matrices;

\item[] ${\mathbb C}_{{\rm H}}^{m}$ stands for the sets of all $m\times m$ complex Hermitian matrices;

\item[] ${\mathbb C}_{{\rm SH}}^{m}$ stand for the sets of all $m\times m$ complex Hermitian matrices
 and complex  skew-Hermitian matrices;

\item[]  the symbols $A^{*}$,  $r(A)$ and ${\mathscr R}(A)$ stand for the conjugate transpose,
the rank and the range (column space) of a matrix
$A \in \mathbb C^{m \times n}$, respectively;

\item[]  $[\, A, \, B\,]$ denotes a row block matrix consisting of $A$ and $B$;

\item[]  the Moore--Penrose inverse of $A\in {\mathbb
C}^{m\times n}$, denoted by $A^{\dag}$, is defined to be the unique
solution $X$ satisfying the four matrix equations $AXA = A,$ $XAX =
X,$ $(AX)^{*} = AX$ and $(XA)^{*} =XA$;

\item[] the symbols $E_A$ and $F_A$ stand for $E_A = I_m - AA^{\dag}$ and $F_A = I_n-A^{\dag}A$, their ranks are
given by $r(E_A) = m - r(A)$ and $r(F_A) = n - r(A)$;

\item[] $i_{+}(A)$ and $i_{-}(A)$, called the partial inertia of  $A \in {\mathbb C}_{{\rm H}}^{m}$,
 are defined to be  the numbers of the positive and negative eigenvalues of $A$ counted with
multiplicities, respectively, where $r(A) = i_{+}(A) + i_{-}(A);$

\item[]  $A \succcurlyeq 0$ ($A \succ  0$) means that $A$ is Hermitian positive semi-definite (positive definite);

\item[]  two $A,  \, B \in {\mathbb C}_{{\rm H}}^{m}$ are said to satisfy the inequality $A \succcurlyeq B$
($A \succ  B)$ in the L\"owner partial ordering if $A - B$ is Hermitian positive semi-definite (positive definite);

\item[] a positive  semi-definite matrix $A$ of order $m$ is said to be a contraction if all its
eigenvalues are less then or equal to 1, i.e., $0\preccurlyeq A \preccurlyeq I_m$, to be a strict contraction if all
its eigenvalues are less then 1, i.e., $0\preccurlyeq A \prec I_m$.
\end{enumerate}
A well-known property of the Moore--Penrose inverse is $(A^{\dag})^* =
(A^*)^{\dag}$. In particular $AA^{\dag} =A^{\dag}A$ if $A = A^{*}$.
We shall repeatedly use them in the latter part of this paper. One
of the most important applications of generalized inverses is to
derive some closed-form formulas for calculating ranks and inertias of matrices,
as well as general solutions of matrix equations; see Lemmas
2.1--2.9 below. Results on the Moore--Penrose inverse can be found,
e.g., in \cite{BG,Ber,Ho}.

The L\"owner partial ordering for matrices,  as a natural extension of inequalities for real numbers,
is one of the most useful concepts in matrix theory for  characterizing relations between two
complex  Hermitian (real symmetric) matrices of the same size, while a main object of study in
core matrix theory is to compare Hermitian matrices in the L\"owner partial ordering and
to establish various possible matrix inequalities. This subject was extensively studied by
many authors, and numerous matrix inequalities in the L\"owner partial ordering were
established in the literature. In the investigation of the L\"owner partial ordering
between two Hermitian matrices, a challenging task is to solve matrix inequalities that
involve unknown matrices. This topic can generally be stated as follows: \\

\medskip

\noindent {\bf Problem.} \ For a given matrix-valued function $\phi(X)$ that satisfies $\phi(X) = \phi^*(X),$
 where $X$ is a variable matrix,  establish  necessary and sufficient conditions for the matrix inequality
\begin{equation}
\phi(X) \succcurlyeq  0,  \ \ \phi(X) \succ  0, \ \ \phi(X) \preccurlyeq 0, \ \ \phi(X) \prec  0
\label{10a}
\end{equation}
to hold, respectively, and find solutions $X$ of the matrix inequalities.

\medskip

A matrix-valued function for complex matrices is a map between matrix spaces, which can generally be written as
$Y = \phi(X)$ for $Y \in {\mathbb C}^{m \times n}$ and $X \in {\mathbb C}^{p \times q},$  or briefly,
$f:  {\mathbb C}^{m \times n} \rightarrow {\mathbb C}^{p \times q},$  where ${\mathbb C}^{m \times n}$
and ${\mathbb C}^{p \times q}$ are two two complex matrix spaces.  As usual, linear  matrix-valued
functions as common representatives of various matrix-valued functions are extensively
studied from theoretical and applied points of view.  When $\phi(X)$ in (\ref{10a}) is a linear matrix-valued function,
it is usually called a linear matrix inequality (LMI) in the literature. A systematic work on LMIs and their
applications in system and control theory can be found, e.g., in
\cite{BGFB,SIG}. LMIs in the L\"owner partial ordering are usually taken as convex constraints to unknown matrices
and vectors in mathematical programming and optimization
theory.

This paper aims at solving the following three groups of LMIs of fundamental type:
\begin{align}
& AXB \succcurlyeq C \,(\succ C, \ \preccurlyeq C, \  \prec C),
\label{10b}
\\
& AXA^{*} \succcurlyeq B \, (\succ  B, \ \preccurlyeq B, \  \prec B ),
\label{10c}
\\
& AX + (AX)^{*} \succcurlyeq B \, (\succ  B, \ \preccurlyeq B, \  \prec B).
\label{10d}
\end{align}
They are the simplest cases of various types of LMIs and are the starting point of many
advanced study on complicated LMIs.

Recall that any Hermitian nonnegative definite (positive definite) matrix
$M$ can be written as $M = UU^*$ for certain (nonsingular) matrix $U$. Hence,
the mechanism of a matrix inequality in the L\"owner partial ordering can
be explained by certain matrix equation that involves an unknown quadratic term.
In fact, any  matrix inequality $\phi(X) \succcurlyeq 0$ positive semi-definiteness
 (matrix inequality $\phi(X) \succ 0$ for positive definiteness) can equivalently be relaxed
 to
\begin{equation}
\phi(X) - UU^* =0
\label{10aa}
\end{equation}
for certain (nonsingular) matrix $U$. Due to the
non-commutativity of matrix algebra, there are no general methods
for finding analytical solutions of quadratic matrix equations, so
that it is hard to solve for the unknown matrices $X$ and $U$ from the
equation in (\ref{10aa}) for a general $\phi(X)$. However, for the
three fundamental LMIs in (\ref{10b})--(\ref{10d}), we are
able to establish their analytical solutions by using the relaxed
matrix equation in (\ref{10aa}), and ordinary operations of the given
matrices and their generalized inverses.

Matrix equations and matrix inequalities in the L\"owner partial ordering have been main objects
of study in matrix theory and their applications. Many new theories and
methods were developed in the investigations of matrix equations and inequalities.
In particular, the concept of generalized inverses of matrices was introduced when
Penrose considered general solutions of the matrix equations $AX  = B$ and $AXB = C$, cf. \cite{Pen}.
The three matrix equations  associated with (\ref{10b})--(\ref{10c}) are
\begin{align}
 AXB = C, \ \  AXA^{*}=  B, \ \ AX + (AX)^{*} = B,
\label{10f}
\end{align}
which were extensively studied from theoretical and practical points of view, while the three
matrix-valued functions
\begin{align}
& \phi_1(X) = C - AXB, \ \ \
\label{10g}
\\
& \phi_2(X)  = B - AXA^{*},
\label{10g1}
\\
& \phi_3(X) = B - AX - (AX)^{*}  \ \ \ 
\label{10h}
\end{align}
associated with (\ref{10b})--(\ref{10d}) were recently considered in
\cite{LT-nlaa,LT-jota,Ti-se,T-laa10,T-laa11,TC,TL}.  Because (\ref{10b})--(\ref{10h})
are some simplest cases of matrix equations, matrix inequalities
and matrix functions, they have been attractive objects of study in matrix theory and
applications. In fact, it is remarkable that simply knowing when the
 LMIs in (\ref{10b})--(\ref{10d}) are feasible  gives some deep insights into the relations
 between both sides of the LMIs.

This paper is organized as follows. In Section 2, we give a group of known results on matrix
equations, as well as some expansion formulas for calculating (extremal)
ranks and inertias of matrices.  In Section 3, we solve for the inequality in
(\ref{10b}), and discuss various algebraic properties of the LMI and its solution.
In particular, we shall give a group of closed-form formulas for calculating
the extremal ranks and inertias of $D - AXB$ subject to
$AXB \succcurlyeq C$, and use the formulas to  establish necessary and sufficient
conditions for the two-sided matrix inequality  $D \succcurlyeq AXB \succcurlyeq C$ to
be solvable. In Sections 4 and 5, we establish necessary and sufficient conditions
for the LMIs in (\ref{10c}) and (\ref{10d}) to be feasible, respectively, and derive
general solutions in closed-forms of these LMIs. In Section 7, we give a group of formulas
for calculating the extremal ranks and inertias of $A - BX - XB^*$ subject to $BXB^*=C$,
and use the formulas to characterize the existence of Hermitian
matrix $X$ that satisfies $BX + XB^* \succcurlyeq A\, ( BX + XB^* \succ  A)$
subject to $BXB^*=C$. Some further research problems are presented in Section 8.

\section{Preliminaries}
\renewcommand{\theequation}{\thesection.\arabic{equation}}
\setcounter{section}{2}
\setcounter{equation}{0}

In this section, we present some known or new results on solving matrix equations,
as well as formulas for calculating ranks and
inertias of matrices, which will be used in the latter part of
this paper.

\begin{lemma} [\cite{KM}]\label{T25}
Let $A, \, B\in \mathbb C^{m\times n}$ be given$.$ Then$,$ the following hold$.$
\begin{enumerate}
\item[{\rm(a)}]  The matrix equation
\begin{equation}
AX = B
\label{25}
\end{equation}
has a Hermitian solution $X \in \mathbb C_{{\rm H}}^{n}$  if and only if ${\mathscr R}(B) \subseteq
{\mathscr R}(A)$ and $AB^{*} =BA^{*}.$ In this case$,$ the general
Hermitian solution of {\rm (\ref{25})} can be written in the following parametric form
\begin{equation}
X = A^{\dag}B + (A^{\dag}B)^{*} - A^{\dag}BA^{\dag}A + F_AWF_A,
\label{26}
\end{equation}
where $W\in \mathbb C_{{\rm H}}^{n}$ is arbitrary$.$

\item[{\rm(b)}] The matrix equation
\begin{equation}
AXX^{*} = B
\label{27}
\end{equation}
has a solution for $XX^*$ if and only if ${\mathscr R}(B) \subseteq {\mathscr
R}(A),$ $AB^{*} \succcurlyeq 0$ and $r(AB^{*}) =r(B).$ In this case$,$ the
general solution of {\rm (\ref{27})} can be written  in the following parametric form
\begin{equation}
XX^{*} = B^{*}(AB^{*})^{\dag}B   + F_AWW^{*}F_A,
\label{28}
\end{equation}
where $W\in \mathbb C^{n \times n}$ is arbitrary$.$
\end{enumerate}
\end{lemma}

\begin{lemma}[\cite{Pen}] \label{T27}
Let $A \in \mathbb C^{m\times n},$  $B\in \mathbb C^{p\times q}$ and
$C\in \mathbb C^{m\times q}$ be given$.$ Then$,$ the matrix equation
\begin{equation}
AXB = C
\label{213}
\end{equation}
has a solution if and only if ${\mathscr R}(C) \subseteq {\mathscr
R}(A)$ and ${\mathscr R}(C^{*}) \subseteq {\mathscr R}(B^{*}),$ or
equivalently$,$ $E_AC = 0$ and $CF_B =0.$  In this case$,$ the
general solution of {\rm (\ref{213})} can be written in the
following parametric forms
\begin{align}
X &= A^{\dag}CB^{\dag} + W - A^{\dag}AWBB^{\dag}, \label{214}
\\
X &= A^{\dag}CB^{\dag} + F_AU_1 +  U_2E_B,
\label{214a}
\end{align}
respectively$,$ where $W, \,  U_1,  \, U_2 \in \mathbb C^{n \times p}$ are arbitrary$.$
\end{lemma}

\begin{lemma} \label{T26}
Let $A \in {\mathbb C}^{m \times n}$ and $B \in {\mathbb C}_{{\rm H}}^{m}$ be given$.$ Then$,$ the following hold$.$
\begin{enumerate}
\item[{\rm(a)}] {\rm \cite{Gr}} The matrix equation
\begin{equation}
AXA^{*} = B
\label{29}
\end{equation}
has a solution $X \in \mathbb C_{{\rm H}}^{n}$ if and only if
${\mathscr R}(B) \subseteq {\mathscr R}(A),$ or equivalently$,$
$AA^{\dag}B = B.$ In this case$,$ the general Hermitian solution of {\rm (\ref{29})} can be written
 in the following parametric forms
\begin{align}
X & = A^{\dag}B(A^{\dag})^{*} +  U - A^{\dag}AUA^{\dag}A,
\label{z42}
\\
X & = A^{\dag}B(A^{\dag})^{*} +  F_AV + V^{*}F_A,
\label{z43}
\end{align}
respectively$,$ where $U \in \mathbb C_{{\rm H}}^{n}$ and $ V \in \mathbb C^{n \times n}$ are arbitrary$.$

\item[{\rm(b)}] {\rm \cite{Gr,KM}} There exists an $X\in \mathbb C^{n\times n}$  such that
\begin{equation}
AXX^{*}A^{*} = B
\label{211pp}
\end{equation}
if and only if  $B \succcurlyeq 0$ and ${\mathscr R}(B) \subseteq {\mathscr R}(A).$ In this case$,$ the general solution of {\rm (\ref{211pp})}
can be written as
\begin{equation}
XX^{*} = 
A^{\dag}B(A^{\dag})^*  + F_AVB(A^{\dag})^* +  A^{\dag}BV^*F_A +  F_AWW^*F_A,
\label{212pp}
\end{equation}
where $V\in \mathbb C^{n \times m}$  and $W\in \mathbb C^{n \times n}$ are arbitrary$.$

\item[{\rm(c)}] {\rm \cite{Ba}} Under  $A, \, B \in {\mathbb C}^{m \times m},$ there exists an
$X\in \mathbb C^{m \times m}$  such that
\begin{equation}
AXX^{*}A^{*} = B
\label{211}
\end{equation}
if and only if  $B \succcurlyeq 0$ and ${\mathscr R}(B) \subseteq {\mathscr
R}(A).$ In this case$,$ the general solution of {\rm (\ref{211})}
can be written as
\begin{equation}
XX^{*} =(\, A^{\dag}B^{\frac{1}{2}} + F_AV\,)(\, A^{\dag}B^{\frac{1}{2}} + F_AV\,)^*,
\label{212}
\end{equation}
where $V\in \mathbb C^{m \times m}$ is arbitrary$.$
\end{enumerate}
\end{lemma}

\begin{lemma} [\cite{TL}]\label{T24}
Let $A \in \mathbb C^{m\times n}$ and $B\in \mathbb C_{{\rm
H}}^{m}$ be given$.$ Then$,$ the following hold$.$
\begin{enumerate}
\item[{\rm(a)}] There exists an $X \in \mathbb C^{n\times m}$ such that
\begin{equation}
AX + (AX)^{*} = B
\label{21}
\end{equation}
if and only if $E_ABE_A  =0.$ In this case$,$ the general solution of {\rm (\ref{21})} can be written
in the following parametric form
\begin{equation}
X = \frac{1}{2}A^{\dag}B( \, 2I_m - AA^{\dag} \,) + VA^{*} + F_AW,
\label{22}
\end{equation}
where both $V \in \mathbb C_{{\rm SH}}^{n}$ and $W\in \mathbb C^{n
\times m}$ are arbitrary$.$

\item[{\rm(b)}] There exists an $X \in \mathbb C^{n\times m}$ such that
\begin{equation}
AX + (AX)^{*} = BB^{*}
\label{23}
\end{equation}
if and only if ${\mathscr R}(B) \subseteq {\mathscr R}(A).$ In this case$,$ the general solution
of {\rm (\ref{23})} can be written as
\begin{equation}
X = \frac{1}{2} A^{\dag}BB^{*} + VA^{*} + F_AW,
\label{24}
\end{equation}
where both $V  \in \mathbb C_{{\rm SH}}^{n}$ and $W\in \mathbb C^{n
\times m}$ are arbitrary$.$
\end{enumerate}
\end{lemma}

\begin{lemma}\label{T25K}
  Let $A_1 \in {\mathbb C}^{m \times p}, \ B_1  \in
{\mathbb C}^{q \times n}, \ A_2 \in {\mathbb C}^{m \times r},$  $B_2 \in
 {\mathbb C}^{s \times n}$ and $C \in {\mathbb C}^{m \times n}$ be given$.$
Then$,$ the following hold$.$
\begin{enumerate}
\item[{\rm (a)}]{\rm \cite{OZ}}  There exist $X\in {\mathbb C}^{p\times q}$ and
$Y\in {\mathbb C}^{r\times s}$ such that
\begin{equation}
A_1XB_1 + A_2YB_2  = C
\label{15}
\end{equation}
if and only if the following four rank equalities
\begin{align}
r[ \, C,  \, A_1,\, A_2 \, ] = r[ \,  A_1, \,  A_2 \,  ], & \ \
r\!\left[\!\!\begin{array}{c} C \\ B_1 \\ B_2  \end{array} \!\!\right] =  r\!
\left[\!\! \begin{array}{c} B_1 \\ B_2 \end{array} \!\!\right]\!,
\label{16}
\\
r\!\left[\!\! \begin{array}{cc} C & A_1  \\ B_2 &  0 \end{array}\!\!\right]
 = r(A_1) + r(B_2),  & \ \
r \!\left[\!\! \begin{array}{cc}  C & A_2  \\ B_1 &  0
\end{array}\!\!\right] = r( A_2) + r(B_1)
\label{17}
\end{align}
hold$,$ or equivalently$,$
\begin{equation}
 [\,  A_1, \,  A_2 \,][\,A_1, \,  A_2 \,  ]^{\dag} C = C, \ \
 C\left[\!\! \begin{array}{c} B_1 \\ B_2 \end{array}
\!\!\right]^{\dag}\left[\!\! \begin{array}{c} B_1 \\ B_2
\end{array} \!\!\right] =C, \ \
E_{A_1}CF_{B_2} = 0,  \ \  E_{A_2}CF_{B_1} = 0. \
\label{18}
\end{equation}

\item[{\rm (b)}]{\rm \cite{Ti-lama00}}  Under {\rm(\ref{16})}
and {\rm (\ref{17})}$,$ the general solutions
 of {\rm (\ref{15})} can be decomposed as
\begin{equation}
 X =  X_0 + X_1X_2 +  X_3 \ \ and  \ \ Y =
Y_0 - Y_1Y_2 + Y_3,
 \label{19}
\end{equation}
where $X_0 $ and $ Y_0$ are a pair of special solutions of {\rm (\ref{15})}$,$
$X_1, \, X_2, \, X_3$ and $Y_1, \, Y_2, \, Y_3$ are the general solutions of
the following four homogeneous matrix equations
\begin{equation}
A_1X_1 + A_2Y_1 = 0,  \ \ X_2B_1 + Y_2B_2 = 0,  \ \
A_1X_3B_1 = 0, \ \  A_2Y_3B_2 = 0.
\label{110}
\end{equation}
By using generalized inverses of matrices$,$  {\rm (\ref{19})} can
be written  in the following parametric forms
\begin{align}
  X & = X_0 + [\,  I_p, \,  0 \, ]F_GWE_H
\!\left[\!\! \begin{array}{c} I_q  \\ 0
 \end{array} \!\!\right] + F_{A_1}W_1 + W_2E_{B_1},
\label{111}
\\
Y & = Y_0 - [\,  0, \,   I_r \, ] F_GWE_H \!\left[\!\!
\begin{array}{c} 0  \\ I_s
\end{array}\!\! \right] + F_{A_2}W_3 + W_4E_{B_2},
\label{112}
\end{align}
where $G = [ \, A_1, \, A_2 \, ], \ H = \left[\!\!\begin{array}{c} B_1
\\ B_2
\end{array} \!\!\right]\!,$ the five matrices  $W, \,  W_1, \, W_2, \, W_3$ and $W_4$ are arbitrary$.$
\end{enumerate}
\end{lemma}

Lemmas \ref{T25}--\ref{T25K} show that general solutions of some simple matrix equations can be written
as analytical forms composed by the given matrices and their generalized inverses, as well as
 arbitrary matrices. These analytical formulas can be easily used to establish various algebraic
 properties of the solutions of the equations, such as, their ranks, ranges, uniqueness, definiteness, etc.

In order to simplify various matrix expression involving generalized inverse of matrices and
arbitrary matrices, we need some formulas for ranks and inertias of matrices. The following is
 obvious from the definitions of rank and inertia.

\begin{lemma} \label{T20}
Let $A \in {\mathbb C}^{m\times m},$ $B \in {\mathbb C}^{m\times
n},$ and $C \in {\mathbb C}_{{\rm H}}^{m}.$ Then$,$ the following hold$.$
\begin{enumerate}
\item[{\rm (a)}] $A$ is nonsingular if and only if $r(A) = m.$

\item[{\rm (b)}] $B =0$ if and only if $r(B) = 0.$

\item[{\rm (c)}]  $C \succ 0$ $(C \prec 0)$ if and only if $i_{+}(C) = m$ $(i_{-}(C) = m)$,

\item[{\rm (d)}]  $C \succcurlyeq 0$ $(C \preccurlyeq 0)$ if and only if $i_{-}(C) = 0$ $(i_{+}(C) = 0)$.
\end{enumerate}
\end{lemma}

\begin{lemma}  \label{T27K}
Let ${\cal S}$ be a set consisting of  matrices over ${\mathbb
C}^{m\times n},$ and let ${\cal H}$ be a set consisting of Hermitian
matrices over ${\mathbb C}_{{\rm H}}^{m}.$ Then$,$ the following hold$.$
\begin{enumerate}
\item[{\rm (a)}]  Under $m =n,$ ${\cal S}$ has a nonsingular matrix if and only if
$\max_{X\in {\cal S}} r(X) = m.$

\item[{\rm (b)}] Under $m =n,$  all $X\in {\cal S}$ are nonsingular if and only if
 $\min_{X\in {\cal S}} r(X) = m.$

\item[{\rm (c)}] $0\in {\cal S}$ if and only if
$\min_{X\in {\cal S}} r(X) = 0.$

\item[{\rm (d)}] ${\cal S} = \{ 0\}$ if and only if
$\max_{X\in {\cal S}} r(X) = 0.$

\item[{\rm (e)}] ${\cal H}$ has a matrix $X \succ  0$  $(X \prec  0)$ if and only if
$\max_{X\in {\cal H}} i_{+}(X) = m  \ \left(\max_{X\in
{\cal H}} i_{-}(X) = m \right)\!.$

\item[{\rm (f)}] All $X\in {\cal H}$ satisfy $X \succ 0$ $(X \prec  0)$ if and only if
$\min_{X\in {\cal H}} i_{+}(X) = m \ \left(\min_{X\in {\cal H}}
i_{-}(X) = m\, \right)\!.$

\item[{\rm (g)}] ${\cal H}$ has a matrix  $X \succcurlyeq 0$ $(X \preccurlyeq 0)$ if and only if
$\min_{X \in {\cal H}} i_{-}(X) = 0 \ \left(\min_{X\in
{\cal H}} i_{+}(X) = 0 \,\right)\!.$

\item[{\rm (h)}] All $X\in {\cal H}$ satisfy $X \succcurlyeq 0$ $(X \preccurlyeq 0)$  if and only if
$\max_{X\in {\cal H}} i_{-}(X) = 0 \ \left(\max_{X\in {\cal H}}
i_{+}(\,X) = 0\, \right)\!.$
\end{enumerate}
\end{lemma}

The question of whether a given matrix function is nonnegative definite or positive definite
everywhere is ubiquitous in mathematics and applications. Lemma \ref{T27K}(e)--(h) show that if certain explicit formulas for calculating the global
extremal inertias of a given Hermitian matrix function are established,
we can use them, as demonstrated in Sections 2, 3 and 5 below, to derive necessary
and sufficient conditions for the Hermitian matrix function to be definite or semi-definite.

\begin{lemma} [\cite{MS}]\label{T28}
Let $A \in \mathbb C^{m \times n}, \ B \in \mathbb C^{m \times k},$
$C \in \mathbb C^{l \times n}$ and $D \in \mathbb C^{l \times k}$. Then$,$ the following hold$.$
the following rank expansion formulas hold
\begin{align}
r[\, A, \, B \,] & = r(A) + r(E_AB) = r(B) + r(E_BA),
\label{215}
\\
r \!\left[\!\! \begin{array}{c}  A  \\ C  \end{array}
\!\!\right] & = r(A) + r(CF_A) = r(C) + r(AF_C),
\label{216}
\\
r\!\left[\!\! \begin{array}{cc}  A  & B  \\ C &  0
\end{array} \!\!\right] & = r(B) + r(C) + r(E_BAF_C),
\label{217}
\\
r\!\left[\!\! \begin{array}{cc}  AA^*  & B  \\ B^* &  0
\end{array} \!\!\right] & = r[\, A, \ B\,] + r(B),
\label{217a}
\\
r\!\left[\!\! \begin{array}{cc} A & B \\ C & D
\end{array} \!\!\right] & = r(A) + r\!\left[\!\! \begin{array}{cc} 0 & E_AB \\
CF_A &  D - CA^{\dag}B\end{array} \!\!\right]\!.
\label{218}
\end{align}
If $\R(B) \subseteq \R(A)$  and  $\R(C^*) \subseteq \R(A^*),$ then
\begin{align}
r\!\left[\!\! \begin{array}{cc} A & B \\ C & D
\end{array} \!\!\right] & = r(A) + r(\, D - CA^{\dag}B\,).
\label{218a}
\end{align}
\end{lemma}

\begin{lemma} [\cite{T-laa10}] \label{T23}
Let $A \in {\mathbb C}_{{\rm H}}^{m},$ $B \in \mathbb C^{m\times n},$
 $D \in  {\mathbb C}_{{\rm H}}^{n},$  and define
$$
M_1 = \left[\!\!\begin{array}{cc}  A  & B  \\ B^{*}  & 0 \end{array}
\!\!\right],   \ \ \  M_2 = \left[\!\!\begin{array}{cc}  A  & B  \\ B^{*}  & D \end{array}
\!\!\right]\!.
$$
 Then$,$ the partial inertias of $M_1$ and $M_2$ can be expanded as
\begin{align}
& i_{\pm}(M_1) = r(B) + i_{\pm}(E_BAE_B),
\label{21x}
\\
& i_{\pm}(M_2)  =i_{\pm}(A) + i_{\pm}\!\left[\!\!\begin{array}{cc} 0 &  E_AB
 \\
 B^{*}E_A & D - B^{*}A^{\dag}B \end{array}\!\!\right]\!.
\label{24pp}
\end{align}
In particular$,$
\begin{enumerate}
\item[{\rm(a)}] If $A \succcurlyeq 0,$ then
\begin{align}
i_{+}(M_1)= r[\, A,  \, B \,], \ \ i_{-}(M_1) = r(B).
\label{111j}
\end{align}

\item[{\rm(b)}] If $A \preccurlyeq 0,$ then
\begin{align}
i_{+}(M_1) = r(B), \ \ i_{-}(M_1) = r[\, A,  \, B \,].
\label{112j}
\end{align}

\item[{\rm(c)}] If ${\mathscr R}(B) \subseteq {\mathscr R}(A),$ then
\begin{align}
i_{\pm}(M_2) = i_{\pm}(A) + i_{\pm}(\,  D - B^*A^{\dag}B \,).
\label{113j}
\end{align}

\item[{\rm (d)}] $i_{\pm}(M_1) = m \Leftrightarrow  i_{\mp}(E_BAE_B) =0 \  and  \ r(E_BAE_B) =r(E_B).$

\item[{\rm(e)}] $M_2 \succcurlyeq 0$ $\Leftrightarrow$ $A \succcurlyeq 0,$ ${\mathscr R}(B) \subseteq
{\mathscr R}(A)$ and $D - B^{*}A^{\dag}B \succcurlyeq 0$ $\Leftrightarrow$
$D \succcurlyeq 0,$ ${\mathscr R}(B^{*}) \subseteq {\mathscr R}(D)$ and $ A - BC^{\dag}B^{*} \succcurlyeq 0.$

\item[{\rm(f)}] $M_2 \succ  0$ $\Leftrightarrow A \succ 0$ and $D - B^{*}A^{-1}B \succ 0$
$\Leftrightarrow$ $D \succ 0$  and $A - BD^{-1}B^{*} \succ 0.$

\item[{\rm(g)}] Under $A \succcurlyeq 0$ and $A_1 \succcurlyeq 0,$  the inequality $A \succcurlyeq A_1$ holds if and only if
${\mathscr R}(A_1) \subseteq {\mathscr R}(A)$ and $A_1 - A_1A^{\dag}A_1 \succcurlyeq 0.$
\end{enumerate}
\end{lemma}

\begin{lemma} \label{T21}
Let $A, \, B \in \mathbb C_{{\rm H}}^{m}$ and $P\in \mathbb
C^{m\times n}.$
\begin{enumerate}
\item[{\rm(a)}] If $A \succcurlyeq B,$ then $P^{*}AP \succcurlyeq P^{*}BP.$

\item[{\rm(b)}] $A \succcurlyeq 0$ if and only if $A^{\dag} \succcurlyeq 0.$

\item[{\rm(c)}] If $ I_m - A \succcurlyeq 0,$ then $ I_m - PP^{\dag}APP^{\dag} \succcurlyeq 0.$

\item[{\rm(d)}] If $ I_m - A \succ  0,$ then $ I_m - PP^{\dag}APP^{\dag} \succ  0.$
\end{enumerate}
\end{lemma}

\noindent {\bf Proof.} \ Result (a) is obvious from the definition of the nonnegative
definiteness of Hermitian matrix. Result (b) is obvious from similarity decomposition of
$A$ and the definition of the Moore--Penrose inverse of a matrix. If $A$ is Hermitian, then we
can find  by Lemma \ref{T23}(c), $*$-congruence transformation and (\ref{218a}) that
\begin{align*}
i_{\pm}(\, I_m - PP^{\dag}APP^{\dag} \,)  & =
i_{\pm}\!\left[\!\!\begin{array}{cc} A & APP^{\dag} \\ PP^{\dag}A & I_m \end{array} \!\!
\right] - i_{\pm}(A) = i_{\pm}\!\left[\!\!\begin{array}{cc} A  - APP^{\dag}A & 0 \\ 0 & I_m \end{array} \!\!
\right] - i_{\pm}(A)
\\
& = i_{\pm}(I_m) + i_{\pm}(\, A  - APP^{\dag}A \,) - i_{\pm}(A)
\\
& = i_{\pm}(I_m) + i_{\pm}\!\left[\!\!\begin{array}{cc} P^*P & P^*A \\ AP & A \end{array} \!\!\right] -i_{\pm}(PP^*)
 - i_{\pm}(A)
\\
& = i_{\pm}(I_m)+ i_{\pm}\!\left[\!\!\begin{array}{cc} P^*P  - P^*AP & 0 \\ 0 & A \end{array} \!\!\right]
- i_{\pm}(PP^*) - i_{\pm}(A)
\\
& = i_{\pm}(I_m) + i_{\pm}[\, P^*(\, I_m - A\,)P\,] - i_{\pm}(PP^*),
\end{align*}
namely
\begin{align}
i_{+}(\, I_m - PP^{\dag}APP^{\dag} \,) & = m - r(P) + i_{+}[\, P^*(\, I_m - A\,)P\,],
\label{21pp}
\\
i_{-}(\, I_m - PP^{\dag}APP^{\dag} \,) & =  i_{-}[\, P^*(\, I_m - A\,)P\,].
\label{22pp}
\end{align}
If $A \preccurlyeq I_m,$ then (\ref{22pp}) reduces to
$$
i_{-}(\, I_m - PP^{\dag}APP^{\dag} \,) =  i_{-}[\, P^*(\, I_m - A\,)P\,] =0.
$$
Hence, (c) follows by Lemma \ref{T20}(d). If $A \prec  I_m,$ then $P^*(\, I_m - A\,)P \succcurlyeq 0$
and $i_{+}[\, P^*(\, I_m - A\,)P\,] = r[\, P^*(\, I_m - A\,)P\,] = r(P)$. Thus,
 (\ref{21pp}) reduces to
$$
 i_{+}(\, I_m - PP^{\dag}APP^{\dag} \,)  = m - r(P) + r(P) =m.
$$
Hence, (d) follows by Lemma \ref{T20}(c).  \qquad $\Box$

\medskip

The following results on the extremal ranks of $A - BXC$ with respect to a variable
matrix $X$ were shown in \cite{DG} by using  restricted singular value decompositions (RSVDs) of matrices,
and in \cite{Ti-se,TC} by using generalized inverses of matrices.

\begin{lemma} \label{T12}
 Let $A \in \mathbb C^{m \times n}, \, B \in \mathbb C^{m \times k}$ \
 and $C \in \mathbb C^{l \times n}$ be given$.$ Then the global maximum and minimum ranks
of $A - BXC$ with respect to $X \in \mathbb C^{ k\times l}$ are given by
\begin{align}\label{1.1}
\max_{X \in \mathbb C^{ k\times l}} \!\!r(\, A - BXC \,) & =  \min  \left\{ r[ \, A, \,  B \, ], \ \
 r \!\left[\!\! \begin{array}{c}  A  \\ C  \end{array}\!\right] \right\}\!,
\\
\label{1.2}
\min_{X \in \mathbb C^{ k\times l}} \!\!r(\, A - BXC \,)  & =  r[\, A, \, B \,] +
 r \!\left[\!\! \begin{array}{c} A \\ C \end{array} \!\!\right] -
 r \!\left[\!\! \begin{array}{cc} A & B \\ C & 0 \end{array} \!\!\right]\!.
\end{align}
In particular$,$
\begin{align}\label{12a}
\max_{X \in \mathbb C^{ k\times n}} \!\!r(\, A - BX \,) & =  \min
\left\{ r[ \, A, \,  B \, ], \ \ n \right\}\!,
\\
\label{12b}
 \min_{X \in \mathbb C^{ k\times n}} \!\!r(\, A - BX \,)  &
= r[\, A, \, B \,] - r(B).
\end{align}
\end{lemma}

\begin{lemma}[\cite{LT-jota,T-laa11}]\label{T212xx}
Let $A \in {\mathbb C}_{{\rm H}}^{m}$, $B \in {\mathbb C}^{m\times n}$ and $C \in {\mathbb C}^{p \times m}$ be
 given$.$ Then$,$ the extremal ranks and inertias of $A-BXC-(BXC)^{*}$ are given by
\begin{align}
& \max_{X\in {\mathbb C}^{n \times p}}\!\!r[\,A-BXC-(BXC)^{*}\,] = \min \left\{
r[\,A,\, B,\, C^{*}\,],
 \ \ r\!\left[\begin{array}{cc} A & B
\\B^{*} & 0
\end{array}\right]\!, \ \ r\!\left[\!\begin{array}{cc} A & C^{*}
\\ C & 0
\end{array}\!\right] \right\}\!,
\label{129}
\\
& \min_{X\in {\mathbb C}^{n \times p}}\!\!r[\,A-BXC-(BXC)^{*}\,]   = 2r[\,A,\,
B,\, C^{*}\,] + \max\{\, s_{+} + s_{-},  \ t_{+} + t_{-}, \ s_{+} +
t_{-}, \  s_{-} + t_{+} \, \},
\label{130}
\\
& \max_{X\in {\mathbb C}^{n \times p}}\!\!i_{\pm}[\,A-BXC-(BXC)^{*}\,]  =
\min\!\left\{i_{\pm}\!\left[\!\begin{array}{ccc} A
 & B  \\  B^{*}  & 0
   \end{array}\!\right], \ \ i_{\pm}\!\left[\!\begin{array}{ccc} A
 & C^{*}  \\  C  & 0
   \end{array}\!\right] \right\}\!,
\label{131}
\\
& \min_{X\in {\mathbb C}^{n \times p}}\!\!i_{\pm}[\,A-BXC-(BXC)^{*}\,] = r[\,A,\,
B,\, C^{*}\,] + \max\{\, s_{\pm}, \ \ t_{\pm} \, \},
\label{132}
\end{align}
where
\begin{align*}
s_{\pm}  = i_{\pm}\!\left[\!\!\begin{array}{cc} A & B \\ B^{*} & 0\end{array}\!\!\right]
 - r\!\left[\begin{array}{ccc} A & B  & C^{*} \\ B^{*} & 0 & 0
\end{array}\!\!\right]\!, \ \
 t_{\pm}  =i_{\pm}\!\left[\!\!\begin{array}{cc} A & C^{*} \\ C & 0\end{array}\!\!\right]
 - r\!\left[\begin{array}{ccc} A & B  & C^{*} \\ C & 0 & 0
\end{array}\!\!\right]\!.
\end{align*}
In particular,
\begin{align}
\max_{X\in {\mathbb C}^{n \times m}}\!\!r[\,A - BX - (BX)^{*}\,] & = \min \left\{m,
 \ \ r\!\left[\begin{array}{cc} A & B
\\B^{*} & 0
\end{array}\right] \right\}\!,
\label{137}
\\
\min_{X\in {\mathbb C}^{n \times m}}\!\!r[\,A - BX - (BX)^{*}\,] & = r\!\left[\begin{array}{cc} A & B
\\B^{*} & 0
\end{array}\right]-2r(B),
\label{138}
\\
\max_{X\in {\mathbb C}^{n \times m}}\!\!i_{\pm}[\,A  - BX - (BX)^{*}\,] &
= i_{\pm}\!\left[\!\begin{array}{ccc} A & B  \\  B^{*}  & 0
   \end{array}\!\right]\!,
\label{139}
\\
\min_{X\in {\mathbb C}^{n \times m}}\!\!i_{\pm}[\,A -  BX - (BX)^{*}\,]&
= i_{\pm}\!\left[\!\begin{array}{ccc} A
 & B  \\  B^{*}  & 0
   \end{array}\!\right] - r(B).
\label{140}
\end{align}
\end{lemma}

The matrices $X$ that satisfy (\ref{129})--(\ref{140}) (namely, the
global maximizers and minimizers of the objective rank and inertia
functions) are not necessarily unique and their expressions were
also given in \cite{LT-jota,T-laa11} by using certain simultaneous decomposition
of the three given matrices and their generalized inverses.

We also need the  following results on the ranks and inertias of the quadratic matrix-valued
functions
$$
A \pm (\, BX + C \,)(\, BX + C \,)^* = A \pm (BXX^*B^*  + BXC^* + CX^*B^* + CC^*)
$$
and their consequences.

\begin{lemma} [\cite{T-laa12}] \label{T212F}
Let $A \in \mathbb C_{{\rm H}}^{m}$ and   $B\in {\mathbb C}^{m \times k}$ and
$C \in {\mathbb C}^{m \times n}$ be given$,$ and let
$$
G_1 = \!\left[\!\! \begin{array}{cc}  A + CC^* & B \\ B^* & 0 \end{array} \!\!\right]\!, \ \
G_2 = \!\left[\!\! \begin{array}{cc}  A - CC^* & B \\ B^* & 0 \end{array} \!\!\right]\!, \ \
G_3 = \left[\!\! \begin{array}{ccc}  A & B & C \\ B^* & 0 & 0 \end{array} \!\!\right]\!.
$$
Then$,$ the following hold$.$
\begin{enumerate}
\item[{\rm(a)}]  The extremal ranks and inertias of
$\phi_1(X) =A + (\, BX + C \,)(\, BX + C \,)^*$  are given by
\begin{align}
& \max_{X\in {\mathbb C}^{k \times n}}\!\!r[\,\phi_1(X)\,] =
\min\left\{r[\,  A, \,  B, \, C \,], \ \  r(G_1), \ \ r(A) + n \right\},
\label{243j}
\\
&\min_{X\in {\mathbb C}^{k \times n}}\!\!r[\,\phi_1(X)\,] = 2r[\,A,
\,B, \, C\,]  + \max\{ \, h_1, \ \ h_2, \ \ h_3, \ \ h_4 \, \},
\label{244j}
\\
& \max_{X\in {\mathbb C}^{k \times n}}\!\!i_{+}[\,\phi_1(X)\,]
 = \min\left\{\, i_{+}(G_1), \ \   i_{+}(A)  +n  \right\},
\label{245j}
\\
& \max_{X\in {\mathbb C}^{k \times n}}\!\!i_{-}[\,\phi_1(X)\,]
 = \min\left\{\, i_{-}(G_1), \ \   i_{-}(A)  \right\},
\label{246j}
\\
&  \min_{X\in {\mathbb C}^{k \times n}} \!\!i_{+} [\,\phi_1(X)\,] =
r[\,A, \,B, \, C\,]  + \max \left\{ i_{+}(G_1) - r(G_3),  \  \
i_{+}(A)  - r[\,  A, \, B \,] \right\},
\label{247j}
\\
&  \min_{X\in {\mathbb C}^{k \times n}} \!\!i_{-} [\,\phi_1(X)\,] =
r[\,A, \,B, \, C\,]  + \max \left\{ i_{-}(G_1) - r(G_3),  \ \
i_{-}(A)  - r[\,  A, \, B \,] - n \right\},
\label{248j}
\end{align}
where
\begin{align*}
 h_1 & = r(G_1) - 2r(G_3), \ \ \ h_2  =  r(A) -  2r[\,A, \,  B\,] - n,
\\
h_3 & =  i_{-}(G_1) - r(G_3) + i_{+}(A) - r[\,A, \,  B\,],
\\
h_4 &=  i_{+}(G_1) - r(G_3) + i_{-}(A) - r[\,A, \,  B\,] -n.
\end{align*}

\item[{\rm(b)}]  The extremal ranks and inertias of
$\phi_2(X) = A - (\, BX + C \,)(\, BX + C \,)^*$  are given by
\begin{align}
&\max_{X\in {\mathbb C}^{k \times n}}\!\!r[\,\phi_2(X)\,] =
\min\left\{r[\,  A, \,  B, \, C \,], \ \  r(G_2), \ \ r(A) + n
\right\}, \label{249j}
\\
&\min_{X\in {\mathbb C}^{k \times n}}\!\!r[\,\phi_2(X)\,] = 2r[\,A,
\,B, \, C\,]  + \max\{ \, h_5, \ \ h_6, \ \ h_7, \ \ h_8 \, \},
\label{250j}
\\
& \max_{X\in {\mathbb C}^{k \times n}}\!\!i_{+}[\,\phi_2(X)\,]
 = \min\left\{\, i_{+}(G_2), \ \   i_{+}(A) \, \right\},
\label{251j}
\\
& \max_{X\in {\mathbb C}^{k \times n}}\!\!i_{-}[\,\phi_2(X)\,]
 = \min\left\{\, i_{-}(G_2), \ \   i_{-}(A)  +n \, \right\},
\label{252j}
\\
&  \min_{X\in {\mathbb C}^{k \times n}} \!\!i_{+} [\,\phi_2(X)\,]  =
r[\,A, \,B, \, C\,]  + \max \left\{ i_{+}(G_2) - r(G_3),  \ \
i_{+}(A)  - r[\,  A, \, B \,] - n \right\}, \label{253j}
\\
&  \min_{X\in {\mathbb C}^{k \times n}} \!\!i_{-}[\,\phi_2(X)\,]  =
r[\,A, \,B, \, C\,]  + \max \left\{ i_{-}(G_2) - r(G_3),
 \ \ i_{-}(A)  - r[\,  A, \, B \,] \right\},
\label{254j}
\end{align}
where
\begin{align*}
 h_5 & = r(G_2) - 2r(G_3), \ \ \  h_6  =  r(A) -  2r[\,A, \,  B\,] - n,
\\
h_7 & =  i_{+}(G_2) - r(G_3) + i_{-}(A) - r[\,A, \,  B\,],
\\
h_8 &=  i_{-}(G_2) - r(G_3) + i_{+}(A) - r[\,A, \,  B\,] -n.
\end{align*}
\end{enumerate}
\end{lemma}

When $C =0$,  Lemma \ref{T212F} reduces to the following result.

\begin{corollary} [\cite{T-mcm}] \label{T213F}
Let  $ A \in \mathbb C_{{\rm H}}^{m}$ and  $B \in {\mathbb C}^{m \times n}$ be given$,$  and let
$M = \left[\!\! \begin{array}{cccc}  A & B
 \\ B^* & 0 \end{array} \!\!\right]\!.$ Then$,$ the following hold$.$
\begin{enumerate}
\item[{\rm(a)}]  The extremal ranks and partial inertias of
$A \pm  BXX^*B^*$  are given by
\begin{align}
& \max_{X \in {\mathbb C}^{n \times n}} \!\!r(\, A + BXX^*B^* \,)  = r[\, A, \, B\,],
\label{278}
\\
& \min_{X \in {\mathbb C}^{n \times n}}\!\!r(\, A + BXX^*B^* \,)   = i_{+}(A) + r[\, A, \, B\,] - i_{+}(M),
\label{279}
\\
& \max_{X \in {\mathbb C}^{n \times n}}\!\!i_{+}(\, A + BXX^*B^* \,)   = i_{+}(M),
\label{280}
\\
& \max_{X \in {\mathbb C}^{n \times n}}\!\!i_{-}(\, A + BXX^*B^* \,)   = i_{-}(A),
\label{282}
\\
& \min_{X\in {\mathbb C}^{n}_{{\rm H}}}\!\!i_{+}(\, A + BXX^*B^* \,)  = i_{+}(A),
\label{281}
\\
& \min_{X \in {\mathbb C}^{n \times n}}\!\!i_{-}(\, A + BXX^*B^* \,)   =r[\, A, \, B\,] - i_{+}(M),
 \label{283}
\end{align}
and
\begin{align}
& \max_{X \in {\mathbb C}^{n \times n}}\!\!r(\, A - BXX^*B^* \,)  = r[\, A, \, B\,],
\label{284}
\\
&\min_{X \in {\mathbb C}^{n \times n}}\!\!r(\, A - BXX^*B^* \,)  =  i_{-}(A) + r[\, A, \, B\,] - i_{-}(M),
\label{285}
\\
& \max_{X \in {\mathbb C}^{n \times n}}\!\!i_{+}(\, A - BXX^*B^* \,)   = i_{+}(A),
\label{286}
\\
& \max_{X \in {\mathbb C}^{n \times n}}\!\!i_{-}(\, A - BXX^*B^* \,)   = i_{-}(M),
\label{288}
\\
& \min_{X \in {\mathbb C}^{n \times n}}\!\!i_{+}(\, A - BXX^*B^* \,)  =  r[\,A, \, B\,] -  i_{-}(M),
\label{287}
\\
&  \min_{X \in {\mathbb C}^{n \times n}}\!\!i_{-}(\, A - BXX^*B^* \,)  = i_{-}(A).
\label{289}
\end{align}

\item[{\rm(b)}] If $A \succcurlyeq 0,$ then
\begin{align}
& \max_{X \in {\mathbb C}^{n \times n}} \!\!r(\, A + BXX^*B^* \,)  = r[\, A, \, B\,],
\label{290}
\\
& \min_{X \in {\mathbb C}^{n \times n}}\!\!r(\, A + BXX^*B^* \,)  = r(A),
\label{291}
\end{align}
and
\begin{align}
& \max_{X \in {\mathbb C}^{n \times n}}\!\!r(\, A - BXX^*B^* \,)  = r[\, A, \, B\,],
\label{292}
\\
&\min_{X \in {\mathbb C}^{n \times n}}\!\!r(\, A - BXX^*B^* \,)  =  r[\, A, \, B\,] - r(B),
\label{293}
\\
& \max_{X \in {\mathbb C}^{n \times n}}\!\!i_{+}(\, A - BXX^*B^* \,)   = r(A),
\label{294}
\\
& \max_{X \in {\mathbb C}^{n \times n}}\!\!i_{-}(\, A - BXX^*B^* \,)   = r(B),
\label{296}
\\
& \min_{X \in {\mathbb C}^{n \times n}}\!\!i_{+}(\, A - BXX^*B^* \,)  =  r[\,A, \, B\,] -  r(B),
\label{295}
\\
&  \min_{X \in {\mathbb C}^{n \times n}}\!\!i_{-}(\, A - BXX^*B^* \,)  = 0.
\label{297}
\end{align}

\end{enumerate}
\end{corollary}

\section{General solutions $AXB \succcurlyeq \,(\succ, \, \preccurlyeq, \, \prec ) \, C$ and their properties}
\renewcommand{\theequation}{\thesection.\arabic{equation}}
\setcounter{section}{3}
\setcounter{equation}{0}

A necessary condition for (\ref{10b}) to hold is $AXB = (AXB)^*$. In such a case, the
matrix $X$ satisfying  $AXB = (AXB)^*$ is called a symmetrizer of
$AXB$; see \cite{BK}. In this section, we derive an analytical presentation for the general solution of
the LMI in (\ref{10b}) by using the given matrices and their generalized inverses,
and establish various algebraic properties of the LMI.

\begin{theorem} \label{T51}
Let $A \in \mathbb C^{m\times p},$ $B \in \mathbb C^{q\times m}$ and
$C \in \mathbb C_{{\rm H}}^{m}$ be given$,$ and define
$M = [\, E_A, \, F_B \,].$ Then$,$ the following hold$.$
\begin{enumerate}
\item[{\rm(a)}] There exists an $X \in \mathbb C^{p\times q}$ such that
 \begin{equation}
AXB \succcurlyeq C
\label{52w}
\end{equation}
if and only if
\begin{equation}
M^{*}CM \preccurlyeq 0 \ \ and \ \ \R(M^{*}CM) = \R(M^*C).
 \label{53w}
\end{equation}
In this case$,$ the general solution of {\rm (\ref{52w})} and the corresponding
$AXB$ can be written  in the following parametric forms
\begin{align}
& X = A^{\dag}CB^{\dag} - A^{\dag}CM(M^{*}CM)^{\dag}M^{*}CB^{\dag} +
A^{\dag}E_MUU^{*}E_MB^{\dag}+ W - A^{\dag}AWBB^{\dag},
 \label{54w}
\\
& AXB = C - CM(M^{*}CM)^{\dag}M^{*}C + E_MUU^{*}E_M,
 \label{55w}
\end{align}
where $U\in \mathbb C^{m\times m}$ and $W\in \mathbb C^{p\times q}$ are arbitrary$.$

\item[{\rm(b)}] There exists an $X \in \mathbb C^{p\times q}$ such that
\begin{equation}
AXB  \succ  C
\label{56w}
\end{equation}
if and only if
\begin{equation}
M^{*}CM \preccurlyeq 0 \ \ and \ \  r\!\left[\!\! \begin{array}{cccc} A  & 0 & C
\\ 0 & B^{*} & C \end{array} \!\!\right] = m +  r[\, A, \, B^{*} \,].
\label{57w}
\end{equation}
 In this case$,$ the general solution of {\rm (\ref{56w})} can be written as
{\rm (\ref{54w})}$,$ in which $U\in \mathbb C^{m\times m}$ is any
matrix such that $r[\, CM, \ E_MU\,] = m,$ and $W\in \mathbb C^{p\times q}$ is arbitrary$.$
\end{enumerate}
\end{theorem}

\noindent {\bf Proof.} \ Inequality (\ref{10b}) is obviously equivalent to the following linear-quadratic matrix equation
\begin{equation}
AXB  = C + YY^{*}.
\label{514w}
\end{equation}
By Lemma \ref{T27}, this equation is solvable for $X$ if and only if
\begin{equation}
 E_A(\, C + YY^{*}\,)= 0 \ \ {\rm and} \ \  (\, C + YY^{*}\,)F_B =0,
 \label{515w}
\end{equation}
that is,
\begin{equation}
\left[\!\! \begin{array}{c} E_A \\ F_B \end{array} \!\!\right]YY^{*}
= - \left[\!\! \begin{array}{c} E_AC \\ F_BC \end{array}
\!\!\right]\!.
\label{516w}
\end{equation}
By Lemma \ref{T25}(b), this quadratic matrix equation is solvable for $YY^{*}$ if and only if
$$
\left[\!\! \begin{array}{c} E_A \\ F_B \end{array} \!\!\right]\!C[\, E_A, \, F_B \,] \preccurlyeq 0  \ \ {\rm and}
 \ \ r\left( \left[\!\! \begin{array}{c} E_A \\ F_B \end{array} \!\!\right]\!C[\, E_A, \, F_B \,] \right) =
r(C[\, E_A, \, F_B \,]),
$$
establishing (\ref{53w}). Under (\ref{53w}), the general solution of (\ref{514w}) can be written as
$$
YY^{*} = -CM(M^{*}CM)^{\dag}M^{*}C + E_MUU^{*}E_M,
$$
where $U\in \mathbb C^{m\times m}$ is arbitrary. Substituting the $YY^{*}$ into
(\ref{514w}) gives
\begin{equation}
AXB = C - CM(M^{*}CM)^{\dag}M^{*}C + E_MUU^{*}E_M.
 \label{517w}
\end{equation}
By Lemma \ref{T27}, the general solution of (\ref{517w}) is
$$
X  = A^{\dag}CB^{\dag} - A^{\dag}CM(M^{*}CM)^{\dag}M^{*}CB^{\dag} +
A^{\dag}E_MUU^{*}E_MB^{\dag} + W - A^{\dag}AWBB^{\dag},
$$
establishing (\ref{54w}) and (\ref{55w}).

It can be seen from (\ref{517w}) that (\ref{56w}) holds if and only if
\begin{equation}
- CM(M^{*}CM)^{\dag}M^{*}C + E_MUU^{*}E_M \succ 0
 \label{518w}
\end{equation}
for some $U$. Under (\ref{53w}), we have
$$
r[\, -CM(M^{*}CM)^{\dag}M^{*}C + E_MUU^{*}E_M \,]= r[\, -CM(M^{*}CM)^{\dag}M^{*}C, \,  E_MUU^{*}E_M \,]
= r[\, CM, \,  E_MU \,].
$$
Also by (\ref{216}),
\begin{align}
& r(CM)  = r(C[\, E_A, \, F_B \,]) = r\!\left[\!\! \begin{array}{cc} C  & C \\ A^* & 0 \\ 0 & B
\end{array} \!\!\right] - r(A) - r(B) = r\!\left[\!\! \begin{array}{cccc} A  & 0 & C
\\ 0 & B^{*} & C \end{array} \!\!\right] - r(A) - r(B),
\label{518w1}
\\
& r(M)  = r[\, E_A, \, F_B \,]   = r\!\left[\!\! \begin{array}{cccc} A  & 0 & I_m
\\ 0 & B^{*} & I_m \end{array} \!\!\right] - r(A) - r(B) = m +  r[\, A, \, B^* \,]  - r(A) - r(B).
\label{518w2}
\end{align}
Hence,
\begin{align}
& \max_{U} r[\, -CM(M^{*}CM)^{\dag}M^{*}C + E_MUU^{*}E_M  \,]  = r[\, CM, \, E_M  \,] \nb
\\
& = r (MM^{\dag}CM)  + r(E_M) = r(CM) + m - r(M) = r\!\left[\!\! \begin{array}{cccc} A  & 0 & C
\\ 0 & B^{*} & C
\end{array} \!\!\right] - r[\, A, \, B^{*} \,].
\label{518w3}
\end{align}
Thus, (\ref{518w}) is equivalent to (\ref{57w}).
\qquad $\Box$

\medskip

The following  result can  be shown similarly.

\begin{corollary} \label{T51a}
Let $A \in \mathbb C^{m\times p},$ $B \in \mathbb C^{q\times m}$ and
$C \in \mathbb C_{{\rm H}}^{m}$ be given$,$ and define
$M = [\, E_A, \, F_B \,].$ Then$,$ the following hold$.$
\begin{enumerate}
\item[{\rm(a)}] There exists an $X \in \mathbb C^{p\times q}$ such that
\begin{equation}
AXB  \preccurlyeq C
 \label{58w}
\end{equation}
if and only if
\begin{equation}
M^{*}CM \succcurlyeq 0 \ \ and \ \ \R(M^{*}CM) = \R(M^*C).
 \label{59w}
\end{equation}
 In this case$,$ the general solution of
 {\rm (\ref{58w})} and the corresponding $AXB$ can be written in the following parametric forms
\begin{align}
& X  = A^{\dag}CB^{\dag} - A^{\dag}CM(M^{*}CM)^{\dag}M^{*}CB^{\dag} -
A^{\dag}E_MUU^{*}E_MB^{\dag} + W - A^{\dag}AWBB^{\dag},
 \label{510w}
\\
& AXB = C - CM(M^{*}CM)^{\dag}M^{*}C - E_MUU^{*}E_M,
 \label{511w}
\end{align}
where $U\in \mathbb C^{m\times m}$ and $W\in \mathbb C^{p\times q}$ are arbitrary$.$

\item[{\rm(b)}] There exists an $X \in \mathbb C^{p\times q}$ such that
\begin{equation}
AXB \prec  C
\label{512w}
\end{equation}
if and only if
\begin{equation}
M^{*}CM \succcurlyeq 0 \ \ and \ \  r\!\left[\!\! \begin{array}{cccc}  A  & 0 & C
\\ 0 & B^{*} & C \end{array} \!\!\right] = m +  r[\, A, \, B^{*} \,].
\label{513w}
\end{equation}
In this case$,$ the general
solution of {\rm (\ref{512w})} can be written as {\rm (\ref{510w})}$,$
in which $U$ is any matrix such that $r[\, CM, \ E_MU\,] = m,$ and
$W\in \mathbb C^{p\times q}$ is
arbitrary$.$
\end{enumerate}
\end{corollary}

We next establish some algebraic properties of the fixed parts in (\ref{54w}) and (\ref{510w}).

\begin{corollary} \label{T52}
Let $A \in \mathbb C^{m\times p},$ $B \in \mathbb C^{q\times m}$ and
$C \in \mathbb C_{{\rm H}}^{m}$ be given$,$
 and define
\begin{align}
& \widehat{X} = A^{\dag}CB^{\dag} - A^{\dag}CM(M^{*}CM)^{\dag}M^{*}CB^{\dag}, \ \
  M = [\, E_A, \, F_B \,],  \ \  N = \!\left[\!\! \begin{array}{cccc} C & C  & A & 0 \\
  C & C  & 0 & B^{*} \\ A^{*} & 0 & 0 & 0 \\  0 & B  & 0 & 0  \end{array} \!\!\right]\!.
\label{520w}
\end{align}
\begin{enumerate}
\item[{\rm(a)}] Under the condition that {\rm (\ref{52w})} has a solution$,$
the $\widehat{X}$ in {\rm (\ref{520w})} satisfies $A\widehat{X}B \succcurlyeq C,$ and
\begin{align}
& i_{\pm}(A\widehat{X}B) = r(A) + r(B) + i_{\pm}(C) -  i_{\pm}(N),
\label{521w}
\\
& r(\widehat{X}) =  r(A\widehat{X}B) =  2r(A) + 2r(B) + r(C) -  r(N),
\label{522w}
\\
&i_{+}(\, A\widehat{X}B - C\,) = r(\, A\widehat{X}B - C\,) = 2r\!\left[\!\! \begin{array}{cccc} A & 0 & C  \\
 0 & B^{*}  & C \end{array} \!\!\right] - r(N).
\label{523w}
\end{align}

\item[{\rm(b)}] Under the condition that {\rm (\ref{58w})} has a solution$,$
the $\widehat{X}$ in {\rm (\ref{520w})} satisfies $A\widehat{X}B \preccurlyeq C,$ and
\begin{align}
& i_{\pm}(A\widehat{X}B) = r(A) + r(B) + i_{\pm}(C) -  i_{\pm}(N),
\label{524w}
\\
& r(\widehat{X}) =  r(A\widehat{X}B) =  2r(A) + 2r(B) + r(C) -  r(N),
\label{525w}
\\
&i_{-}(\, A\widehat{X}B - C\,) =  r(\,A\widehat{X}B -C\,)  = 2r\!\left[\!\! \begin{array}{cccc} A & 0 & C  \\
 0 & B^{*}  & C \end{array} \!\!\right] - r(N).
\label{526w}
\end{align}
\end{enumerate}
\end{corollary}

\noindent {\bf Proof.} \ Under the condition that (\ref{52w}) has a solution$,$
set $U = W = 0$ in (\ref{54w}). Then we see that the $\widehat{X}$ in (\ref{520w})
is a solution of $AXB \succcurlyeq C.$ Also note from (\ref{517w}) that
\begin{align}
A\widehat{X}B =  C - CM(M^{*}CM)^{\dag}M^{*}C, \ \ \  C - A\widehat{X}B =  CM(M^{*}CM)^{\dag}M^{*}C \preccurlyeq 0.
\label{527w}
\end{align}
In this case, applying (\ref{113j}) and (\ref{218a}) to (\ref{527w}), we obtain
\begin{align}
& r(\widehat{X}) = r[\,  A^{\dag}CB^{\dag} - A^{\dag}CM(M^{*}CM)^{\dag}M^{*}CB^{\dag} \,]
= r\!\left[\!\begin{array}{cc} M^{*}CM &  M^{*}CB^{\dag} \\
 A^{\dag}CM   & A^{\dag}CB^{\dag} \end{array}\!\right] - r(M^{*}CM),
\label{528w}
\\
& i_{\pm}(A\widehat{X}B) = i_{\pm}[\, C - CM(M^{*}CM)^{\dag}M^{*}C \,]
= i_{\pm}\!\left[\!\begin{array}{cc} M^{*}CM &  M^{*}C \\
 CM   & C \end{array}\!\right] - i_{\pm}(M^{*}CM),
\label{529w}
\\
&  r(\, C - A\widehat{X}B \,) = r[\, CM(M^{*}CM)^{\dag}M^{*}C \,]
= r\!\left[\!\begin{array}{cc} M^{*}CM &  M^{*}C \\
 CM   & 0 \end{array}\!\right] - r(M^{*}CM).
\label{530w}
\end{align}
Applying elementary matrix operations, congruence matrix operations and (\ref{217}), we obtain
\begin{align}
r\!\left[\!\begin{array}{cc} M^{*}CM &  M^{*}CB^{\dag}
\\
 A^{\dag}CM   & A^{\dag}CB^{\dag} \end{array}\!\right] & = r\!\left(\left[\!\begin{array}{cc} E_A \\ F_B \\
 A^{\dag} \end{array}\!\right]C [\, E_A, \,F_B, \,  B^{\dag} \,] \right)
  = r\!\left(\left[\!\begin{array}{cc} I_m
  \\
  F_B
  \\
 A^{\dag} \end{array}\!\right]C [\, E_A, \,I_m, \,  B^{\dag} \,] \right)  = r(C),
\label{531w}
\\
 i_{\pm}\!\left[\!\begin{array}{cc} M^{*}CM &  M^{*}C
 \\
 CM   & C \end{array}\!\right] & = i_{\pm}\!\left[\!\begin{array}{cc} 0 &  0 \\
 0  & C \end{array}\!\right] = i_{\pm}(C),
 \label{532w}
 \\
 r\!\left[\!\begin{array}{cc} M^{*}CM &  M^{*}C \\ CM   & 0 \end{array}\!\right]
& = 2r(M^*C) = 2r\!\left[\!\! \begin{array}{cccc} A & 0 & C  \\
 0 & B^{*}  & C \end{array} \!\!\right]- 2r(A) - 2r(B),
\label{533w}
\\
i_{\pm}(M^{*}CM) & = i_{\pm}\!\left(\left[\!\begin{array}{cc} E_A \\ F_B \end{array}\!\right]C [\, E_A, \,F_B \,]
\right) =  i_{\pm}(N) -  r(A) - r(B),
\label{534w}
 \\
 r(M^{*}CM) & = r(N) -  2r(A) - 2r(B).
\label{535w}
\end{align}
Substituting these formulas into  (\ref{528w})--(\ref{530w})  yields (\ref{521w})--(\ref{523w}).
Results (b) can be shown similarly. \qquad $\Box$

\begin{corollary} \label{T53}
Let $A \in \mathbb C^{m\times p},$ $B \in \mathbb C^{q\times m}$ and
$C \in \mathbb C_{{\rm H}}^{m}$ be given$,$  and let $M$ and $N$ be of the forms in
 {\rm (\ref{520w})}$.$ Also assume that {\rm (\ref{52w})} is feasible$,$ and define
\begin{equation}
{\cal S}_1 =\{\, X\in {\mathbb C}^{p\times q} \ | \  AXB \succcurlyeq  C \, \}.
\label{536w}
\end{equation}
Then$,$ the following hold$.$
\begin{enumerate}
\item[{\rm(a)}]  The minimal matrices of $AXB$ and $AXB - C$ subject
  to $X \in {\cal S}_1$ in the L\"owner partial ordering are given by
\begin{align}
&\min_{\succcurlyeq}\{ AXB \ | \  X \in {\cal S}_1 \}  = C - CM(M^{*}CM)^{\dag}M^{*}C,
\label{537w}
\\
& \min_{\succcurlyeq}\{ AXB - C \ | \  X \in {\cal S}_1 \}  =  - CM(M^{*}CM)^{\dag}M^{*}C.
\label{538w}
\end{align}

\item[{\rm(b)}]  The  extremal ranks and partial inertias of $AXB$ and $AXB - C$ subject
  to $X \in {\cal S}_1$ are given by
\begin{align}
&\max_{X \in {\cal S}_1} r(AXB)  =  \max_{X \in {\cal S}_1} i_{+}(AXB)  =r(A) + r(B) - r[\, A, \, B^*\,],
\label{539w}
\\
&\min_{X \in {\cal S}_1} r(AXB)   =  \min_{X \in {\cal S}_1} i_{+}(AXB)  = r(A) + r(B) + i_{+}(C) -  i_{+}(N),
\label{540w}
\\
& \max_{X \in {\cal S}_1} i_{-}(AXB)  =  r(A) + r(B) + i_{-}(C) -  i_{-}(N),
\label{541w}
\\
& \min_{X \in {\cal S}_1} i_{-}(AXB)  = 0,
\label{542w}
\\
& \max_{X \in {\cal S}_1} r(\, AXB - C \,)  =  r(N) - r(A) - r(B) - r[\, A, \, B^*\,],
\label{543w}
\\
& \min_{X \in {\cal S}_1} r(\, AXB - C \,)  = r\!\left[\!\! \begin{array}{cccc} A  & 0 & C
\\ 0 & B^{*} & C
\end{array} \!\!\right] - r(A) - r(B).
\label{544w}
\end{align}
\end{enumerate}
In consequence$,$ the following hold$.$
\begin{enumerate}
\item[{\rm(c)}] There exists an $X \in \mathbb C^{p\times q}$ such that
$ AXB \succ 0$  and $ AXB \succcurlyeq C$ if and only if $r[\, A, \, B^*\,]  = r(A) + r(B) - m.$

\item[{\rm(d)}] There exists an $X \in \mathbb C^{p\times q}$ such that
$ 0 \succ  AXB \succcurlyeq C$ if and only if $C\prec 0$ and $r\!\left[\!\! \begin{array}{cccc} A  & 0 & C
\\ 0 & B^{*} & C \end{array} \!\!\right] =  r(A) +  r(B).$

\item[{\rm(e)}] There exists an $X \in \mathbb C^{p\times q}$ such that
$ 0 \succcurlyeq AXB \succcurlyeq C$ if and only if $C\preccurlyeq 0.$

\item[{\rm(f)}] There always exists an $X \in \mathbb C^{p\times q}$ such that
$ AXB \succcurlyeq 0$ and $ AXB \succcurlyeq C.$
\end{enumerate}
\end{corollary}

\noindent {\bf Proof.} \  From (\ref{517w}), $AXB$ and $AXB - C$ subject to $X \in {\cal S}_1$
can be written as
\begin{align}
& AXB = C - CM(M^{*}CM)^{\dag}M^{*}C + E_MUU^{*}E_M =  A\widehat{X}B  + E_MUU^{*}E_M,
\label{554w}
\\
& AXB  -  C =  - CM(M^{*}CM)^{\dag}M^{*}C + E_MUU^{*}E_M = A\widehat{X}B - C + E_MUU^{*}E_M.
 \label{555w}
\end{align}
Hence,
\begin{align}
AXB \succcurlyeq C - CM(M^{*}CM)^{\dag}M^{*}C, \ \ \ AXB - C \succcurlyeq  - CM(M^{*}CM)^{\dag}M^{*}C
\label{556w}
\end{align}
hold for any $U\in {\mathbb C}^{m \times m}$, which  implies  (\ref{537w}) and  (\ref{538w}).

Applying elementary matrix operations, congruence matrix operations and (\ref{217}), we obtain
\begin{align}
 r(E_M) & = m - r[\,E_A, \,F_B \,] = m - r(E_A) - r(F_BA)  = r(A) + r(B) - r[\, A, \, B^*\,],
\label{557w}
\\
 r[\, E_M, \, A\widehat{X}B \,] & = r(E_M) + r(A\widehat{X}BM)
  = r(E_M) + r[\, A\widehat{X}BE_A, \, A\widehat{X}BF_B \,] = r(E_M) \nb
  \\
  & = r(A) + r(B) - r[\, A, \, B^*\,],
\label{558w}
\\
 i_{\pm}\!\left[\!\! \begin{array}{cc} A\widehat{X}B & E_M \\  E_M
 & 0 \end{array}\!\!\right] & =  i_{\pm}\!\left[\!\! \begin{array}{cc} 0 & E_M \\  E_M
 & 0 \end{array}
\!\!\right]  = r(E_M) = r(A) + r(B) - r[\, A, \, B^*\,].
\label{559w}
\end{align}
Applying (\ref{278})--(\ref{283}) to (\ref{554w}) and (\ref{555w}) and simplifying by
(\ref{557w})--(\ref{559w}), we obtain
 \begin{align*}
\max_{X \in {\cal S}_1} r(AXB) & = \max_{U\in {\mathbb C}^{m \times m}} r(\, A\widehat{X}B + E_MUU^{*}E_M \,)
 = r[\, E_M, \, A\widehat{X}B \,] = r(A) + r(B) - r[\, A, \, B^*\,],
\\
\min_{X \in {\cal S}_1} r(AXB) & = \min_{U\in {\mathbb C}^{m \times m}} r(\, A\widehat{X}B +  E_MUU^{*}E_M  \,) \nb
\\
&= i_{+}(A\widehat{X}B) +  r[\,  E_M, \, A\widehat{X}B \,] - i_{+}\!\left[\!\! \begin{array}{cc}
A\widehat{X}B &  E_M  \\  E_M & 0 \end{array}
\!\!\right]  = r(A) + r(B) + i_{+}(C) -  i_{+}(N),
\\
 \max_{X \in {\cal S}_1} i_{+}(AXB) & = \max_{U\in {\mathbb C}^{m \times m}} i_{+}(\, A\widehat{X}B +
  E_MUU^{*}E_M  \,)
 =  i_{+}\!\left[\!\! \begin{array}{cc} A\widehat{X}B &  E_M \\  E_M & 0 \end{array}
\!\!\right] =  r(A) + r(B) - r[\, A, \, B^*\,],
\\
 \min_{X \in {\cal S}_1} i_{+}(AXB) & = \min_{U\in {\mathbb C}^{m \times m}} i_{+}(\, A\widehat{X}B +
  E_MUU^{*}E_M  \,)  = i_{+}(A\widehat{X}B) = r(A) + r(B) + i_{+}(C) -  i_{+}(N),
\\
 \max_{X \in {\cal S}_1} i_{-}(AXB) & =
 \max_{U\in {\mathbb C}^{m \times m}} i_{-}(\, A\widehat{X}B +   E_MUU^{*}E_M  \,)
 =  i_{-}(A\widehat{X}B) =  r(A) + r(B) + i_{-}(C) -  i_{-}(N),
\\
 \min_{X \in {\cal S}_1} i_{-}(AXB) & =
 \min_{U\in {\mathbb C}^{m \times m}} i_{-}(\, A\widehat{X}B +  E_MUU^{*}E_M\,)
 = r[\, E_M, \, A\widehat{X}B\,] -  i_{+}\!\left[\!\! \begin{array}{cc} A\widehat{X}B & E_M \\  E_M
 & 0 \end{array}
\!\!\right] = 0,
\end{align*}
establishing (\ref{539w})--(\ref{544w}). Note from  (\ref{53w}), that
\begin{align*}
r(\, AXB - C \,) & = r[\, - CM(M^{*}CM)^{\dag}M^{*}C + E_MUU^{*}E_M  \,]
\\
&= r[\,  - CM(M^{*}CM)^{\dag}M^{*}C, \,  E_MUU^{*}E_M \,] = r[\, CM, \, E_MU \,].
\end{align*}
Hence, we can find from (\ref{215}), (\ref{216}), (\ref{533w}) and (\ref{535w}) that
\begin{align*}
\max_{X \in {\cal S}_1} r(\, AXB - C \,) & =
\max_{U\in {\mathbb C}^{m \times m}} r[\, CM, \, E_MU  \,]
\\
& = r[\,  CM, \, E_M \,] = r(M^*CM) + r(E_M) =  r(N) - r(A) - r(B) - r[\, A, \, B^*\,],
\\
\min_{X \in {\cal S}_1} r(\, AXB - C \,) & =
\min_{U\in {\mathbb C}^{m \times m}} r[\,  CM, \, E_MU  \,] = r(CM) = r\!\left[\!\! \begin{array}{cccc} A  & 0 & C
\\ 0 & B^{*} & C
\end{array} \!\!\right] - r(A) - r(B),
\end{align*}
establishing (\ref{543w}) and (\ref{544w}). Result (b) can be shown similarly. \qquad $\Box$

\begin{corollary} \label{T53a}
Let $A \in \mathbb C^{m\times p},$ $B \in \mathbb C^{q\times m}$ and
$C \in \mathbb C_{{\rm H}}^{m}$ be given$,$  and let $M$ and $N$ be of the forms in
 {\rm (\ref{520w})}$.$ Also assume that {\rm (\ref{58w})} is feasible$,$ and define
\begin{equation}
{\cal S}_2 =\{\, X\in {\mathbb C}^{p\times q} \ | \  AXB \preccurlyeq  C \, \}.
\label{545w}
\end{equation}
Then$,$ the following hold$.$
\begin{enumerate}
\item[{\rm(a)}]  The maximal matrices of $AXB$ and $AXB - C$ subject to $X \in {\cal S}_2$
in the L\"owner partial ordering are given by
\begin{align}
&\max_{\succcurlyeq}\{ AXB \ | \  X \in {\cal S}_2 \}  = C - CM(M^{*}CM)^{\dag}M^{*}C,
\label{546w}
\\
& \max_{\succcurlyeq}\{ AXB - C \ | \  X \in {\cal S}_2 \}  =   - CM(M^{*}CM)^{\dag}M^{*}C.
\label{547w}
\end{align}

\item[{\rm(b)}] The  extremal ranks and partial inertias of $AXB$ and  $AXB -C$ subject to
$X \in {\cal S}_2$ are given by
\begin{align}
&\max_{X \in {\cal S}_2} r(AXB)  = \max_{X \in {\cal S}_2} i_{-}(AXB)    =r(A) + r(B) - r[\, A, \, B^*\,],
\label{548w}
\\
& \min_{X \in {\cal S}_2} r(AXB)  = \min_{X \in {\cal S}_2} i_{-}(AXB) =  r(A) + r(B) + i_{-}(C) -  i_{-}(N),
\label{549w}
\\
&\max_{X \in {\cal S}_2} i_{+}(AXB)  = r(A) + r(B) + i_{+}(C) -  i_{+}(N),
\label{550w}
\\
&\min_{X \in {\cal S}_2} i_{+}(AXB)   = 0,
\label{551w}
\\
&\max_{X \in {\cal S}_2} r(\, AXB - C\, )  =  r(N) - r(A) - r(B) - r[\, A, \, B^*\,],
\label{552w}
\\
& \min_{X \in {\cal S}_2} r(\, AXB - C \,)  = r\!\left[\!\! \begin{array}{cccc} A  & 0 & C
\\ 0 & B^{*} & C
\end{array} \!\!\right] - r(A) - r(B).
\label{553w}
\end{align}
\end{enumerate}
In consequence$,$
\begin{enumerate}
\item[{\rm(c)}] There exists an $X \in \mathbb C^{p\times q}$ such that
$ AXB \prec 0$  and $ AXB \preccurlyeq C$ if and only if $r[\, A, \, B^*\,]  = r(A) + r(B) - m.$

\item[{\rm(d)}] There exists an $X \in \mathbb C^{p\times q}$ such that
$ 0 \prec  AXB \preccurlyeq C$ if and only if $C\succ 0$ and $r\!\left[\!\! \begin{array}{cccc} A  & 0 & C
\\ 0 & B^{*} & C \end{array} \!\!\right] =  r(A) +  r(B).$

\item[{\rm(e)}] There exists an $X \in \mathbb C^{p\times q}$ such that
$ 0 \preccurlyeq AXB \preccurlyeq C$ if and only if $C\succcurlyeq 0.$

\item[{\rm(f)}] There always exists an $X \in \mathbb C^{p\times q}$ such that
$AXB \preccurlyeq 0$ and $AXB \preccurlyeq C.$
\end{enumerate}
\end{corollary}

In what follows, we give some consequences of Theorem \ref{T51} for different choice of
$C$ in (\ref{10b}).

\begin{theorem} \label{T45hh}
Let $A \in \mathbb C^{m\times p},$ $B \in \mathbb C^{q\times m}$ and
$C \in \mathbb C_{{\rm H}}^{m}$ be given$,$ and assume that $AXB = C$ is consistent$.$ Then$,$ the following hold$.$
\begin{enumerate}
\item[{\rm(a)}]  The general solution of $AXB \succcurlyeq C$  and the corresponding
$AXB$ can be written in the following parametric forms
\begin{align}
& X = A^{\dag}CB^{\dag}  + A^{\dag}E_MUU^{*}E_MB^{\dag}+ W - A^{\dag}AWBB^{\dag},
 \label{462hh}
\\
& AXB = C + E_MUU^{*}E_M,
 \label{463hh}
\end{align}
where $M = [\, E_A, \, F_B \,],$ and $U\in \mathbb C^{m\times m}$ and $W\in \mathbb C^{p\times q}$ are arbitrary$.$

\item[{\rm(b)}] There exists an $X \in \mathbb C^{p\times q}$ such that $AXB \succ  C$
if and only if $r(A) = r(B) = m.$ In this case$,$ the general solution of $AXB \succ  C$  and
the corresponding $AXB$ can be written as
\begin{align}
& X = A^{\dag}CB^{\dag}  + A^{\dag}UB^{\dag}+ W - A^{\dag}AWBB^{\dag},
 \label{464hh}
\\
& AXB = C + U,
 \label{465hh}
\end{align}
where $0\prec  U$ and $W\in \mathbb C^{p\times q}$ are arbitrary$.$

\item[{\rm(c)}] The general solution of $AXB \preccurlyeq C$  and the corresponding
$AXB$ can be written in the following parametric forms
\begin{align}
& X  = A^{\dag}CB^{\dag}  - A^{\dag}E_MUU^{*}E_MB^{\dag} + W - A^{\dag}AWBB^{\dag},
 \label{466h}
\\
& AXB = C - E_MUU^{*}E_M,
 \label{467hh}
\end{align}
where $U\in \mathbb C^{m\times m}$ and $W\in \mathbb C^{p\times q}$ are arbitrary$.$

\item[{\rm(d)}] There exists an $X \in \mathbb C^{p\times q}$ such that $AXB \prec  C$
if and only if $r(A) = r(B) = m.$ In this case$,$  the general solution of $AXB \succ  C$  and
the corresponding $AXB$ can be written in the following parametric forms
\begin{align}
& X = A^{\dag}CB^{\dag}  - A^{\dag}UB^{\dag}+ W - A^{\dag}AWBB^{\dag},
 \label{468hh}
\\
& AXB = C - U,
 \label{469hh}
\end{align}
where $0\prec  U$ and $W\in \mathbb C^{p\times q}$ are arbitrary$.$
\end{enumerate}
\end{theorem}

\begin{corollary} \label{T46hh}
Let $A \in \mathbb C^{m\times p},$ $B \in \mathbb C^{q\times m}$ and $C \in \mathbb C^{m \times m}$
be given$,$ and let $M = [\, E_A, \, F_B \,].$ Then$,$ the following hold$.$
\begin{enumerate}
\item[{\rm(a)}] The inequality
\begin{equation}
AXB  \succcurlyeq - CC^{*}
\label{515}
\end{equation}
is always feasible$;$ the general solution of {\rm (\ref{515})} and the corresponding $AXB$
can be written in the following parametric forms
\begin{align}
& X = -A^{\dag}CC^{*}B^{\dag} +
A^{\dag}C(M^{*}C)^{\dag}(M^{*}C)C^{*}B^{\dag} +
A^{\dag}E_MUU^{*}E_MB^{\dag} + W - A^{\dag}AWBB^{\dag},
 \label{516}
 \\
& AXB = - CC^{*} + C(M^{*}C)^{\dag}(M^{*}C)C^{*} + E_MUU^{*}E_M,
 \label{516a}
\end{align}
where $U\in \mathbb C^{m\times m}$ and $W\in \mathbb C^{p\times q}$ are arbitrary$.$

\item[{\rm(b)}] There exists an $X \in \mathbb C^{p\times q}$ such that
\begin{equation}
AXB  \succ  - CC^{*}
 \label{517}
\end{equation}
if and only if $r\!\left[\!\! \begin{array}{cccc} A  & 0 & C
\\ 0 & B^{*} & C \end{array} \!\!\right] = m +  r[\, A, \, B^{*} \,].$
In this case$,$ the general solution of {\rm (\ref{517})}
can be written as {\rm (\ref{516})}$,$ in which $U\in {\mathbb C}^{m\times m}$ is any matrix such that
$r[\, CC^*M, \,  E_MU \,] = m,$ and $W\in
\mathbb C^{p\times q}$ is arbitrary$.$

\item[{\rm(c)}] The inequality
\begin{equation}
AXB \preccurlyeq CC^{*}
 \label{518}
\end{equation}
is always feasible$;$ the general solution of {\rm (\ref{518})} and the corresponding $AXB$
can be written in the following parametric forms
\begin{align}
& X = A^{\dag}CC^{*}B^{\dag} - A^{\dag}C(M^{*}C)^{\dag}(M^{*}C)C^{*}B^{\dag} -
A^{\dag}E_MUU^{*}E_MB^{\dag}+ W - A^{\dag}AWBB^{\dag},
 \label{519}
 \\
& AXB = CC^{*} - C(M^{*}C)^{\dag}(M^{*}C)C^{*} + E_MUU^{*}E_M,
 \label{519a}
\end{align}
where $U\in \mathbb C^{m\times m}$ and $W\in \mathbb C^{p\times q}$ are arbitrary$.$

\item[{\rm(d)}] There exists an $X \in \mathbb C^{p\times q}$ such that
\begin{equation}
AXB  \prec  CC^{*}
\label{520}
\end{equation}
if and only if $r\!\left[\!\!\begin{array}{cccc} A  & 0 & C
\\ 0 & B^{*} & C \end{array} \!\!\right] = m +  r[\, A, \, B^{*} \,].$
 In this case$,$ the general solution of {\rm (\ref{520})}
can be written as {\rm (\ref{519})}$,$ in which $U$ is any matrix
such that $r[\,CC^*M, \,  E_MU \,] = m,$ and $W\in \mathbb C^{p\times q}$ is arbitrary$.$
\end{enumerate}
\end{corollary}

\begin{corollary}\label{T55}
Let $A \in \mathbb C^{m \times p}, \, B \in \mathbb C^{q \times m}$ and $C \in \mathbb C^{m \times m}$ be
given$,$ and let $M = [\, E_A, \, F_B \,].$
Then$,$ the following hold$.$
\begin{enumerate}
\item[{\rm (a)}] There exists an $X \in \mathbb C^{p\times q}$ such that
\begin{equation}
AXB \succcurlyeq CC^{*}
 \label{521}
\end{equation}
if and only if
\begin{equation}
\R(C) \subseteq \R(A) \ \ and \ \ \R(C) \subseteq \R(B^{*}).
 \label{522}
\end{equation}
In this case$,$ the general solution of {\rm (\ref{521})} and the corresponding $AXB$ can be written in the following parametric forms
\begin{align}
& X = A^{\dag}CC^{*}B^{\dag} + A^{\dag}E_MUU^{*}E_MB^{\dag}+ W -
A^{\dag}AWBB^{\dag},
 \label{523}
 \\
& AXB =CC^{*} + E_MUU^{*}E_M,
 \label{523a}
\end{align}
where $U\in \mathbb C^{m\times m}$ and $W\in \mathbb C^{p\times q}$ are arbitrary$.$

\item[{\rm (b)}] There exists an $X \in \mathbb C^{p\times q}$  such that
\begin{equation}
AXB \succ  CC^{*}
 \label{524}
\end{equation}
if and only if $r(A) = r(B) =m.$ In  this case$,$ the general solution of {\rm (\ref{524})} can be written as
{\rm (\ref{523})}$,$ in which $U\in \mathbb C^{q\times q}$ is any matrix with $r(E_{M}U) = m,$ and
 $W\in \mathbb C^{p\times q}$ is arbitrary$.$

\item[{\rm (c)}] There exists an $X \in \mathbb C^{p\times q}$ such that
\begin{equation}
AXB  \preccurlyeq - CC^{*}
\label{525}
\end{equation}
if and only if  {\rm (\ref{522})} holds$.$ In this case$,$  the general solution of {\rm (\ref{525})}
can be written in the following parametric forms
\begin{align}
& X = - A^{\dag}CC^{*}B^{\dag} - A^{\dag}E_MUU^{*}E_MB^{\dag}+ W -
A^{\dag}AWBB^{\dag},
 \label{526}
 \\
 & AXB = - CC^{*} - E_MUU^{*}E_M,
 \label{526a}
\end{align}
where $U\in \mathbb C^{m\times m}$ and $W\in \mathbb C^{p\times q}$ are arbitrary$.$

\item[{\rm (d)}] There exists an $X \in \mathbb C^{p\times q}$ such that
\begin{equation}
AXB \prec  - CC^{*}
 \label{527}
\end{equation}
if and only if $r(A) = r(B) =m.$ In  this case$,$  the general solution of {\rm (\ref{527})} can be written as
{\rm (\ref{526})}$,$ in which $U\in \mathbb C^{m \times m}$ is any matrix with $r(E_{M}U) = m,$ and
$W\in \mathbb C^{p\times q}$ is arbitrary$.$
\end{enumerate}
\end{corollary}

\begin{corollary}\label{T56}
Let $A \in \mathbb C^{m \times p}, \, B \in \mathbb C^{q \times
m}$ be given$,$ and let $M = [\, E_A, \, F_B \,].$ Then$,$ the following hold$.$
\begin{enumerate}
\item[{\rm (a)}] The general solution of
\begin{equation}
AXB \succcurlyeq 0
 \label{528}
\end{equation}
 and the corresponding $AXB$ can be written in the following parametric forms
\begin{align}
&  X=  A^{\dag}E_MUU^{*}E_MB^{\dag}+ W - A^{\dag}AWBB^{\dag},
 \label{529}
 \\
&  AXB = E_MUU^{*}E_M,
 \label{529xx}
\end{align}
where $U\in \mathbb C^{m\times m}$ and $W\in \mathbb C^{p\times q}$ are arbitrary$.$

\item[{\rm (b)}] There exists an $X \in \mathbb C^{p\times q}$ such that
\begin{equation}
AXB \succ  0
 \label{530}
\end{equation}
if and only if $r(A) = r(B) =m.$ In  this case$,$  the general solution of {\rm (\ref{530})}
can be written as {\rm (\ref{529})}$,$ in which $U\in \mathbb C^{m \times m}$ is any matrix such that
$r(E_{M}U) = m,$ and $W\in \mathbb C^{p\times q}$ is arbitrary$.$
\end{enumerate}
\end{corollary}

We next establish a group of formulas for calculating the ranks and inertias of
$AXB - D$ subject to  (\ref{52w}), and use the results obtained to derive necessary
and sufficient conditions for the following two-sides inequality
\begin{equation}
D \succcurlyeq  AXB \succcurlyeq C
\label{478jj}
\end{equation}
and their variations to hold.

\begin{corollary} \label{T48jj}
Let $A \in \mathbb C^{m\times p},$ $B \in \mathbb C^{q\times m}$ and
$C, \, D \in \mathbb C_{{\rm H}}^{m}$ be given$,$  and let ${\cal S}_1$ be of the forms in
 {\rm (\ref{536w})}$,$ and define
\begin{align}
& K_1 =\left[\!\begin{array}{cccccc} C & C & C & A & 0 \\  C & C  & C & 0 & B^*
\\
C & C & C-D & 0 & 0
\\
A^* & 0 & 0 & 0 & 0
\\
0 & B &  0 & 0 & 0
\end{array}\!\right]\!,   \ \  K_2 =\!\left[\!\! \begin{array}{cccc} D & D  & A & 0 \\
  D & D  & 0 & B^{*} \\ A^{*} & 0 & 0 & 0 \\  0 & B  & 0 & 0  \end{array} \!\!\right]\!,  \ \
K_3 = \left[\!\! \begin{array}{cccc} A  & 0 & D
\\ 0 & B^{*} & D
\end{array} \!\!\right]\!.
\label{479jj}
\end{align}
Then$,$ the extremal ranks and partial inertias of $AXB - D$ subject
  to $X \in {\cal S}_1$ are given by
\begin{align}
&\max_{X \in {\cal S}_1} r(\, AXB - D \,)  = r(K_3) - r[\, A, \, B^*\,],
\label{480jj}
\\
&\min_{X \in {\cal S}_1} r(\, AXB - D \,)  =  i_{+}(K_1) +  r(K_3) - r(K_2),
\label{381jj}
\\
&\max_{X \in {\cal S}_1} i_{+}(\, AXB - D \,)  =i_{-}(K_2) - r[\, A, \, B^*\,],
\label{482jj}
\\
&\max_{X \in {\cal S}_1} i_{-}(\, AXB - D \,)  =i_{-}(K_1) - i_{-}(K_2),
\label{483jj}
\\
&\min_{X \in {\cal S}_1} i_{+}(\, AXB - D \,)  =i_{+}(K_1) - i_{+}(K_2),
\label{484jj}
\\
&\min_{X \in {\cal S}_1} i_{-}(\, AXB - D \,)  = r(K_3) -  i_{-}(K_2).
\label{485jj}
\end{align}
In consequence$,$ the following hold$.$
\begin{enumerate}
\item[{\rm(a)}] There exists an $X \in \mathbb C^{p\times q}$ such that
$ AXB \succ D$  and $ AXB \succcurlyeq C$ if and only if $i_{-}(K_2) = r[\, A, \, B^*\,] + m.$

\item[{\rm(b)}] There exists an $X \in \mathbb C^{p\times q}$ such that
$ D \succ  AXB \succcurlyeq C$ if and only if $D\succ  C$ and $i_{-}(K_1) = i_{-}(K_2) + m.$

\item[{\rm(c)}] There exists an $X \in \mathbb C^{p\times q}$ such that
$ D \succcurlyeq AXB \succcurlyeq C$ if and only if $D\succcurlyeq C$ and  $i_{+}(K_1) = i_{+}(K_2).$

\item[{\rm(d)}] There exists an $X \in \mathbb C^{p\times q}$ such that
$ AXB \succcurlyeq C$ and $ AXB \succcurlyeq D$ f and only if
 $r(K_3) = i_{-}(K_2).$

\end{enumerate}

\end{corollary}

\noindent {\bf Proof.} \  From (\ref{517w}), $AXB - D$ subject to $X \in {\cal S}_1$
can be written as
\begin{align}
AXB  -  D =  C - D - CM(M^{*}CM)^{\dag}M^{*}C + E_MUU^{*}E_M = A\widehat{X}B  - D + E_MUU^{*}E_M.
 \label{486jj}
\end{align}
Applying elementary matrix operations, congruence matrix operations and (\ref{217}), we obtain
\begin{align}
 r(E_M) & = m - r[\,E_A, \,F_B \,] = m - r(E_A) - r(F_BA)  \nb
 \\
 & = r(A) + r(B) - r[\, A, \, B^*\,],
\label{487jj}
\\
 r[\, E_M, \, A\widehat{X}B - D\,] & = r(E_M) + r[\,(A\widehat{X}B - D)M \,] =
 r(E_M) + r[\,(\, B^*\widehat{X}^*A^* -  D \,)E_A, \, (\, A\widehat{X}B - D \,)F_B \,] \nb
 \\
 & =  r(E_M) + r[\, DE_A, \, DF_B \,] =  r\!\left[\!\! \begin{array}{cccc} A  & 0 & D
\\ 0 & B^{*} & D
\end{array} \!\!\right] - r[\, A, \, B^*\,] \nb
\\
& = r(K_3) - r[\, A, \, B^*\,],
\label{488jj}
\\
 i_{\pm}\!\left[\!\! \begin{array}{cc} A\widehat{X}B - D & E_M \\  E_M
 & 0 \end{array}\!\!\right] & =  i_{\pm}[\,M(A\widehat{X}B - D)M \,] + r(E_M)
 = i_{\mp}(MDM) + r(E_M)  \nb
 \\
 & = i_{\mp}(K_2) - r[\, A, \, B^*\,],
\label{489jj}
\\
i_{\pm}(\, A\widehat{X}B - D\,) &  = i_{\pm}[\, C  - D - CM(M^{*}CM)^{\dag}M^{*}C \,]
= i_{\pm}\!\left[\!\begin{array}{cc} M^{*}CM &  M^{*}C \\
 CM   & C-D \end{array}\!\right] - i_{\pm}(M^{*}CM) \nb
 \\
& =  i_{\pm}\!\left[\!\begin{array}{cccccc} C & C & C & A & 0 \\  C & C  & C & 0 & B^*
\\
C & C & C-D & 0 & 0
\\
A^* & 0 & 0 & 0 & 0
\\
0 & B &  0 & 0 & 0
\end{array}\!\right] - i_{\pm}\!\left[\!\begin{array}{cccccc} C & C & A & 0 \\  C & C  & 0 & B^*
\\
A^* & 0 & 0 & 0
\\
0 & B &  0 & 0
\end{array}\!\right]  \nb
\\
& =  i_{\pm}(K_1) - i_{\pm}(K_2).
\label{490jj}
\end{align}
Applying (\ref{278})--(\ref{283}) to (\ref{486jj}) and simplifying by
(\ref{487jj})--(\ref{490jj}), we obtain
 \begin{align*}
\max_{X \in {\cal S}_1} r(\,AXB - D\,) & = \max_{U\in {\mathbb C}^{m \times m}} r(\, A\widehat{X}B - D + E_MUU^{*}E_M \,)
 = r[\, E_M, \, A\widehat{X}B  - D\,]
 \\
 & = r(K_3) - r[\, A, \, B^*\,],
\\
\min_{X \in {\cal S}_1} r(\,AXB - D\,) & = \min_{U\in {\mathbb C}^{m \times m}} r(\, A\widehat{X}B - D +  E_MUU^{*}E_M  \,) \nb
\\
&= i_{+}(\, A\widehat{X}B - D\,) +  r[\,  E_M, \, A\widehat{X}B - D \,] - i_{+}\!\left[\!\! \begin{array}{cc}
A\widehat{X}B  - D &  E_M  \\  E_M & 0 \end{array}
\!\!\right] \nb
\\
& = i_{+}(K_1) +  r(K_3) - r(K_2),
\\
 \max_{X \in {\cal S}_1} i_{+}(\,AXB - D\,) & = \max_{U\in {\mathbb C}^{m \times m}} i_{+}(\, A\widehat{X}B - D +
  E_MUU^{*}E_M  \,)
 =  i_{+}\!\left[\!\! \begin{array}{cc} A\widehat{X}B  - D &  E_M \\  E_M & 0 \end{array}
\!\!\right]  \nb
\\
&   =  i_{-}(K_2) - r[\, A, \, B^*\,],
\\
 \max_{X \in {\cal S}_1} i_{-}(\,AXB - D\,) & =
 \max_{U\in {\mathbb C}^{m \times m}} i_{-}(\, A\widehat{X}B  - D +   E_MUU^{*}E_M  \,)
  =  i_{-}(\, A\widehat{X}B - D\,)
  \\
  &  =  i_{-}(K_1) - i_{-}(K_2),
 \\
 \min_{X \in {\cal S}_1} i_{+}(\,AXB - D\,) & = \min_{U\in {\mathbb C}^{m \times m}} i_{+}(\, A\widehat{X}B - D +
  E_MUU^{*}E_M  \,)  = i_{+}(\,A\widehat{X}B - D\,)
  \\
  &  = i_{+}(K_1) - i_{+}(K_2),
\\
 \min_{X \in {\cal S}_1} i_{-}(\,AXB - D\,) & =
 \min_{U\in {\mathbb C}^{m \times m}} i_{-}(\, A\widehat{X}B - D +  E_MUU^{*}E_M\,) \nb
 \\
 & = r[\, E_M, \, A\widehat{X}B - D \,] -  i_{+}\!\left[\!\! \begin{array}{cc} A\widehat{X}B  - D & E_M \\  E_M
 & 0 \end{array}
\!\!\right]
\\
& = r(K_3) -  i_{-}(K_2),
\end{align*}
as required for  (\ref{480jj})--(\ref{485jj}). \qquad $\Box$

\section{General Hermitian solution of the LMI $AXA^{*} \succcurlyeq \,(\succ, \, \preccurlyeq, \, \prec ) \, B$ and its properties}
\renewcommand{\theequation}{\thesection.\arabic{equation}}
\setcounter{section}{4}
\setcounter{equation}{0}

The LMIs in (\ref{10d}) are the simplest case of all LMIs with symmetric pattern.
Due to the importance of matrix inequalities in the L\"owner partial ordering,
any contribution on this type of LMIs is valuable from both
theoretical and practical points of view.  Some previous work
on solvability and general solutions of (\ref{10d}) and
their applications in system and control theory were
given in \cite{SIG} by using SVDs of matrices. In a recent paper
\cite{T-laa10},  necessary and sufficient conditions for the
LMIs in (\ref{10d}) to hold were obtained by using some expansion formulas
for the inertia of the matrix function $B - AXA^{*}$, while general Hermitian solution of $AXA^{*} \preccurlyeq B$
was established in  \cite{T-mia}.  In this section, we reconsider (\ref{10d})
and give a group of complete conclusions on Hermitian solutions of  the LMIs and their algebraic properties.

\begin{theorem} \label{T41}
Let $A \in \mathbb C^{m\times n}$ and  $B\in \mathbb C_{{\rm H}}^{m}$ be given$,$  and let
$N = \left[\!\!\begin{array}{cc}  B & A
\\ A^* & 0 \end{array} \!\!\right].$ Then$,$
\begin{enumerate}
\item[{\rm(a)}]  The following statements are equivalent$:$
\begin{enumerate}
\item[{\rm(i)}] There exists an $X \in {\mathbb C}_{{\rm H}}^{n}$ such that
\begin{equation}
AXA^{*} \succcurlyeq B.
\label{42}
\end{equation}

\item[{\rm(ii)}] $E_ABE_A \preccurlyeq 0$ and  $\R(E_ABE_A)= \R(E_AB).$

\item[{\rm(iii)}] $i_{+}(N) = r(A)$ and $i_{-}(N) = r[\, A, \, B\,].$
\end{enumerate}
In this case$,$ the general Hermitian solution of {\rm (\ref{42})} and  the corresponding $AXA^*$ can be written in the following parametric forms
\begin{align}
& X  = A^{\dag}B(A^{\dag})^* - A^{\dag}BE_A(E_ABE_A)^{\dag}E_AB(A^{\dag})^* + UU^*
 + W - A^{\dag}AWA^{\dag}A,
\label{43}
\\
& AXA^* = B - BE_A(E_ABE_A)^{\dag}E_AB +  AUU^*A^*,
\label{43a}
\end{align}
where $U\in \mathbb C^{n\times n}$ and $W \in {\mathbb C}_{{\rm H}}^{n}$ are arbitrary$.$

\item[{\rm(b)}] There exists an $X \in {\mathbb C}_{{\rm H}}^{n}$ such that
\begin{equation}
AXA^* \succ  B
\label{44}
\end{equation}
if and only if
\begin{equation}
E_ABE_A \preccurlyeq 0 \ \ \ and  \ \  r(E_ABE_A)= r(E_A).
\label{45}
\end{equation}
In this case$,$ the general Hermitian solution of {\rm (\ref{42})} can be
written as {\rm (\ref{43})}$,$ in which $U$ is any matrix such that
$r[\, BE_A, \  AU\,] =m,$  say$,$  $U = I_n,$  and $W  \in {\mathbb C}_{{\rm H}}^{n}$ is arbitrary$.$

\item[{\rm(c)}] {\rm \cite{T-mia}} The following statements are equivalent$:$
\begin{enumerate}
\item[{\rm(i)}]  There exists an $X \in {\mathbb C}_{{\rm H}}^{n}$ such that
\begin{equation}
AXA^{*} \preccurlyeq B.
\label{46}
\end{equation}

\item[{\rm(ii)}] $E_ABE_A \succcurlyeq 0$ and $\R(E_ABE_A)= \R(E_AB).$

\item[{\rm(iii)}] $i_{+}(N) = r[\, A, \, B\,]$ and $i_{-}(N) = r(B).$
\end{enumerate}
In this case$,$ the general Hermitian solution of  {\rm (\ref{46})} and the corresponding
$AXA^*$ can be written in the following parametric forms
\begin{align}
& X =  A^{\dag}B(A^{\dag})^* - A^{\dag}BE_A(E_ABE_A)^{\dag}E_AB(A^{\dag})^* - UU^*
 + W - A^{\dag}AWA^{\dag}A,
\label{47}
\\
& AXA^* = B - BE_A(E_ABE_A)^{\dag}E_AB -  AUU^*A^*,
\label{47a}
\end{align}
where $U \in \mathbb C^{n \times n}$ and $W \in {\mathbb C}_{{\rm H}}^{n}$ are arbitrary$.$

\item[{\rm(d)}] There exists an $X \in {\mathbb C}_{{\rm H}}^{n}$ such that
\begin{equation}
AXA^* \prec B
\label{48}
\end{equation}
if and only if {\rm (\ref{45})} holds$.$ In this case$,$ the general Hermitian solution of
{\rm (\ref{48})} can be written as {\rm (\ref{47})}$,$ in which $U$ is any matrix such that
$r[\, BE_A, \  AU\,] =m,$ say$,$  $U = I_n,$ and $W \in {\mathbb C}_{{\rm H}}^{n}$
is arbitrary$.$
\end{enumerate}
\end{theorem}

\noindent {\bf Proof.} \
Inequality (\ref{42}) can be relaxed to the following quadratic matrix equation
\begin{equation}
AXA^* = B + YY^*.
\label{49}
\end{equation}
 By Lemma \ref{T26}(a), (\ref{49}) is solvable for $X$ if and only if $E_A(B + YY^*) = 0,$ that is,
\begin{equation}
E_AYY^* = -E_AB.
\label{410}
\end{equation}
By Lemma \ref{T25}(b), (\ref{410}) is solvable for $YY^*$ if and only if $E_ABE_A \preccurlyeq 0$
and $r(E_ABE_A)= r(E_AB)$, establishing the equivalence (i) and (ii) in (a).  The
equivalence (ii) and (iii) in (a) follows from (\ref{21x}) and $i_{-}(E_ABE_A) \preccurlyeq
r(E_ABE_A) \preccurlyeq r(E_AB)$. In this case, the general solution of (\ref{410}) can be written as
$$
YY^* =  -BE_A(E_ABE_A)^{\dag}E_AB +  AA^{\dag}UU^*AA^{\dag},
$$
where $U$ is an arbitrary matrix. Substituting the $YY^*$ into (\ref{49}) gives
\begin{equation}
AXA^* = B - BE_A(E_ABE_A)^{\dag}E_AB +  AA^{\dag}UU^*AA^{\dag}.
\label{411}
\end{equation}
By Lemma \ref{T26}(a), the general Hermitian solution of (\ref{411}) can be written as
$$
X = A^{\dag}B(A^{\dag})^* - A^{\dag}BE_A(E_ABE_A)^{\dag}E_AB(A^{\dag})^* + A^{\dag}UU^*(A^{\dag})^*
 + W - A^{\dag}AWA^{\dag}A,
$$
where $U \in \mathbb C^{m \times m}$ and $W \in {\mathbb C}_{{\rm H}}^{n}$ are arbitrary$.$ Replacing
 $A^{\dag}UU^*(A^{\dag})^*$ with  $UU^*$ gives (\ref{43}), which is also the general
 solution of  (\ref{42}).

 It can be seen from (\ref{43a}) that (\ref{44}) holds if and only if
\begin{equation}
-BE_A(E_ABE_A)^{\dag}E_AB +  AUU^*A^*\succ 0
\label{414yy}
\end{equation}
for some $U$. Under (ii) in (a),  we have
\begin{equation*}
r[\, -BE_A(E_ABE_A)^{\dag}E_AB +  AUU^*A^* \,] =
r[\, -BE_A(E_ABE_A)^{\dag}E_AB,  \  AUU^*A^* \,] = r[\, BE_A,  \ AU \,].
\end{equation*}
Hence,
\begin{align*}
\max_{U} r[\, -BE_A(E_ABE_A)^{\dag}E_AB +  AUU^*A^* \,]
 = \max_{U} r[\, BE_A,  \ AU \,]= r[\, BE_A, \,   A \,] = r(E_ABE_A) + r(A),
\end{align*}
so that (\ref{44}) holds if and only if $r(E_ABE_A) + r(A) = m.$ Thus (b) follows.
Results (c) and (d) can be shown similarly. \qquad $\Box$

\medskip

Concerning the constant term in (\ref{43}), we have the consequence.

\begin{corollary} \label{T42J}
Let  $A \in \mathbb C^{m\times n}$ and  $B\in \mathbb C_{{\rm
H}}^{m}$  be given$,$ and let
\begin{equation}
\widehat{X} = A^{\dag}B(A^{\dag})^* - A^{\dag}BE_A(E_ABE_A)^{\dag}E_AB(A^{\dag})^*.
\label{415yy}
\end{equation}
Then$,$ the following hold$.$
\begin{enumerate}
\item[{\rm(a)}] Under the condition that {\rm (\ref{42})} is feasible$,$
$\widehat{X}$ is a Hermitian solution of {\rm (\ref{42})}$,$ and
\begin{align}
& i_{+}(\widehat{X}) = i_{+}(A\widehat{X}A^*) = i_{+}(B),
\label{416yy}
\\
& i_{-}(\widehat{X}) =  i_{-}(A\widehat{X}A^*) = r(A) + i_{-}(B) - r[\, A, \, B \,],
\label{417yy}
\\
& r(\widehat{X}) =  r(A\widehat{X}A^*) = r(A) + r(B) - r[\, A, \, B \,],
\label{418yy}
\\
& i_{-}(\,B - A\widehat{X}A^*\,) = r(\,B - A\widehat{X}A^*\,) =
 r(B) - r(\,A\widehat{X}A^*\,) = r[\, A, \, B\,] - r(A).
 \label{419yy}
\end{align}

\item[{\rm(b)}] Under the condition that {\rm (\ref{46})} is feasible$,$
$\widehat{X}$ is a Hermitian solution of ${\rm (\ref{46})},$ and
\begin{align}
& i_{+}(\widehat{X}) =  i_{+}(A\widehat{X}A^*) =
r(A) + i_{+}(B)  - r[\, A, \, B \,],
\label{420yy}
 \\
& i_{-}(\widehat{X}) =  i_{-}(A\widehat{X}A^*) = i_{-}(B),
\label{421yy}
\\
& r(\widehat{X}) = r(A\widehat{X}A^*) = r(A) + r(B) - r[\, A, \, B \,],
\label{422yy}
\\
& i_{+}(\,B - A\widehat{X}A^*\,) = r(\,B - A\widehat{X}A^*\,)
= r(B) - r(\,A\widehat{X}A^*\,) = r[\, A, \, B\,] - r(A).
\label{423yy}
\end{align}
\end{enumerate}
\end{corollary}

\noindent {\bf Proof.} \ Under the condition that (\ref{42}) has a solution$,$
set $U = W = 0$ in  (\ref{43}),  we see that $\widehat{X}$ in  (\ref{415yy})
is a Hermitian solution of $AXA^* \succcurlyeq B.$ In this case, applying (\ref{113j}) to
(\ref{415yy}) and simplifying  by congruence matrix operations,  we obtain
\begin{align}
i_{\pm}(\widehat{X}) & = i_{\pm}[\,  A^{\dag}B(A^{\dag})^* - A^{\dag}BE_A(E_ABE_A)^{\dag}E_AB(A^{\dag})^* \,] \nb
\\
& = i_{\pm}\!\left[\!\begin{array}{cc} E_ABE_A &  E_AB(A^{\dag})^* \\
A^{\dag}BE_A & A^{\dag}B(A^{\dag})^* \end{array}\!\right] - i_{\pm}(E_ABE_A) \nb
\\
& = i_{\pm}\!\left[\!\begin{array}{cc} B &  B(A^{\dag})^* \\
A^{\dag}B & A^{\dag}B(A^{\dag})^* \end{array}\!\right] - i_{\pm}(E_ABE_A)  \nb
\\
&= i_{\pm}\!\left[\!\begin{array}{cc} B &  0 \\
0 & 0 \end{array}\!\right] - i_{\pm}(E_ABE_A)  = i_{\pm}(B) - i_{\pm}(E_ABE_A),
\label{424yy}
\\
i_{\pm}(A\widehat{X}A^*) & = i_{\pm}[\, B - BE_A(E_ABE_A)^{\dag}E_AB \,] \nb
\\
& = i_{\pm}\!\left[\!\begin{array}{cc} E_ABE_A &  E_AB \\
BE_A & B \end{array}\!\right] - i_{\pm}(E_ABE_A) \nb
\\
& = i_{\pm}\!\left[\!\begin{array}{cc} 0 &  0 \\
0 & B \end{array}\!\right] - i_{\pm}(E_ABE_A)  = i_{\pm}(B) - i_{\pm}(E_ABE_A).
\label{425yy}
\end{align}
In consequence,
\begin{align*}
&i_{+}(\widehat{X}) = i_{+}(A\widehat{X}A^*)  = i_{+}(B), \ \
\\
& i_{-}(\widehat{X}) = i_{-}(A\widehat{X}A^*) = i_{-}(B) -
i_{-}(E_ABE_A) = i_{-}(B) -  r(E_AB) = i_{-}(B) + r(A) - r[\, A, \, B \,],
\end{align*}
establishing (\ref{416yy}), (\ref{417yy}) and (\ref{418yy}).
 Applying (\ref{113j}) and simplifying  by congruence matrix operations,
we obtain
\begin{align}
i_{\pm}(\,B - A\widehat{X}A^*\,) & = i_{\pm}[\,BE_A(E_ABE_A)^{\dag}E_AB \,] \nb
\\
& = i_{\pm}\!\left[\!\begin{array}{cc} - E_ABE_A &  E_A B \\
BE_A &  0 \end{array}\!\right] - i_{\mp}(E_ABE_A) \nb
\\
& = i_{\pm}\!\left[\!\begin{array}{cc} 0 &  E_AB \\
BE_A & 0 \end{array}\!\right] - i_{\mp}(E_ABE_A)  = r(E_AB) - i_{\mp}(E_ABE_A).
\label{426yy}
\end{align}
In consequence,
\begin{align*}
& i_{+}(\, B - A\widehat{X}A^* \,) = r(E_AB) - i_{-}(E_ABE_A) =r(E_AB) - r(E_ABE_A) = 0,
\\
& i_{-}(\,B - A\widehat{X}A^*\,) = r(E_AB) - i_{+}(E_ABE_A) =r[\, A, \, B\,] - r(A),
\end{align*}
establishing  (\ref{419yy}).   Result (b) can be shown similarly.
\qquad $\Box$

\begin{corollary} \label{T43J}
Let $A \in {\mathbb C}^{m\times n}$ and  $B\in \mathbb C_{{\rm H}}^{m}$ be given$.$
Then$,$ the following hold$.$
\begin{enumerate}
\item[{\rm(a)}] Under the condition that {\rm (\ref{42})} is feasible$,$ define
\begin{equation}
{\cal S}_1 =\{\, X\in {\mathbb C}_{{\rm H}}^{n} \ | \  AXA^* \succcurlyeq  B \, \}.
\label{427yy}
\end{equation}
Then$,$ the minimal matrices of $AXA^*$ and $AXA^* - B$ subject
  to $X \in {\cal S}_1$ in the L\"owner partial ordering are given by
\begin{align}
&\min_{\succcurlyeq}\{ AXA^* \ | \  X \in {\cal S}_1 \}  =  B - BE_A(E_ABE_A)^{\dag}E_AB,
\label{428yy}
\\
& \min_{\succcurlyeq}\{ AXA^* - B \ | \  X \in {\cal S}_1 \}  =  - BE_A(E_ABE_A)^{\dag}E_AB,
\label{429yy}
\end{align}
while the  extremal ranks and partial inertias of $AXA^{*}$ and $AXA^{*} -B$ subject
  to $X \in {\cal S}_1$ are given by
\begin{align}
&\max_{X \in {\cal S}_1} r(AXA^*)  = \max_{X \in {\cal S}_1} i_{+}(AXA^*)  = r(A),
\label{430yy}
\\
& \min_{X \in {\cal S}_1} r(AXA^*)  = \min_{X \in {\cal S}_1} i_{+}(AXA^*) = i_{+}(B),
\label{431yy}
\\
&\max_{X \in {\cal S}_1} i_{-}(AXA^*) = r(A) + i_{-}(B) - r[\, A, \, B \,],
\label{432yy}
\\
&\min_{X \in {\cal S}_1} i_{-}(AXA^*)  = 0,
\label{433yy}
\\
&\max_{X \in {\cal S}_1} r(\,AXA^{*} - B\,)  = r[\, A, \, B \,],
\label{434yy}
\\
& \min_{X \in {\cal S}_1} r(\, AXA^{*} - B  \,)  = r[\, A, \, B \,] - r(A).
\label{435yy}
\end{align}

\item[{\rm(b)}]  Under the condition that {\rm (\ref{46})} is feasible$,$ and define
\begin{equation}
{\cal S}_2 =\{\, X\in {\mathbb C}_{{\rm H}}^{n} \ | \  AXA^* \preccurlyeq  B \, \}.
\label{436yy}
\end{equation}
Then$,$ the maximal matrices of $AXA^*$ and $AXA^* - B$ subject to $X \in {\cal S}_2$
in the L\"owner partial ordering are given by
\begin{align}
&\max_{\succcurlyeq}\{ AXA^* \ | \  X \in {\cal S}_2 \}  =  B - BE_A(E_ABE_A)^{\dag}E_AB,
\label{437yy}
\\
& \max_{\succcurlyeq}\{ AXA^* - B \ | \  X \in {\cal S}_2 \}  =  - BE_A(E_ABE_A)^{\dag}E_AB,
\label{438yy}
\end{align}
while
the  extremal ranks and partial inertias of $AXA^{*}$ and  $AXA^{*} -B$ subject to $X \in {\cal S}_2$ are given by
\begin{align}
&\max_{X \in {\cal S}_2} r(AXA^*)  = \max_{X \in {\cal S}_2} i_{-}(AXA^*)  = r(A),
\label{439yy}
\\
& \min_{X \in {\cal S}_2} r(AXA^*)  = \min_{X \in {\cal S}_2} i_{-}(AXA^*) = i_{-}(B),
\label{440yy}
\\
&\max_{X \in {\cal S}_2} i_{+}(AXA^*)  = r(A) + i_{+}(B) - r[\, A, \, B \,],
\label{441yy}
\\
&\min_{X \in {\cal S}_2} i_{+}(AXA^*)  = 0,
\label{442yy}
\\
&\max_{X \in {\cal S}_2} r(AXA^{*} - B)  = r[\, A, \, B \,],
\label{443yy}
\\
& \min_{X \in {\cal S}_2} r(\, AXA^{*} - B  \,)  = r[\, A, \, B \,] - r(A).
\label{444yy}
\end{align}
\end{enumerate}
\end{corollary}

\begin{corollary} \label{T54hh}
Let $A \in \mathbb C^{m\times n}$ and  $B\in \mathbb C_{{\rm
H}}^{m}$ be given$,$  and assume that $AXA^{*} = B$ is consistent$.$ Then$,$ the following hold$.$
\begin{enumerate}
\item[{\rm(a)}]  The general Hermitian solution of $AXA^{*} \succcurlyeq B$ and the corresponding $AXA^*$ can be written
in the following parametric forms
\begin{align}
& X  = A^{\dag}B(A^{\dag})^*  + UU^* + W - A^{\dag}AWA^{\dag}A,
\label{545ff}
\\
& AXA^* = B +  AUU^*A^*,
\label{546ff}
\end{align}
where $U\in \mathbb C^{n\times n}$ and $W \in {\mathbb C}_{{\rm H}}^{n}$ are arbitrary$.$

\item[{\rm(b)}] There exists an $X \in {\mathbb C}_{{\rm H}}^{n}$ such that $AXA^* \succ  B$
if and only if $r(A) =m.$ In this case$,$ the general Hermitian solution  $AXA^{*} \succcurlyeq B$ and the corresponding $AXA^*$
can be written in the following parametric forms
\begin{align}
& X  = A^{\dag}B(A^{\dag})^*  + UU^* + W - A^{\dag}AWA^{\dag}A,
\label{547ff}
\\
& AXA^* = B +  AUU^*A^*,
\label{548ff}
\end{align}
where $U\in \mathbb C^{n\times n}$ is any matrix such that $r(AU) = m$ and $W \in {\mathbb C}_{{\rm H}}^{n}$ is arbitrary$.$

\item[{\rm(c)}]
 The general Hermitian solution of $AXA^{*} \preccurlyeq B$ and the corresponding $AXA^*$ can be written
 in the following parametric forms
\begin{align}
& X =  A^{\dag}B(A^{\dag})^*  - UU^* + W - A^{\dag}AWA^{\dag}A,
\label{549ff}
\\
& AXA^* = B -  AUU^*A^*,
\label{550ff}
\end{align}
where $U \in \mathbb C^{n \times n}$ and $W \in {\mathbb C}_{{\rm H}}^{n}$ are arbitrary$.$

\item[{\rm(d)}] There exists an $X \in {\mathbb C}_{{\rm H}}^{n}$ such that $AXA^* \prec B$ if
and only if $r(A) =m.$ In this case$,$
the general Hermitian solution of $AXA^{*} \prec  B$ and the corresponding $AXA^*$ can be written in the following parametric forms
\begin{align}
& X =  A^{\dag}B(A^{\dag})^*  - UU^* + W - A^{\dag}AWA^{\dag}A,
\label{551ff}
\\
& AXA^* = B -  AUU^*A^*,
\label{552ff}
\end{align}
where $U\in \mathbb C^{n\times n}$ is any matrix such that $r(AU) = m$ and  and $W \in {\mathbb C}_{{\rm H}}^{n}$ is arbitrary$.$
\end{enumerate}
\end{corollary}

\begin{theorem} \label{T44J}
Let $A \in \mathbb C^{m\times n}$ and  $B\in \mathbb C^{m\times m}$ be given$.$ Then$,$ the following hold$.$
\begin{enumerate}
\item[{\rm(a)}] The inequality
\begin{equation}
AXA^* \succcurlyeq - BB^*
\label{413}
\end{equation}
is always feasible$;$ the general Hermitian solution of {\rm (\ref{413})} and the corresponding $AXA^*$ can be written in the following parametric forms
\begin{align}
& X = A^{\dag}B(E_AB)^{\dag}(E_AB)B^*(A^{\dag})^* - A^{\dag}BB^*(A^{\dag})^*  + UU^* +
W - A^{\dag}AWA^{\dag}A,
\label{414}
\\
& AXA^* = B(E_AB)^{\dag}(E_AB)B^* - BB^*  + AUU^*A^*,
\label{414a}
\end{align}
where $U\in \mathbb C^{n \times n}$ and $W \in {\mathbb C}_{{\rm H}}^{n}$ are arbitrary$.$

\item[{\rm(b)}] There exists an $X \in \mathbb C^{n\times n}$ such
that
\begin{equation}
AXA^* \succ  -BB^*
\label{415}
\end{equation}
if and only if  $r[\,A, \,B \,]= m.$ In this case$,$ the general Hermitian solution of {\rm (\ref{415})}  can be
written as {\rm (\ref{414})}$,$ in which $U$ is any matrix such that
$r(AU)=r(A),$ say$,$ $U = I_n,$  and $W \in {\mathbb C}_{{\rm H}}^{n}$ is arbitrary$.$

\item[{\rm(c)}]  The inequality
\begin{equation}
AXA^* \preccurlyeq BB^*
\label{416}
\end{equation}
is always feasible$;$ the general Hermitian solution of  {\rm (\ref{416})}  and the corresponding $AXA^*$ can be written in the following parametric forms
\begin{align}
& X = A^{\dag}BB^*(A^{\dag})^* -A^{\dag}B(E_AB)^{\dag}(E_AB)B^*(A^{\dag})^* - UU^* +
W - A^{\dag}AWA^{\dag}A,
\label{417}
\\
& AXA^* = BB^* -  B(E_AB)^{\dag}(E_AB)B^* - AUU^*A^*,
\label{417a}
\end{align}
where $U\in \mathbb C^{n \times n}$ and $W\in {\mathbb C}_{{\rm H}}^{n}$ are arbitrary$.$

\item[{\rm(d)}] There exists an $X \in \mathbb C^{n\times n}$ such
that
\begin{equation}
AXA^* \prec  BB^*
\label{418}
\end{equation}
if and only if  $r[\,A, \,B \,]= m.$ In this case$,$ the general Hermitian solution of  {\rm (\ref{418})}  can be
written as {\rm (\ref{417})}$,$ in which $U$ is any matrix such that
$r(AU) =r(A),$ say$,$ $U = I_n,$ and  $W\in {\mathbb C}_{{\rm H}}^{n}$ is arbitrary$.$
\end{enumerate}
\end{theorem}

\begin{corollary} \label{T45}
Let $A \in \mathbb C^{m\times n}$ and  $B \in \mathbb C^{m\times m}$ be given$.$ Then$,$ the following hold$.$
\begin{enumerate}
\item[{\rm(a)}] There exists an $X \in \mathbb C^{n\times n}$ such
that
\begin{equation}
AXA^* \succcurlyeq BB^*
\label{419}
\end{equation}
if and only if ${\mathscr R}(B) \subseteq {\mathscr R}(A).$ In this case$,$
the general Hermitian solution and the corresponding $AXA^*$ can be written in the following parametric forms
\begin{align}
& X = A^{\dag}BB^*(A^{\dag})^* + UU^* +  W - A^{\dag}AWA^{\dag}A,
\label{420}
\\
& AXA^* = BB^* + AUU^*A^*,
\label{420a}
\end{align}
where $U\in \mathbb C^{n\times n}$ and $W\in {\mathbb C}_{{\rm H}}^{n}$ are arbitrary$.$

\item[{\rm(b)}] There exists an $X \in \mathbb C^{n\times n}$ such
that
\begin{equation}
AXA^* \succ  BB^*
\label{421}
\end{equation}
if and only if  $r(A)= m.$ In this case$,$ the general Hermitian solution of {\rm (\ref{421})}
can be written as {\rm (\ref{420})}$,$  in which
$U\in {\mathbb C}^{n\times n}$ is any matrix such that $r(AU) =m,$ and
 $W\in {\mathbb C}_{{\rm H}}^{n}$ is arbitrary$.$

\item[{\rm(c)}] There exists an $X \in {\mathbb C}_{{\rm H}}^{n}$ such that
\begin{equation}
AXA^*  \preccurlyeq - BB^*
\label{422}
\end{equation}
if and only if  ${\mathscr R}(B) \subseteq {\mathscr R}(A).$ In this case$,$
the general Hermitian solution of {\rm (\ref{422})} and the corresponding $AXA^*$ can be written in the following parametric forms
\begin{align}
& X = -A^{\dag}BB^*(A^{\dag})^*  - UU^* +  W - A^{\dag}AWA^{\dag}A,
\label{423}
\\
& AXA^* = -BB^*  - AUU^*A^*,
\label{423a}
\end{align}
where $U\in \mathbb C^{n\times n}$ and  $W\in {\mathbb C}_{{\rm H}}^{n}$ are arbitrary$.$

\item[{\rm(d)}] There exists an  $X \in {\mathbb C}_{{\rm H}}^{n}$ such that
\begin{equation}
AXA^* \prec   - BB^*
\label{424}
\end{equation}
if and only if  $r(A) =m.$ In this case$,$ the general Hermitian solution of {\rm (\ref{424})} can be written as
 {\rm (\ref{423})}$,$  in which $U\in \mathbb C^{n\times n}$ is any matrix such that  $r(AU) =m,$ and
$W\in  {\mathbb C}_{{\rm H}}^{n}$ is arbitrary$.$
\end{enumerate}
\end{corollary}

\begin{corollary} [\cite{T-laa10}] \label{T46}
Let $A \in \mathbb C^{m\times n}$ be given$.$ Then$,$ the following hold$.$
\begin{enumerate}
\item[{\rm(a)}] The general solution of  $AXA^{*} \succcurlyeq 0$ and the corresponding $AXA^*$  can be written as
can be written in the following parametric forms
\begin{align}
&X = UU^{*} + W - A^{\dag}AWA^{\dag}A,
\label{425}
\\
& AXA^* = -AUU^{*}A^*,
\label{425xx}
\end{align}
where $U\in \mathbb C^{n\times n}$ and  $W \in {\mathbb C}_{{\rm H}}^{n}$ are arbitrary$.$

\item[{\rm(b)}] There exists an $X \in  {\mathbb C}_{{\rm H}}^{n}$ such that $AXA^{*} \succ  0$
 if and only if  $r(A)= m.$ In this case$,$ the general Hermitian solution of $AXA^{*} \succ  0$ can be written
as {\rm (\ref{425})}$,$ in which $U\in \mathbb C^{n\times n}$ and $W\in  {\mathbb C}_{{\rm H}}^{n}$
are arbitrary$.$

\item[{\rm(c)}] The general Hermitian solution of  $AXA^{*} \preccurlyeq 0$ and the corresponding $AXA^*$  can be written
in the following parametric forms
\begin{align}
& X = -UU^{*} +  W - A^{\dag}AWA^{\dag}A,
\label{426}
\\
& AXA^* = -AUU^{*}A^*,
\label{426xx}
\end{align}
where $U\in {\mathbb C}^{n\times n}$ and  $W\in {\mathbb C}_{{\rm H}}^{n}$ are arbitrary$.$

\item[{\rm(d)}] There exists an $X \in {\mathbb C}_{{\rm H}}^{n}$
such that $AXA^{*} \prec 0$ if and only if $r(A)= m.$
In this case$,$ the general Hermitian solution of  $AXA^{*} \prec  0$ can be
written as {\rm (\ref{426})}$,$
 in which $U\in \mathbb C^{n\times n}$ is any matrix such that $r(AU) =m,$ and
 $W\in  {\mathbb C}_{{\rm H}}^{n}$ is arbitrary$.$
\end{enumerate}
\end{corollary}

\begin{theorem} \label{T47}
Let $A \in \mathbb C^{m\times n}$ and  $B \in \mathbb C^{m\times m}$ be given$.$
 Then$,$ the following hold$.$
\begin{enumerate}
\item[{\rm(a)}]  There exists an $X \in \mathbb C^{n\times n}$ such
that
\begin{equation}
AXX^{*}A^{*} \succcurlyeq BB^{*}
\label{427}
\end{equation}
if and only if  ${\mathscr R}(B) \subseteq {\mathscr R}(A).$ In this case$,$  a solution of {\rm (\ref{427})}  and the corresponding $AXX^*A^*$  can be written
in the following parametric forms
\begin{align}
& XX^{*}  = [\, A^{\dag}(\,BB^*  + AUU^{*}A^* \,)^{1/2} + F_AW \,][\, A^{\dag}(\, BB^*  + AUU^{*}A^*\,)^{1/2} + F_AW \,]^*,
\label{428}
\\
& AXX^{*}A^*  = BB^*  + AUU^{*}A^*,
\label{428xx}
\end{align}
where $U\in {\mathbb C}^{m\times m}$ and  $W\in {\mathbb C}^{n\times m}$ are arbitrary$.$

\item[{\rm(b)}]  There exists an $X \in \mathbb C^{n\times n}$ such
that
\begin{equation}
AXX^{*}A^{*} \succ  BB^{*}
 \label{429}
\end{equation}
if and only if  $r(A) =m.$ In this case$,$ a  solution of {\rm (\ref{429})}
can be written as
\begin{align}
& XX^{*}  = [\, A^{\dag}(\,BB^*  + UU^{*} \,)^{1/2} + F_AW \,][\, A^{\dag}(\, BB^*  + UU^{*}\,)^{1/2} + F_AW \,]^*,
\label{429jj}
\\
& AXX^{*}A^*  = BB^*  + UU^{*},
\label{429jjj}
\end{align}
where $U\in {\mathbb C}^{m\times m}$ is any matrix with $r(U) =m,$ and $W\in {\mathbb C}^{n\times m}$ is arbitrary$.$
\end{enumerate}
\end{theorem}

An application to partitioned matrices is given below.

\begin{corollary} \label{T48}
Let
\begin{equation}
\phi(X) =\left[\!\! \begin{array}{cc}  AXA^{*} & B
\\
B^{*} & CC^{*}\end{array} \!\!\right]\!,
\label{430}
\end{equation}
where $A \in \mathbb C^{m\times n}$, $B \in \mathbb C^{m\times p}$ and
 $C \in \mathbb C^{p\times p}$ are given$.$ Then$,$ the following hold$.$
\begin{enumerate}
\item[{\rm(a)}] There exists an $X \in {\mathbb C}_{{\rm H}}^{n}$
such that $\phi(X) \succcurlyeq 0$  if and only if
\begin{equation}
{\mathscr R}(B) \subseteq {\mathscr R}(A) \ \ and \ \  {\mathscr
R}(B^{*}) \subseteq {\mathscr R}(C).
 \label{431}
\end{equation}
In this case$,$ the general solution of $\phi(X) \succcurlyeq 0$ can be written in the following parametric form
\begin{equation}
X = A^{\dag}B(CC^{*})^{\dag}B^{*}(A^{\dag})^{*} + UU^{*} + W -
A^{\dag}AWA^{\dag}A,
 \label{432}
\end{equation}
where $U\in \mathbb C^{n\times n}$ and  $W\in {\mathbb C}_{{\rm H}}^{n}$  are arbitrary$.$

\item[{\rm(b)}]  There exists an $X \in \mathbb C^{n\times n}$ such that $\phi(X) \succ  0$ in {\rm (\ref{430})}
if and only if
\begin{equation}
r(A) = m \ \ and  \ \  r(C) =p.
\label{433}
\end{equation}
In this case$,$ the general solution of $\phi(X) \succ  0$ can be written in the following parametric form
\begin{equation}
X = A^{\dag}B(CC^{*})^{-1}B^{*}(A^{\dag})^{*} + UU^{*} + W -
A^{\dag}AWA^{\dag}A,
 \label{434}
\end{equation}
where $U\in \mathbb C^{n\times n}$ is any matrix such $r(AU) =m$ and $W\in {\mathbb C}_{{\rm H}}^{n}$  is arbitrary$.$
\end{enumerate}
\end{corollary}

\noindent {\bf Proof.} \ It  is easily seen from Lemma \ref{T23}(e) and (f) that
\begin{align}
& \phi(X) \succcurlyeq 0 \Leftrightarrow {\mathscr R}(B) \subseteq {\mathscr R}(A), \ \
 {\mathscr R}(B^{*}) \subseteq {\mathscr R}(C) \ \  {\rm and} \ \ AXA^{*} \succcurlyeq B(CC^{*})^+B^{*},
\label{435}
\\
& \phi(X) \succ  0 \Leftrightarrow r(A) =m, \ \  r(C) =p \ \  {\rm and} \ \
AXA^{*} \succ  B(CC^{*})^+B^{*}.
\label{436}
\end{align}
Solving the two inequalities in (\ref{435}) and (\ref{436}) by Theorem \ref{T41} leads to (a) and (b). \qquad $\Box$

\medskip

We next solve $AXX^{*}A^{*} \preccurlyeq BB^{*}$. It is obvious that the
inequality has a trivial solution $X =0$. However, the inequality
may only have zero solution in some cases. For example, the
inequality
$$
\left[\!\! \begin{array}{cc} x^2 & 0
\\
0 &  0 \end{array} \!\!\right] \preccurlyeq \left[\!\! \begin{array}{cc} 1 & 1
\\
1 & 1\end{array} \!\!\right]
$$
only has a solution $x=0$.

\begin{theorem} \label{T49}
Let $A \in \mathbb C^{m\times n}$ and $B \in \mathbb C^{m\times m}$ be given$.$ Then$,$ the following hold$.$
\begin{enumerate}
\item[{\rm(a)}] There exists an $X \in \mathbb C^{n\times n}$ such that both $AX \neq 0$ and
\begin{equation}
AXX^{*}A^{*} \preccurlyeq BB^{*}
\label{437}
\end{equation}
if and only if
\begin{equation}
 {\mathscr R}(A) \cap {\mathscr R}(B) \neq \{0\}.
\label{438}
\end{equation}
In this case$,$ a solution of {\rm (\ref{437})} and the corresponding $AXX^*A^*$ can be written in the following parametric forms
\begin{align}
& XX^{*} =  [\, A^{\dag}(BF_{B_1}VF_{B_1}B^{*})^{1/2} + F_AW \,] [\, A^{\dag}(BF_{B_1}VF_{B_1}B^{*})^{1/2} + F_AW \,]^*,
\label{439}
\\
& AXX^{*}A^{*} = BF_{B_1}VF_{B_1}B^{*},
\label{439xx}
\end{align}
where $B_1 = E_AB,$ $V$ is any matrix satisfying $0\prec  V \preccurlyeq I_m,$
 and $W \in \mathbb C^{n\times m}$ is arbitrary$.$  The rank of {\rm (\ref{439xx})} is
 \begin{equation}
\max_{AXX^{*}A^{*} \preccurlyeq BB^{*}} r(AXX^{*}A^*) = r(A) + r(B)- r[\, A, \ B \,].
\label{439a}
\end{equation}

\item[{\rm(b)}] There exists an $X \in \mathbb C^{n\times n}$ such that $AX \neq 0$ and
\begin{equation}
 AXX^{*}A^{*} \prec  BB^{*}
\label{440}
\end{equation}
if and only if
\begin{equation}
A \neq 0 \ \  and  \ \ r(B) = m.
\label{441}
\end{equation}
In this case$,$ a solution of {\rm (\ref{440})} can be written as {\rm (\ref{439})}$,$  
in
which $V$ is any matrix satisfying $0\prec  V \prec  I_m,$ and
$W \in \mathbb C^{n\times m}$ is arbitrary$.$

\item[{\rm(c)}] Under the condition ${\mathscr R}(B) \subseteq  {\mathscr R}(A),$
there always exists an $X \in \mathbb C^{n\times n}$ such that both
$AX \neq 0$ and
\begin{equation}
AXX^{*}A^{*} \preccurlyeq BB^{*},
 \label{442}
\end{equation}
 and a solution of {\rm (\ref{442})} and the corresponding $AXX^*A^*$ can be written in the following parametric forms
\begin{align}
& XX^{*} = [\, A^{\dag}(BVB^{*})^{1/2} + F_AW \,] [\, A^{\dag}(BVB^{*})^{1/2} + F_AW \,]^*,
\label{443}
\\
& AXX^{*}A^* = BVB^{*},
\label{443xx}
\end{align}
where $V$ is any matrix satisfying $0\prec  V \preccurlyeq I_m,$
 and $W \in \mathbb C^{n\times m}$ is arbitrary$.$

\item[{\rm(d)}] Under the condition ${\mathscr R}(B) \subseteq  {\mathscr R}(A),$
there exists an $X \in \mathbb C^{n\times n}$ such that $AX \neq 0$
and
\begin{equation}
 AXX^{*}A^{*} \prec  BB^{*}
\label{444}
\end{equation}
 if and only if $ r(B) = m.$ In this case$,$ a solution of {\rm (\ref{444})}
can be written as {\rm (\ref{443})}$,$ in which $V$ is any matrix
satisfying $0 \prec  V \prec  I_m,$ and $W \in \mathbb C^{n\times m}$ is arbitrary$.$
\end{enumerate}
\end{theorem}

\noindent {\bf Proof.} \
It can be seen from Lemma \ref{T23}(g) that if there exists an $X$ such that $AX \neq 0$ and
(\ref{437}) hold, then  ${\mathscr R}(AX) \subseteq {\mathscr R}(B)$, which obviously
implies that (\ref{438}) holds.  On  the other hand, it can be derived from $E_ABF_{E_AB} =0$ that
\begin{equation}
 AA^{\dag}BF_{E_AB} = BF_{E_AB},
\label{445}
\end{equation}
and from (\ref{215}) and (\ref{216})  that
\begin{equation}
 r(BF_{E_{A}B}) = r\!\left[\!\! \begin{array}{c} B
\\
E_{A}B \end{array} \!\!\right] - r(E_AB) = r(A) + r(B) - r[\, A, \, B\,]
= \dim [{\mathscr R}(A) \cap {\mathscr R}(B)].
\label{446}
\end{equation}
Hence if (\ref{438}) holds, then $BF_{E_{A}B} \neq 0$ and ${\mathscr R}(BF_{E_{A}B}) =
{\mathscr R}(A) \cap {\mathscr R}(B)$ by (\ref{445}) and (\ref{446}). In this case,
$$
AA^{\dag}BF_{E_AB}VF_{E_AB}B^{*}(A^{\dag})^{*}A =
BF_{E_AB}VF_{E_AB}B^{*}.
$$
Thus we can derive from (\ref{439}) and Lemma \ref{T21}(c) that
$$
BB^{*} - AXX^{*}A^{*} = BB^{*} - BF_{E_AB}VF_{E_AB}B^{*} = B(\,I_m -
F_{E_AB}VF_{E_AB} \,)B^{*} \succcurlyeq 0,
$$
that is, (\ref{439}) is  a solution of (\ref{437}).
The two conditions in (\ref{441}) are obvious under the condition that both $AX \neq 0$ and (\ref{440}) hold.
Conversely, if (\ref{441}) holds, we can derive from (\ref{440}) and Lemma \ref{21}(d) that
$I_m  - F_{E_AB}VF_{E_AB}\succ 0$ and
$$
BB^{*} - AXX^{*}A^{*} = BB^{*} - BF_{E_AB}VF_{E_AB}B^{*} = B(\,I_m -
F_{E_AB}VF_{E_AB} \,)B^{*} \succ 0.
$$
Results (c) and (d) are direct consequences of (a) and (b). \qquad $\Box$

\medskip

A direct consequence of Corollary \ref{T48jj} is given below.

\begin{corollary} \label{T511jj}
Let $A \in \mathbb C^{m\times n}$ and $B, \  C \in \mathbb C_{{\rm H}}^{m}$ be given$,$
and let ${\cal S}_1$ be of the form in
 {\rm (\ref{427yy})}$,$ and define
\begin{align}
& K_1 =\left[\!\begin{array}{cccccc} B & B & A
\\
 B & B-C & 0
\\
A^* & 0 & 0
\end{array}\!\right]\!,   \ \  K_2 =\!\left[\!\! \begin{array}{cccc} C  & A
 \\ A^{*} & 0 \end{array} \!\!\right]\!.
\label{581}
\end{align}
Then$,$ the extremal ranks and partial inertias of $AXA^* - C$ subject
  to $X \in {\cal S}_1$ are given by
\begin{align}
&\max_{X \in {\cal S}_1} r(\, AXA^* - C \,)  = r[\, A, \, C\,],
\label{582}
\\
&\min_{X \in {\cal S}_1} r(\, AXA^* - C \,)  =  i_{+}(K_1) -  r(K_2) + r[\, A, \, C\,],
\label{583}
\\
&\max_{X \in {\cal S}_1} i_{+}(\,AXA^* - C \,)  =i_{-}(K_2),
\label{584}
\\
&\max_{X \in {\cal S}_1} i_{-}(\, AXA^* - C \,)  =i_{-}(K_1) - i_{-}(K_2),
\label{585}
\\
&\min_{X \in {\cal S}_1} i_{+}(\, AXA^* - C\,)  =i_{+}(K_1) - i_{+}(K_2),
\label{586}
\\
&\min_{X \in {\cal S}_1} i_{-}(\, AXA^* - C\,)  = r[\, A, \, C\,] -  i_{-}(K_2).
\label{587}
\end{align}
In consequence$,$ the following hold$.$
\begin{enumerate}
\item[{\rm(a)}] There exists an $X \in {\mathbb C}_{{\rm H}}^{n}$ such that
$ AXA^* \succcurlyeq B$  and $ AXA^* \succ  C$ if and only if $i_{-}(K_2) = m.$

\item[{\rm(b)}] There exists an $X \in {\mathbb C}_{{\rm H}}^{n}$ such that
$ C \succ  AXA^* \succcurlyeq B$ if and only if $C\succ  B$ and $i_{-}(K_1) = i_{-}(K_2) + m.$

\item[{\rm(c)}] There exists an $X \in {\mathbb C}_{{\rm H}}^{n}$ such that
$ C \succcurlyeq AXA^* \succcurlyeq B$ if and only if $D\succcurlyeq C$ and  $i_{+}(K_1) = i_{+}(K_2).$

\item[{\rm(d)}] There exists an $X \in {\mathbb C}_{{\rm H}}^{n}$ such that
$ AXA^* \succcurlyeq B$ and $ AXA^* \succcurlyeq C$ f and only if
 $i_{-}(K_2) =  r[\, A, \, C\,].$

\end{enumerate}

\end{corollary}

\section{General solution of $AX +(AX)^{*}  \succcurlyeq  \,(\succ, \, \preccurlyeq, \, \prec ) \, B$ and its properties}
\renewcommand{\theequation}{\thesection.\arabic{equation}}
\setcounter{section}{5}

\setcounter{equation}{0}

The inequality in (\ref{10d}) was approached in \cite{T-laa11} by
using a relation method and their general solutions were given analytically.
In this section, we reconsider this inequality and give some new conclusions
on algebraic properties of its solution.

\begin{theorem} \label{T31}
Let $A \in \mathbb C^{m\times n}$ and  $B \in \mathbb C_{{\rm H}}^{m}$ be given$,$ and let
$M = \left[\!\! \begin{array}{cc} B & A \\ A^* & 0 \end{array} \!\!\right]\!.$
Then$,$ the following hold$.$
\begin{enumerate}
\item[{\rm(a)}] {\rm \cite{T-laa11}} The following statements are equivalent$:$
\begin{enumerate}
\item[{\rm(i)}]  There exists an $X \in \mathbb C^{n\times m}$ such that
\begin{equation}
AX + (AX)^*  \succcurlyeq B.
\label{32}
\end{equation}

\item[{\rm(ii)}] $E_ABE_A \preccurlyeq 0.$

 \item[{\rm(iii)}] $i_{+}(M)= r(A).$
\end{enumerate}
In this case$,$ the general solution of {\rm (\ref{32})} and the corresponding $AX + (AX)^{*}$
can be written in the following parametric forms
\begin{align}
& X = \frac{1}{2}A^{\dag}B\widehat{A} +
\frac{1}{2}A^{\dag}(\, AU + J^{\frac{1}{2}} \,)(\, AU + J^{\frac{1}{2}} \,)^*\widehat{A} + VA^* + F_AW,
\label{34}
\\
& AX + (AX)^{*} = B + (\, AU + J^{\frac{1}{2}} \,)(\, AU + J^{\frac{1}{2}} \,)^*,
\label{35j}
\end{align}
where  $J= -E_ABE_A,$  $\widehat{A} = 2I_m -  AA^{\dag},$ and
$U, \ W\in \mathbb C^{n \times m}$ and $V\in \mathbb C_{{\rm SH}}^{n}$ are arbitrary$.$

\item[{\rm(b)}] {\rm \cite{T-laa11}} There exists an $X \in \mathbb C^{n\times m}$ such that
\begin{equation}
AX + (AX)^* \succ  B
\label{315j}
\end{equation}
if and only if
\begin{equation}
E_ABE_A \preccurlyeq 0 \ \  and  \ \ \R(E_ABE_A) =\R(E_A),
\label{316j}
\end{equation}
 or equivalently$,$ $i_{-}(M) = m.$ In this case$,$ the general solution of {\rm (\ref{315j})} can be
written as {\rm (\ref{34})}$,$ in which $U$ is any matrix such that
$r(\, AU + J^{\frac{1}{2}}  \,) = m,$  say$,$  $U = A^*,$  $V\in \mathbb C_{{\rm SH}}^{n}$ and
$W\in \mathbb C^{n \times m}$ are arbitrary$.$

\item[{\rm(c)}] Under {\rm (a)}$,$ let
\begin{equation}
{\cal S}_1 =\{\,  X \in {\mathbb C}^{n\times m} \ | \  AX + (AX)^{*} \succcurlyeq  B \, \}.
\label{36j}
\end{equation}
 Then$,$ the  extremal ranks and partial inertias of $AX + (AX)^{*}$ and  $AX + (AX)^{*} - B$  subject
  to $X \in {\cal S}_1$ are given by
\begin{align}
\max_{X \in {\cal S}_1} r[\, AX + (AX)^{*} \,] & =  \min\{\, 2 r(A), \ \  r[\,A, \,B\,]\, \},
\label{37j}
\\
 \min_{X \in {\cal S}_1} r[\, AX + (AX)^{*} \,] & = \max\{ \,2r(A) + 2r[\,A, \,B\,] - 2r(M), \ \
 r(B) - m,  \nb
 \\
& \ \ \ \ \ \ \ \ \ \ \ \ i_{+}(B) + r(A) + r[\, A, \,B\,] - r(M), \nb
\\
& \ \ \ \ \ \ \ \ \ \ \ \ i_{-}(B) + r(A) + r[\, A, \,B\,] - r(M) - m \, \},
\label{38j}
\\
\max_{X \in {\cal S}_1} i_{+}[\, AX + (AX)^{*} \,] & = r(A),
\label{39j}
\\
\max_{X \in {\cal S}_1} i_{-}[\, AX + (AX)^{*} \,] & = \min\left\{\, r(A), \ \   i_{-}(B) \, \right\},
\label{310j}
\\
 \min_{X \in {\cal S}_1} i_{+}[\, AX + (AX)^{*} \,]  & =  \max \left\{ r[\,A, \, B\,] + r(A) - r(M),  \
 \ i_{+}(B) \right\},
\label{311j}
\\
\min_{X \in {\cal S}_1} i_{-}[\, AX + (AX)^{*} \,]  & = r[\,A, \, B\,] + r(A) - r(M),
\label{312j}
\\
\max_{X \in {\cal S}_1} r[\, AX + (AX)^{*} - B  \,] & = r(M) - r(A),
\label{313j}
\\
 \min_{X \in {\cal S}_1} r[\, AX + (AX)^{*} - B  \,] & = r(M) - 2r(A).
\label{314j}
\end{align}
In consequence$,$

\item[{\rm(d)}] There exists an $X \in \mathbb C^{n\times m}$ such that
$AX + (AX)^*  \succ  0$ and $AX + (AX)^*  \succcurlyeq B$  if and only if $r(A) = m.$

\item[{\rm(e)}] There exists an $X \in \mathbb C^{n\times m}$ such that
$0 \succ  AX + (AX)^*  \succcurlyeq B$ if and only if $r(A) = m$ and $ B \prec  0.$

\item[{\rm(f)}] There exists an $X \in \mathbb C^{n\times m}$ such that
$AX + (AX)^* \succcurlyeq 0$ and $AX + (AX)^*  \succcurlyeq B$  if and only if $r(N) = r[\,A, \, B\,] + r(A).$

\item[{\rm(g)}] There exists an $X \in \mathbb C^{n\times m}$ such that
$ 0 \succcurlyeq AX + (AX)^* \succcurlyeq B$  if and only if $ B \preccurlyeq 0.$

\end{enumerate}
\end{theorem}

\noindent {\bf Proof.} \ Inequality (\ref{32}) can be relaxed to the following quadratic matrix equation
\begin{equation}
AX + (AX)^{*} = B + YY^{*},
\label{332j}
\end{equation}
where $Y \in {\mathbb C}^{m \times m}$.  From Lemma \ref{T24}(a), there exists an $X$
that satisfies (\ref{332j}) if and only if $YY^{*}$ satisfies $E_A(\,B + YY^{*}\,)E_A =0$, that is,
\begin{equation}
E_AYY^{*}E_A = -E_ABE_A = J.
\label{333j}
\end{equation}
Further by Lemma \ref{T26}(c), there exists a $YY^{*}$ that satisfies (\ref{333j})
if and only if (ii) of (a) holds, in which case, the general solution of (\ref{333j})
can be written as
\begin{align}
YY^{*} = (\, AU + J^{\frac{1}{2}}\,)(\, AU + J^{\frac{1}{2}} \,)^*,
\label{334j}
\end{align}
where  $U\in {\mathbb C}^{n \times m}$ is arbitrary. Substituting this $YY^{*}$
into (\ref{332j}) gives
\begin{equation}
AX + (AX)^{*} = B + (\, AU + J^{\frac{1}{2}} \,)(\, AU + J^{\frac{1}{2}} \,)^*.
\label{335j}
\end{equation}
Applying Lemma \ref{T24}(a) to this equation, we obtain (\ref{34}).

Setting (\ref{313j}) equal to $m$ gives $r(M) - r(A) =m$, i.e.,
$r(E_ABE_A) =r(E_A)$ by (\ref{21x}), which is further equivalent to (\ref{316j}).
The equivalence of $i_{-}(M) =m$ and (\ref{316j}) follows
from (\ref{21x}) and $i_{-}(E_ABE_A) \leqslant  r(E_ABE_A) \leqslant r(E_A)$.

Applying Lemma \ref{T212F}(a) to (\ref{335j}), we obtain
\begin{align}
&\max_{U\in {\mathbb C}^{n \times m}}\!\!\!r[\, B + (\, AU + J^{\frac{1}{2}} \,)(\, AU + J^{\frac{1}{2}} \,)^*\,]
= \min\left\{r[\,  A, \,  B, \, J^{\frac{1}{2}} \,], \  r\!\left[\!\! \begin{array}{cc}  B + J & A
 \\ A^* & 0 \end{array} \!\!\right]\!, \  r(B) + m \right\},
\label{336j}
\\
&\min_{U\in {\mathbb C}^{n \times m}}\!\!\!r[\,B + (\, AU + J^{\frac{1}{2}} \,)(\, AU + J^{\frac{1}{2}} \,)^*\,] =
2r[\,A, \,B, \, J^{\frac{1}{2}}\,]  + \max\{ \, h_1, \ \ h_2, \ \ h_3, \ \ h_4 \, \},
\label{337j}
\\
& \max_{U\in {\mathbb C}^{n \times m}}\!\!\!i_{+}[\, B + (\, AU + J^{\frac{1}{2}} \,)(\, AU + J^{\frac{1}{2}} \,)^*\,]
 = \min\left\{\, i_{+} \!\left[\!\! \begin{array}{cc}  B + J & A
 \\ A^* & 0 \end{array} \!\!\right]\!, \ \   i_{+}(B)  + m\, \right\},
\label{338j}
\\
& \max_{U\in {\mathbb C}^{n \times m}}\!\!\!i_{-}[\, B + (\, AU + J^{\frac{1}{2}} \,)(\, AU + J^{\frac{1}{2}} \,)^*\,]
 = \min\left\{\, i_{-} \!\left[\!\! \begin{array}{cc}  B + J & A
 \\ A^* & 0 \end{array} \!\!\right]\!, \ \   i_{-}(B) \, \right\},
\label{339j}
\\
&  \min_{U\in {\mathbb C}^{n \times m}} \!\!\!i_{+} [\, B + (\, AU + J^{\frac{1}{2}} \,)
(\, AU + J^{\frac{1}{2}} \,)^*\,]  \nb
\\
& = r[\,A, \,B, \, J^{\frac{1}{2}}\,]  + \max \left\{ i_{+}\!\left[\!\! \begin{array}{cc}  B + J & A
 \\ A^* & 0 \end{array} \!\!\right] - r\!\left[\!\!\begin{array}{cccc}  B & A & J^{\frac{1}{2}}
 \\ A^* & 0 & 0  \end{array} \!\!\right],  \  \ i_{+}(B)  - r[\,  A, \, B \,] \right\},
\label{340j}
\\
&  \min_{U\in {\mathbb C}^{n \times m}}\!\!\! i_{-} [\, B + (\, AU + J^{\frac{1}{2}} \,)(\, AU + J^{\frac{1}{2}} \,)^*\,]  \nb
\\
&  = r[\,A, \,B, \, J^{\frac{1}{2}}\,]  + \max \left\{ i_{-}\!\left[\!\! \begin{array}{cc}  B + J & A
 \\ A^* & 0 \end{array} \!\!\right] - r\!\left[\!\!\begin{array}{cccc}  B & A & J^{\frac{1}{2}}
 \\ A^* & 0 & 0  \end{array} \!\!\right],  \ \ i_{-}(B)  - r[\,  A, \, B \,] - m \right\},
\label{341j}
\end{align}
where
\begin{align*}
 h_1 & = r\!\left[\!\! \begin{array}{cc}  B + J & A
 \\ A^* & 0 \end{array} \!\!\right] - 2r\!\left[\!\! \begin{array}{ccc}  B & A & J^{\frac{1}{2}}
 \\ A^* & 0 & 0 \end{array} \!\!\right],
 \\
h_2 & =  r(B) -  2r[\,A, \,  B\,] - m,
\\
h_3 & =  i_{-}\!\left[\!\! \begin{array}{cc}  B + J & A
 \\ A^* & 0 \end{array} \!\!\right] - r\!\left[\!\! \begin{array}{ccc}  B & A & J^{\frac{1}{2}}
 \\ A^* & 0  & 0 \end{array} \!\!\right] + i_{+}(B) - r[\,A, \,  B\,],
\\
h_4 &=  i_{+}\!\left[\!\! \begin{array}{cc}  B + J & A
 \\ A^* & 0 \end{array} \!\!\right] - r\!\left[\!\! \begin{array}{ccc}  B & A & J^{\frac{1}{2}}
 \\ A^* & 0  & 0 \end{array} \!\!\right] + i_{-}(B) - r[\,A, \,  B\,] -m.
\end{align*}
Simplifying the ranks and partial inertias of the block matrices in  (\ref{336j})--(\ref{341j}) gives
\begin{align}
r[\,  A, \,  B, \, J^{\frac{1}{2}} \,] & = r[\,  A, \,  B, \, J \,] =
r[\,  A, \,  B, \, E_ABE_A \,] = r[\,  A, \,  B \,],
\label{342j}
\\
i_{\pm}\!\left[\!\! \begin{array}{cc}  B + J & A
 \\ A^* & 0 \end{array} \!\!\right]&  = i_{\pm}\!\left[\!\! \begin{array}{cc}  B - E_ABE_A  & A
 \\ A^* & 0 \end{array} \!\!\right] = i_{\pm}\!\left[\!\! \begin{array}{cc}  0  & A
 \\ A^* & 0 \end{array} \!\!\right]  = r(A),
\label{343j}
\\
 r\!\left[\!\! \begin{array}{ccc}  B & A & J^{\frac{1}{2}}
 \\ A^* & 0  & 0 \end{array} \!\!\right] & = r\!\left[\!\! \begin{array}{ccc}  B & A & J
 \\ A^* & 0  & 0 \end{array} \!\!\right] = r\!\left[\!\! \begin{array}{ccc}  B & A & E_ABE_A
 \\ A^* & 0  & 0 \end{array} \!\!\right] \nb
 \\
&  = r\!\left[\!\! \begin{array}{ccc}  B & A & BE_A
 \\ A^* & 0  & 0 \end{array} \!\!\right] = r\!\left[\!\! \begin{array}{ccc}  B & A & B
 \\ A^* & 0  & 0 \\ 0 & 0 & A^* \end{array} \!\!\right] - r(A) \nb
 \\
 &  = r\!\left[\!\! \begin{array}{ccc}  B & A & 0
 \\ A^* & 0  & 0 \\ 0 & 0 & A^* \end{array} \!\!\right] - r(A) = r\!\left[\!\! \begin{array}{ccc}  B & A
 \\ A^* & 0 \end{array} \!\!\right].
 \label{344j}
 \end{align}
 Substituting (\ref{342j})--(\ref{344j}) into (\ref{336j})--(\ref{341j}) gives (\ref{37j})--(\ref{312j}).
 It can be seen from (\ref{335j}) that
\begin{equation}
r[\, AX + (AX)^{*} - B \,] = r(\, AU + J^{\frac{1}{2}} \,).
\label{345j}
\end{equation}
Hence, we  derive from (\ref{12a}) and (\ref{12b})  that
\begin{align*}
\max r[\, AX + (AX)^{*} - B \,]&  = \max_{U} r(\, AU + J^{\frac{1}{2}} \,)  =
r[\, A,\, J^{\frac{1}{2}}\,] \nb
 \\
 & = r[\,A, \, E_ABE_A\,] = r(A) + r(E_ABE_A)  \ \ \mbox{(by (\ref{215}))}
\\
\min r[\, AX + (AX)^{*} - B \,]  & = \min_{U} r(\, AU + J^{\frac{1}{2}} \,)  =
r[\, A, \, J^{\frac{1}{2}} \,] - r(A) \nb
\\
&  =  r[\, A, \,  E_ABE_A \,] - r(A)   = r(E_ABE_A) \ \ \mbox{(by (\ref{215}))},
\end{align*}
establishing (\ref{313j}) and (\ref{314j}). \qquad $\Box$

\medskip

 The following results can be shown similarly.

\begin{theorem} \label{T31a}
Let $A \in \mathbb C^{m\times n}$ and  $B, \ C \in \mathbb C_{{\rm
H}}^{m}$ be given$,$ and let
$M = \left[\!\! \begin{array}{cc} B & A \\ A^* & 0 \end{array} \!\!\right]\!.$  Then$,$ the following hold$.$
\begin{enumerate}
\item[{\rm(a)}] {\rm \cite{T-laa11}} The following statements are equivalent$:$
\begin{enumerate}
\item[{\rm(i)}]  There exists an $X \in \mathbb C^{n\times m}$ such that
\begin{equation}
AX + (AX)^*  \preccurlyeq B.
\label{317j}
\end{equation}

\item[{\rm(ii)}] $E_ABE_A \succcurlyeq 0.$

 \item[{\rm(iii)}]   $i_{-}(M)= r(A).$
\end{enumerate}
In this case$,$ the general solution of {\rm (\ref{317j})}
and the corresponding $AX + (AX)^{*}$ can be written in the following parametric forms
\begin{align}
& X  =  \frac{1}{2}A^{\dag}B\widehat{A} - \frac{1}{2}A^{\dag}(\, AU + K^{\frac{1}{2}} \,)
(\, AU + K^{\frac{1}{2}} \,)^*\widehat{A} + VA^* + F_AW,
\label{319j}
\\
& AX + (AX)^{*}  = B - (\, AU + K^{\frac{1}{2}} \,)(\, AU + K^{\frac{1}{2}} \,)^*,
\label{320j}
\end{align}
where  $K= E_ABE_A,$  $\widehat{A} = 2I_m -  AA^{\dag},$  $U, \, W\in \mathbb C^{n \times m}$ and
$V\in \mathbb C_{{\rm SH}}^{n}$ are arbitrary$.$

\item[{\rm(b)}]  {\rm \cite{T-laa11}} There exists an $X \in \mathbb C^{n\times m}$ such that
\begin{equation}
AX + (AX)^* \prec  B
\label{330j}
\end{equation}
if and only if
\begin{equation}
E_ABE_A \succcurlyeq 0 \ \  and  \ \ \R(E_ABE_A) = \R(E_A),
\label{331j}
\end{equation}
 or equivalently$,$ $i_{+}(M) = m.$ In this case$,$ the general solution of {\rm (\ref{330j})}
can be written as {\rm (\ref{319j})}$,$ in which $U$ is any matrix
such that $r(\, AU + K^{\frac{1}{2}}  \,) = m,$  say$,$  $U = A^*,$  $V\in \mathbb C_{{\rm SH}}^{n}$ and
$W\in \mathbb C^{n \times m}$ are arbitrary$.$

\item[{\rm(c)}]  Under {\rm (a)}$,$ let
\begin{equation}
{\cal S}_2 =\{\,  X \in {\mathbb C}^{n\times m} \ | \  AX + (AX)^{*} \preccurlyeq  B \, \}.
\label{321j}
\end{equation}
 Then$,$ the extremal ranks and partial inertias of $AX + (AX)^{*}$ and $AX + (AX)^{*} - B$  subject
  to $X \in {\cal S}_2$ are given by
 \begin{align}
& \max_{X \in {\cal S}_2} r[\, AX + (AX)^{*} \,]  =  \min\{\, 2 r(A), \ \  r[\,A, \,B\,]\, \},
\label{322j}
\\
&  \min_{X \in {\cal S}_2} r[\, AX + (AX)^{*} \,]  = \max\{ \,2r(A) + 2r[\,A, \,B\,] - 2r(N), \ \
 r(B) - m,  \nb
 \\
&  \hspace{4.5cm}  i_{+}(B) + r(A) + r[\, A, \,B\,] - r(N) -m, \ \ i_{-}(B) + r(A) + r[\, A, \,B\,] - r(N)\, \},
\label{323j}
\\
& \max_{X \in {\cal S}_2} i_{+}[\, AX + (AX)^{*} \,]  = \min\left\{\, r(A), \ \   i_{+}(B) \, \right\},
\label{324j}
\\
& \max_{X \in {\cal S}_2} i_{-}[\, AX + (AX)^{*} \,]  = r(A),
\label{325j}
\\
& \min_{X \in {\cal S}_2} i_{+}[\, AX + (AX)^{*} \,]  = r[\,A, \, B\,] + r(A) - r(N),
\label{326j}
\\
& \min_{X \in {\cal S}_2} i_{-}[\, AX + (AX)^{*} \,]  =  \max \left\{ r[\,A, \, B\,] + r(A) - r(N),  \
 \ i_{-}(B) \right\},
\label{327j}
\\
& \max_{X \in {\cal S}_2} r[\, AX + (AX)^{*} - B  \,]  = r(N) - r(A),
\label{328j}
\\
&  \min_{X \in {\cal S}_2} r[\, AX + (AX)^{*} - B  \,] = r(N) - 2r(A).
\label{329j}
\end{align}
In consequence$,$

\item[{\rm(d)}] There exists an $X \in \mathbb C^{n\times m}$ such that
$0 \prec  AX + (AX)^*  \preccurlyeq B$ if and only if $r(A) = m$ and $B \succ 0.$

\item[{\rm(e)}] There exists an $X \in \mathbb C^{n\times m}$ such that
$AX + (AX)^*  \prec  0$ and $AX + (AX)^*  \preccurlyeq B$  if and only if $r(A) = m.$

\item[{\rm(f)}] There exists an $X \in \mathbb C^{n\times m}$ such that
$AX + (AX)^* \preccurlyeq 0$ and $AX + (AX)^*  \preccurlyeq B$  if and only if $r(N) = r[\,A, \, B\,] + r(A).$

\item[{\rm(g)}] There exists an $X \in \mathbb C^{n\times m}$ such that
$ 0 \preccurlyeq AX + (AX)^* \preccurlyeq B$  if and only if $ B \succcurlyeq 0.$
\end{enumerate}
\end{theorem}

Theorem \ref{T31} established identifying conditions for the LMI in (\ref{10b}) to be solvable,
 and gave general expression of the matrix $X$ satisfying (\ref{10b}). In particular, the
general solutions in (\ref{34}) and (\ref{319j}) are represented in closed-forms by using
generalized inverse of the given matrices and arbitrary matrices. Hence, they can be directly used to deal
with various problems on the inequality in (\ref{10b}) and its properties. In what follows,
we present some consequences of Theorem \ref{T31} when $A$ and $B$ satisfy some more conditions.

\begin{corollary} \label{T32}
Let $A \in \mathbb C^{m\times n}$ and  $B \in \mathbb C_{{\rm H}}^{m}$ be given$,$ and assume
that there exists an $X \in \mathbb C^{n\times m}$ such that $AX + (AX)^* = B.$ Then$,$ the following hold$.$
\begin{enumerate}
\item[{\rm(a)}] The general solution $X \in \mathbb C^{n\times m}$ of
\begin{equation}
AX + (AX)^* \succcurlyeq B
\label{322gg}
\end{equation}
and the corresponding $AX + (AX)^{*}$ can be written in the following parametric forms
\begin{align}
& X  = \frac{1}{2}A^{\dag}B(\,2I_m -  AA^{\dag}\,) + UU^*A^* + VA^* + F_AW,
\label{323gg}
\\
& AX + (AX)^{*} = B + 2AUU^*A^*,
\label{324gg}
\end{align}
where $U\in \mathbb C^{n \times n},$ $W\in \mathbb C^{n \times m}$ and
$V\in \mathbb C_{{\rm SH}}^{n}$ are arbitrary$.$

\item[{\rm(b)}] There exists an $X \in \mathbb C^{n\times m}$ such that
\begin{equation}
AX + (AX)^* \succ  B
\label{325gg}
\end{equation}
if and only if $r(A) = m.$ In this case$,$ the general solution  can be
written as {\rm (\ref{323gg})}$,$ in which $U$ is any matrix such that
$r(AU) = m,$ and $W\in \mathbb C^{n \times m}$ and $V\in \mathbb C_{{\rm SH}}^{n}$ are arbitrary$.$

\item[{\rm(c)}]  The general solution $X \in \mathbb C^{n\times m}$ of
\begin{equation}
AX + (AX)^*  \preccurlyeq B
\label{326gg}
\end{equation}
and the corresponding $AX + (AX)^{*}$ can be written in the following parametric forms
\begin{align}
&X  = \frac{1}{2}A^{\dag}B(\,2I_m -  AA^{\dag}\,) - UU^*A^* + VA^* + F_AW,
\label{327gg}
\\
& AX + (AX)^{*} = B - 2AUU^*A^*,
\label{328gg}
\end{align}
where $U\in \mathbb C^{n \times n},$ $W\in \mathbb C^{n \times m}$ and
$V\in \mathbb C_{{\rm SH}}^{n}$ are arbitrary$.$

\item[{\rm(d)}] There exists an $X \in \mathbb C^{n\times m}$ such that
\begin{equation}
AX + (AX)^* \prec  B
\label{329gg}
\end{equation}
if and only if $r(A) = m.$ In this case$,$ the general solution of {\rm (\ref{329gg})} can be
written as {\rm (\ref{327gg})}$,$ in which $U$ is any matrix such that
$r(AU) = m,$ and $W\in {\mathbb C}^{n \times m}$ and $V\in \mathbb C_{{\rm SH}}^{n}$ are arbitrary$.$
\end{enumerate}
\end{corollary}

\begin{corollary} \label{T33}
Let $A \in \mathbb C^{m\times n}$ and  $B \in \mathbb C^{m\times k}$ be given$.$
Then$,$ the following hold$.$
\begin{enumerate}
\item[{\rm(a)}] There exists an $X \in \mathbb C^{n\times m}$ such that
\begin{equation}
AX + (AX)^{*} \succcurlyeq BB^{*}
\label{321}
\end{equation}
if and only if  ${\mathscr R}(B) \subseteq {\mathscr R}(A).$ In this case$,$
the general solution and the corresponding $AX + (AX)^{*}$ can be written as
\begin{align}
& X = \frac{1}{2}A^{\dag}BB^{*} +  UU^{*}A^{*} + VA^{*} + F_AW,
\label{322}
\\
& AX + (AX)^{*} = BB^* + 2AUU^*A^*,
\label{35pp}
\end{align}
where $U \in \mathbb C^{n \times n},$ $V \in \mathbb C_{{\rm
SH}}^{n}$ and $W\in \mathbb C^{n \times m}$ are arbitrary$.$

\item[{\rm(b)}] There exists an $X \in \mathbb C^{n\times m}$ such that
\begin{equation}
AX + (AX)^{*} \succ  BB^{*}
\label{323}
\end{equation}
if and only if  both ${\mathscr R}(B) \subseteq {\mathscr R}(A)$ and $r(A) = m.$ In this case$,$
  the general solution can be written as {\rm (\ref{322})}$,$ in which
$U$ is any matrix with $r(AU) =m,$ and $V \in \mathbb C_{{\rm
SH}}^{n}$ and $W\in \mathbb C^{n \times m}$ are arbitrary$.$

\item[{\rm(c)}] There exists an $X \in \mathbb C^{n\times m}$ such that
\begin{equation}
AX + (AX)^{*}  \preccurlyeq -BB^*
\label{324}
\end{equation}
if and only if  ${\mathscr R}(B) \subseteq {\mathscr R}(A).$ In this case$,$
the general solution  and the corresponding $AX + (AX)^{*}$  can be written as
\begin{align}
& X = -\frac{1}{2}A^{\dag}BB^{*} -  UU^{*}A^{*} + VA^{*} + F_AW,
\label{325}
\\
& AX + (AX)^{*} =  - BB^* - 2AUU^*A^*,
\label{322a}
\end{align}
where $U \in \mathbb C^{n \times n},$ $V\in \mathbb C_{{\rm
SH}}^{n}$ and $W\in \mathbb C^{n \times m}$ are arbitrary$.$

\item[{\rm(d)}] There exists an $X \in \mathbb C^{n\times m}$ such that
\begin{equation}
AX + (AX)^{*}  \prec -  BB^{*}
 \label{326}
\end{equation}
if and only if  both ${\mathscr R}(B) \subseteq {\mathscr R}(A)$ and
$r(A) = m.$ In this case$,$ the general solution can be written as {\rm (\ref{325})}$,$ in which $U$ is
any matrix with $r(AU) =m,$ and $V \in \mathbb C_{{\rm SH}}^{n}$ and
$W\in \mathbb C^{n \times m}$ are arbitrary$.$

\end{enumerate}
\end{corollary}

\begin{corollary} \label{T34}
Let $B \in \mathbb C^{n\times n}$ be given$.$ Then$,$ the following hold$.$
\begin{enumerate}
\item[{\rm(a)}]
The general solution $X\in \mathbb C^{n\times n}$ of
\begin{equation}
X + X^{*} \succcurlyeq BB^{*}
\label{328}
\end{equation}
and the corresponding $X + X^*$ can be written as
\begin{align}
& X = \frac{1}{2}BB^{*} + UU^{*} + V - V^*,
 \label{329}
 \\
& X + X^*  = BB^{*} + 2UU^{*},
\label{329xx}
\end{align}
where $U, \ V \in \mathbb C^{n \times n}$ are arbitrary$.$

\item[{\rm(b)}]
The general solution $X\in \mathbb C^{n\times n}$ of
\begin{equation}
X + X^{*} \succ  BB^{*}
\label{330}
\end{equation}
can be written as {\rm (\ref{329})}$,$  in which $U, \ V \in \mathbb C^{n
\times n}$ are arbitrary with $r(U) = n.$

\item[{\rm(c)}]
The general solution $X\in \mathbb C^{n\times n}$ of
\begin{equation}
X + X^{*} \preccurlyeq BB^{*}
\label{331}
\end{equation}
and the corresponding $X + X^*$ can be written as
\begin{align}
& X = \frac{1}{2}BB^{*} - UU^{*} + V - V^*,
\label{332}
\\
& X + X^*  = BB^{*} - 2UU^{*},
\label{332xx}
\end{align}
where $U, \ V \in {\mathbb C}^{n \times n}$ are
arbitrary$.$

\item[{\rm(d)}]
The general solution $X\in \mathbb C^{n\times n}$ of the inequality
\begin{equation}
X + X^{*} \prec  BB^{*}
 \label{333}
\end{equation}
can be written as {\rm (\ref{332})}$,$ in which $U, \ V \in \mathbb C^{n
\times n}$ are arbitrary with $r(U) = n.$
\end{enumerate}
\end{corollary}

\begin{corollary} \label{T35}
Let $A \in \mathbb C^{m\times n}$ be given$.$ Then$,$ the following hold$.$
\begin{enumerate}
\item[{\rm(a)}] The general solution $X\in \mathbb C^{n\times m}$ of
\begin{equation}
AX + (AX)^{*} \succcurlyeq 0
\label{334}
\end{equation}
and the corresponding $AX + (AX)^*$ can be written as
\begin{align}
& X = UU^{*}A^{*} + VA^{*} + F_AW,
\label{335}
\\
& AX + (AX)^{*}  =  2AUU^{*}A^{*}
\label{335xx}
\end{align}
where $U \in \mathbb C^{n \times n},$ $V \in \mathbb C_{{\rm
SH}}^{n}$ and $W\in \mathbb C^{n \times m}$ are arbitrary$.$

\item[{\rm(b)}] There exists an $X \in \mathbb C^{n\times m}$ such that
\begin{equation}
AX + (AX)^{*} \succ  0
\label{336}
\end{equation}
if and only if $r(A) =m.$ In this case$,$ the general solution can be written as
in {\rm (\ref{335})}$,$  in which
$U \in \mathbb C^{n \times n}$ is any matrix with $r(AU) =m,$
$V\in \mathbb C_{{\rm SH}}^{n}$ and $W\in \mathbb C^{n \times m}$ are
arbitrary$.$

\item[{\rm(c)}] The general solution $X\in \mathbb C^{n\times m}$ of
\begin{equation}
AX + (AX)^{*} \preccurlyeq 0
 \label{337}
\end{equation}
and the corresponding $AX + (AX)^*$ can be written as
\begin{align}
& X = -UU^{*}A^{*} + VA^{*} + F_AW,
\label{338}
\\
& AX + (AX)^{*}  =  - 2AUU^{*}A^{*},
\label{338xx}
\end{align}
where $U \in {\mathbb C}^{n \times n},$ $V \in \mathbb C_{{\rm
SH}}^{n}$ and $W\in \mathbb C^{n \times m}$ are arbitrary$.$

\item[{\rm(d)}] There exists an $X \in \mathbb C^{n\times m}$ such that
\begin{equation}
AX + (AX)^{*} \prec  0
\label{340}
\end{equation}
if and only if $r(A) =m.$ In this case$,$ the general solution
can be written as {\rm (\ref{338})}$,$  in which
$U \in \mathbb C^{n \times n}$ is any matrix with $r(AU) =m,$ and $V
\in \mathbb C_{{\rm SH}}^{n}$ and $W\in \mathbb C^{n \times m}$ are
arbitrary$.$ In particular$,$ if $A$ is square and nonsingular$,$
then the general solution of {\rm (\ref{340})} can be written as
\begin{equation}
X = -UU^{*}A^{*} + VA^{*},
\label{341}
\end{equation}
where $U \in \mathbb C^{n \times n}$ is any matrix with $r(AU) =m,$
and $V\in \mathbb C_{{\rm SH}}^{n}$ is arbitrary$.$
\end{enumerate}
\end{corollary}

As an application  of Theorems \ref{T31} and \ref{T31a}, we next give solutions of
the inequality $(A+B)X + X^{*}(A+B)^{*} \succcurlyeq AB +BA$,
which was considered for $A\succcurlyeq 0$ and $B\succcurlyeq 0$ in
Chan and Kwong \cite{CK}.

\begin{corollary} \label{T36}
Let $A, \, B\in {\mathbb C}^{m \times n}$ be given$.$ Then$,$ there always  exists an $X  \in \mathbb C^{n \times m}$ that satisfies
\begin{equation}
(A+B)X + X^{*}(A+B)^{*} \succcurlyeq AB^* + BA^*.
\label{342}
\end{equation}
The general solution and the corresponding $(A+B)X + X^{*}(A+B)^{*}$ can be written as
 \begin{align}
& X = \frac{1}{2}(A+B)^*  +  \frac{1}{2}\left(\  UU^{*} + V - V^* \right)\!(A+B)^* + F_{(A+B)}W,
 \label{343}
 \\
 & (A+B)X + X^{*}(A+B)^{*} = (A + B)(A + B)^* + (A + B)UU^{*}(A + B)^*,
 \label{V343xx}
\end{align}
where $U, \,V, \, W \in \mathbb C^{n \times n}$ are arbitrary$.$  In particular$,$
there exists an $X \in \mathbb C^{m\times m}$ such that
\begin{equation}
(A+B)X + X^{*}(A+B)^{*} \succ  AB +BA
\label{345}
\end{equation}
if and only if $r(A+B) = m.$

\end{corollary}

We next establish a group of formulas for calculating the ranks and inertias of
$AX + (AX)^{*} - C$ subject to  (\ref{32}), and use the results obtained to derive necessary
and sufficient conditions for the following two-side LMI
\begin{equation}
C \succcurlyeq \  AX + (AX)^{*}  \succcurlyeq B
\label{383rr}
\end{equation}
and their variations to hold.

\begin{theorem} \label{T38}
Let $A \in \mathbb C^{m\times n}$ and  $B, \ C \in \mathbb C_{{\rm
H}}^{m}$ be given$,$ ${\cal S}_1$ be as given in {\rm (\ref{36j})}$,$ and let
\begin{equation}
 N = \left[\!\! \begin{array}{cc} C & A
 \\ A^* & 0 \end{array} \!\!\right]\!,  \ \    K_1 = \left[\!\! \begin{array}{ccc} B & C & A
 \\ A^* & 0  & 0 \\ 0 & A^* &0 \end{array} \!\!\right]\!, \ \
 K_2 =\left[\!\! \begin{array}{ccc} B & C & A \\  A^*  & A^* & 0 \end{array} \!\!\right]\!.
\label{384}
\end{equation}
Then$,$ the extremal ranks and partial inertias of $AX + (AX)^{*} -C$ subject
  to $X \in {\cal S}_1$ are given by
\begin{align}
& \max_{X \in {\cal S}_1} r[\, AX + (AX)^{*}  - C\,]  =  \min\left\{\,
r(K_2)  - r(A), \ \ r(N) \right\},
\label{385}
\\
& \min_{X \in {\cal S}_1} r[\, AX + (AX)^{*}  - C \,]  = \max\{ \, t_1, \ \ t_2, \ \ t_4, \ \ t_4\, \},
\label{386}
\\
& \max_{X \in {\cal S}_1} i_{+}[\, AX + (AX)^{*}  - C \,]  = i_{-}(N),
\label{387}
\\
& \max_{X \in {\cal S}_1} i_{-}[\, AX + (AX)^{*}  - C \,]  = \min\left\{\,  i_{-}(\, B - C\,),  \ \   i_{+}(N) \, \right\},
\label{388}
\\
& \min_{X \in {\cal S}_1} i_{+}[\, AX + (AX)^{*}  - C \,]  =
  \max \{ r(K_2) + i_{-}(N) - r(K_1), \  r(K_2) +  i_{+}(\, B - C \, )  - r[\,  A, \, B - C \,] - r(A) \, \},
\label{389}
\\
& \min_{X \in {\cal S}_1} i_{-}[\, AX + (AX)^{*}  - C \,]  = \max \{ r(K_2) + i_{+}(N) - r(K_1), \  r(K_2) +  i_{-}(\, B - C \, )  - r[\,  A, \, B - C \,] - r(A) - m \, \},
\label{390}
\end{align}
where
\begin{align*}
& t_1  = 2r(K_2) + r(N) - 2r(K_1),
\\
& t_2   = 2r(K_2) + r(\, B - C\,) -  2r[\,A, \,  B - C\,] - 2r(A) - m,
 \\
& t_3  = 2r(K_2) + i_{+}(N) + i_{+}(\, B- C \, )  - r(A) -   r[\, A, \,B -C\,] -  r(K_1),
\\
& t_4  = 2r(K_2) + i_{-}(N) + i_{-}(\, B - C \,) - r(A) -  r[\, A, \,B -C\,] - r(K_1) -m.
\end{align*}
In consequence$,$
\begin{enumerate}
\item[{\rm(a)}] There exists an $X \in \mathbb C^{n\times m}$ such that
$C \succ  AX + (AX)^*  \succcurlyeq B$ if and only if $i_{+}(N) \succcurlyeq m$ and  $ C \succ  B.$

\item[{\rm(b)}] There exists an $X \in \mathbb C^{n\times m}$ such that
$AX + (AX)^*  \succcurlyeq B$ and $AX + (AX)^*  \succ  C$ if and only if $i_{-}(N) \succcurlyeq m.$

\item[{\rm(c)}] There exists an $X \in \mathbb C^{n\times m}$ such that
$AX + (AX)^* \succcurlyeq B$ and $AX + (AX)^*  \succcurlyeq C$  if and only if
$$
r(K_2) + i_{+}(N) - r(K_1) =0 \ \ and \ \
r(K_2) +  i_{-}(\, B - C \, )  = r[\,  A, \, B - C \,] + r(A) + m.
$$

\item[{\rm(d)}] There exists an $X \in \mathbb C^{n\times m}$ such that
$ C \succcurlyeq AX + (AX)^* \succcurlyeq B$  if and only if
$$
C \succcurlyeq B, \   i_{-}(N) = r(A) \ \ and  \ \ r(K_2) = r[\,  A, \, B - C \,] + r(A).
$$
\end{enumerate}
\end{theorem}

\noindent {\bf Proof.} \
Note from (\ref{35j}) that
\begin{equation}
AX + (AX)^{*} -C = B  -C  + (\, AU + J^{\frac{1}{2}} \,)(\, AU + J^{\frac{1}{2}} \,)^*.
\label{391}
\end{equation}
Applying Lemma \ref{T212F}(a) to (\ref{391}), we obtain
\begin{align}
&\max_{U\in {\mathbb C}^{n \times m}}\!\!\!r[\, B  - C + (\, AU + J^{\frac{1}{2}} \,)(\, AU + J^{\frac{1}{2}} \,)^*\,] \nb
\\
& = \min\left\{r[\,  A, \,  B - C, \, J^{\frac{1}{2}} \,], \  r\!\left[\!\! \begin{array}{cc}  B - C + J & A
 \\ A^* & 0 \end{array} \!\!\right]\!, \  r(\,B - C\,) + m \right\},
\label{392}
\\
&\min_{U\in {\mathbb C}^{n \times m}}\!\!\!r[\,B - C + (\, AU + J^{\frac{1}{2}} \,)(\, AU + J^{\frac{1}{2}} \,)^*\,] =
2r[\,A, \, B - C, \, J^{\frac{1}{2}}\,]  + \max\{ \, h_1, \ \ h_2, \ \ h_3, \ \ h_4 \, \},
\label{393}
\\
& \max_{U\in {\mathbb C}^{n \times m}}\!\!\!i_{+}[\, B  - C + (\, AU + J^{\frac{1}{2}} \,)(\, AU + J^{\frac{1}{2}} \,)^*\,] = \min\left\{\, i_{+} \!\left[\!\! \begin{array}{cc} B - C + J & A
 \\ A^* & 0 \end{array} \!\!\right]\!, \ \   i_{+}(\,B - C\,)  +m \, \right\},
\label{394}
\\
& \max_{U\in {\mathbb C}^{n \times m}}\!\!\!i_{-}[\, B - C + (\, AU + J^{\frac{1}{2}} \,)(\, AU + J^{\frac{1}{2}} \,)^*\,]
 = \min\left\{\, i_{-} \!\left[\!\! \begin{array}{cc}  B  - C + J & A
 \\ A^* & 0 \end{array} \!\!\right]\!, \ \   i_{-}(\, B - C\,) \, \right\}\!,
\label{395}
\\
&  \min_{U\in {\mathbb C}^{n \times m}} \!\!\!i_{+}[\, B - C + (\, AU + J^{\frac{1}{2}} \,)
(\, AU + J^{\frac{1}{2}} \,)^*\,]  \nb
\\
& = r[\,A, \,B - C, \, J^{\frac{1}{2}}\,]   + \max \left\{ i_{+}\!\left[\!\! \begin{array}{cc}  B - C  + J & A
 \\ A^* & 0 \end{array} \!\!\right] - r\!\left[\!\!\begin{array}{cccc}  B - C & A & J^{\frac{1}{2}}
 \\ A^* & 0 & 0  \end{array} \!\!\right]\!,  \  i_{+}(\, B - C \, )  - r[\,  A, \, B - C \,] \right\}\!,
\label{396}
\\
&  \min_{U\in {\mathbb C}^{n \times m}} \!\!\!i_{-}[\, B - C + (\, AU + J^{\frac{1}{2}} \,)(\, AU + J^{\frac{1}{2}} \,)^*\,]  \nb
\\
&  = r[\,A, \,B - C, \, J^{\frac{1}{2}}\,]  + \max \left\{ i_{-}\!\left[\!\! \begin{array}{cc}  B - C + J & A
 \\ A^* & 0 \end{array} \!\!\right] - r\!\left[\!\!\begin{array}{cccc}  B  - C & A & J^{\frac{1}{2}}
 \\ A^* & 0 & 0  \end{array} \!\!\right]\!,  \  i_{-}(\,B - C \,)  - r[\,  A, \, B  - C \,] - m \right\}\!,
\label{397}
\end{align}
where
\begin{align*}
&  h_1  = r\!\left[\!\! \begin{array}{cc}  B  - C + J & A
 \\ A^* & 0 \end{array} \!\!\right] - 2r\!\left[\!\! \begin{array}{ccc}  B - C & A & J^{\frac{1}{2}}
 \\ A^* & 0 & 0 \end{array} \!\!\right],
 \\
& h_2  =  r(\, B - C\,) -  2r[\,A, \,  B - C\,] - m,
\\
& h_3  =  i_{-}\!\left[\!\! \begin{array}{cc}  B - C  + J & A
 \\ A^* & 0 \end{array} \!\!\right] - r\!\left[\!\! \begin{array}{ccc}  B - C & A & J^{\frac{1}{2}}
 \\ A^* & 0  & 0 \end{array} \!\!\right] + i_{+}(\, B - C \,) - r[\,A, \,  B - C\,],
\\
& h_4 =  i_{+}\!\left[\!\! \begin{array}{cc}  B - C + J & A
 \\ A^* & 0 \end{array} \!\!\right] - r\!\left[\!\! \begin{array}{ccc}  B - C & A & J^{\frac{1}{2}}
 \\ A^* & 0  & 0 \end{array} \!\!\right] + i_{-}(\, B - C \,) - r[\,A, \,  B - C\,] -m.
\end{align*}
Simplifying the ranks and partial inertias of the block matrices in  (\ref{392})--(\ref{397}) gives
\begin{align}
r[\,  A, \,  B - C, \, J^{\frac{1}{2}} \,] & = r[\,  A, \,  B - C, \, J \,] = r[\,  A, \,  B - C, \, E_ABE_A \,]
=  r[\,  A, \,  B - C, \, BE_A \,] \nb
\\
& =  r\!\left[\!\! \begin{array}{ccc} A & B -C & B \\ 0  & 0  & A^* \end{array} \!\!\right]  - r(A)
=r\!\left[\!\! \begin{array}{ccc} B & C & A \\  A^*  & A^* & 0 \end{array} \!\!\right]  - r(A),
\label{398}
\\
i_{\pm}\!\left[\!\! \begin{array}{cc}  B - C + J & A
 \\ A^* & 0 \end{array} \!\!\right]&  = i_{\pm}\!\left[\!\! \begin{array}{cc}  B - C - E_ABE_A  & A
 \\ A^* & 0 \end{array} \!\!\right] = i_{\mp}\!\left[\!\! \begin{array}{cc} C  & A
 \\ A^* & 0 \end{array} \!\!\right]\!,
\label{399}
\\
 r\!\left[\!\! \begin{array}{ccc}  B - C & A & J^{\frac{1}{2}}
 \\ A^* & 0  & 0 \end{array} \!\!\right] & = r\!\left[\!\! \begin{array}{ccc}  B - C & A & J
 \\ A^* & 0  & 0 \end{array} \!\!\right] = r\!\left[\!\! \begin{array}{ccc}  B - C & A & E_ABE_A
 \\ A^* & 0  & 0 \end{array} \!\!\right] \nb
 \\
&  = r\!\left[\!\! \begin{array}{ccc}  B - C & A & BE_A
 \\ A^* & 0  & 0 \end{array} \!\!\right] = r\!\left[\!\! \begin{array}{ccc}  B - C & A & B
 \\ A^* & 0  & 0 \\ 0 & 0 & A^* \end{array} \!\!\right] - r(A) =  r(K_1) - r(A).
 \label{3100}
 \end{align}
 Substituting (\ref{398})--(\ref{3100}) into (\ref{392})--(\ref{397}) gives (\ref{385})--(\ref{390}).
\qquad $\Box$

\medskip

The following result can be shown similarly.

\begin{theorem} \label{T39}
Let $A \in \mathbb C^{m\times n}$ and  $B, \ C \in \mathbb C_{{\rm H}}^{m}$ be given$,$
${\cal S}_2$ be as given in {\rm (\ref{321j})}$,$ and let
\begin{equation}
N = \left[\!\! \begin{array}{cc} C & A
 \\ A^* & 0 \end{array} \!\!\right]\!,  \ \    K_1 = \left[\!\! \begin{array}{ccc} B & C & A
 \\ A^* & 0  & 0 \\ 0 & A^* &0 \end{array} \!\!\right]\!, \ \
 K_2 =\left[\!\! \begin{array}{ccc} B & C & A \\  A^*  & A^* & 0 \end{array} \!\!\right]\!.
\label{3101}
\end{equation}
Then$,$ the extremal ranks and partial inertias of $AX + (AX)^{*} -C$ subject
  to $X \in {\cal S}_2$ are given by
\begin{align}
& \max_{X \in {\cal S}_2} r[\, AX + (AX)^{*}  - C\,]  =  \min\left\{\,
r(K_2)  - r(A), \ \ r(N) \right\},
\label{3102}
\\
& \min_{X \in {\cal S}_2} r[\, AX + (AX)^{*}  - C \,]  = \max\{ \, t_1, \ \ t_2, \ \ t_4, \ \ t_4\, \},
\label{3103}
\\
& \max_{X \in {\cal S}_2} i_{+}[\, AX + (AX)^{*}  - C \,]  = \min\left\{\,  i_{+}(\, B - C\,),  \ \
 i_{-}(N) \, \right\},
\label{3104}
\\
& \max_{X \in {\cal S}_2} i_{-}[\, AX + (AX)^{*}  - C \,]  = i_{+}(N),
\label{3105}
\\
&  \min_{X \in {\cal S}_2} i_{+}[\, AX + (AX)^{*}  - C \,]   =
  \max \{ r(K_2) + i_{-}(N) - r(K_1),  \  r(K_2) +  i_{+}(\, B - C \, )  - r[\,  A, \, B - C \,] - r(A) -m \, \},
\label{3106}
\\
& \min_{X \in {\cal S}_2} i_{-}[\, AX + (AX)^{*}  - C \,]  = \max \{ r(K_2) + i_{+}(N) - r(K_1), \  r(K_2) +  i_{-}(\, B - C \, )  - r[\,  A, \, B - C \,] - r(A) \, \},
\label{3107}
\end{align}
where
\begin{align*}
& t_1  = 2r(K_2) + r(N) - 2r(K_1),
\\
& t_2   = 2r(K_2) + r(\, B - C\,) -  2r[\,A, \,  B - C\,] - 2r(A) - m,
 \\
& t_3  = 2r(K_2) + i_{-}(N) + i_{-}(\, B - C \, )  - r(A) -   r[\, A, \,B - C\,] -  r(K_1),
\\
& t_4  = 2r(K_2) + i_{+}(N) + i_{+}(\, B - C \,) - r(A) -  r[\, A, \,B - C\,] - r(K_1) -m.
\end{align*}
In consequence$,$ the following hold$.$
\begin{enumerate}
\item[{\rm(a)}] There exists an $X \in \mathbb C^{n\times m}$ such that
$C \prec  AX + (AX)^*  \preccurlyeq B$ if and only if $i_{-}(N) \geqslant m$ and  $ B \succ  C.$

\item[{\rm(b)}] There exists an $X \in \mathbb C^{n\times m}$ such that
$AX + (AX)^*  \preccurlyeq B$ and $AX + (AX)^*  \prec  C$ if and only if $i_{+}(N) \geqslant m.$

\item[{\rm(c)}] There exists an $X \in \mathbb C^{n\times m}$ such that
$AX + (AX)^* \preccurlyeq B$ and $AX + (AX)^*  \preccurlyeq C$  if and only if
$$
r(K_2) + i_{-}(N) - r(K_1) =0 \ \ and \ \
r(K_2) +  i_{+}(\, B - C \, )  = r[\,  A, \, B - C \,] + r(A) + m.
$$

\item[{\rm(d)}] There exists an $X \in \mathbb C^{n\times m}$ such that
$ C \preccurlyeq AX + (AX)^* \preccurlyeq B$  if and only if
$$
C \preccurlyeq B, \   i_{+}(N) = r(A) \ \ and  \ \ r(K_2) = r[\,  A, \, B - C \,] + r(A).
$$
\end{enumerate}
\end{theorem}

\section{The extremal ranks and inertias of $A - BX - XB^*$ subject to $BXB^*=C$}
\renewcommand{\theequation}{\thesection.\arabic{equation}}
\setcounter{section}{6}
\setcounter{equation}{0}

We first establish in this section a group of formulas for calculating the extremal ranks and inertias of
$A - BX - XB^*$ subject to $BXB^*=C$, and use the formulas to  characterize the existence of Hermitian
matrix $X$ satisfying the following inequalities
\begin{align}
BX + XB^* \succcurlyeq \,(\succ, \, \preccurlyeq, \, \prec ) \, A \  \ \ {\rm s.t.}  \ \  BXB^*=C
\label{71}
\end{align}
in the L\"owner partial ordering.

\begin{theorem}\label{T71}
Let $A, \ C \in {\mathbb C}_{{\rm H}}^{m}$ and $B\in {\mathbb C}^{m \times m}$ be given, and
assume that $BXB^*=C$ has a Hermitian solution$.$
Also let
\begin{align}
& M_1 = \left[\!\!\begin{array}{ccc} A & B \\ B^{*} & 0\end{array}\!\!\right]\!, \ \
M_2 = \left[\!\!\begin{array}{ccc} BA & B^2 \\ B^{*} & 0\end{array}\!\!\right]\!, \ \ M_3 = BAB^* - BC - CB^*,
\label{72}
\\
& {\cal S} = \{ X \in \mathbb C_{{\rm H}}^{m} \ | \  BXB^*=C \}.
\label{73}
\end{align}
Then$,$ the following hold$.$
\begin{enumerate}
\item[{\rm(a)}]  The maximal rank of $A - BX - XB^*$ subject to $X \in {\cal S}$ is
\begin{align}
 \max_{X \in {\cal S}}r(\, A - BX - XB^*\,) = \min\!\left\{ \,  m +  r[\, B^2, \ C- BA \,] - r(B), \
r(M_1), \  2m - 2r(B) + r(M_3) \, \right\}.
\label{74}
\end{align}

\item[{\rm(b)}]  The minimal rank of  $A - BX - XB^*$ subject to $X \in {\cal S}$  is
\begin{align}
& \min_{X \in {\cal S}}r(\, A - BX - XB^*\,) = \max\{ \, s_1, \ \ s_2, \ \ s_3, \ \ s_4 \, \},
\label{75}
\end{align}
where
\begin{align*}
& s_1  = 2r[\, B^2, \ C- BA \,] + r(M_1) - 2r(M_2),
\\
& s_2  = 2r[\, B^2, \ C- BA \,] + r(M_3) - 2r[\, B^2, \, BAB^* - CB^*\,],
\\
& s_3  = 2r[\, B^2, \ C- BA \,] + i_{+}(M_1)  - r(M_2) +  i_{-}(M_3) - r[\, B^2, \, CB^* - BAB^*\,],
\\
& s_4  = 2r[\, B^2, \ C- BA \,] + i_{-}(M_1) - r(M_2)  + i_{+}(M_3) - r[\, B^2, \, CB^* - BAB^*\,].
\end{align*}

\item[{\rm(c)}]  The  maximal partial inertia of  $A - BX - XB^*$ subject to $X \in {\cal S}$  is
\begin{align}
 \max_{X \in {\cal S}} i_{\pm}(\, A - BX - XB^*\,) = \min\!\left\{ i_{\pm}(M_1), \
 m - r(B) + i_{\pm}(M_3)  \right\}.
\label{76}
\end{align}

\item[{\rm(d)}] The  minimal partial inertia of  $A - BX - XB^*$ subject to $X \in {\cal S}$  is
\begin{align}
\!\!\!\!\!\min_{X \in {\cal S}}i_{\pm}(\, A - BX - XB^*\,) = r[\, B^2, \ C- BA \,]  + \max\!\left\{ i_{\pm}(M_1)- r(M_2), \  i_{\pm}(M_3) - r[\, B^2, \, CB^* - BAB^* \,]  \right\}.
\label{77}
\end{align}
\end{enumerate}
\end{theorem}

\noindent  {\bf Proof.} \ From Lemma 2.3(b), the general Hermitian solution of $BXB^* = C$ can
be expressed as
\begin{align}
X = B^{\dag}C(B^{\dag})^* + F_BV + V^*F_B,
 \label{78}
\end{align}
where the matrix $V$ is arbitrary. Substituting it  into  $A - BX - XB^*$  yields
\begin{align}
 A - BX - XB^* = A -  C(B^{\dag})^* - B^{\dag}C -  BV^*F_B - F_BVB^*.
 \label{79}
\end{align}
Define
\begin{align}
\phi(V) = G - BV^*F_B - F_BVB^*,
 \label{710}
\end{align}
where $G = A -  C(B^{\dag})^* - B^{\dag}C$.  Applying  (\ref{129})--(\ref{132}) to (\ref{710}) yields
\begin{align}
& \max_{V\in \mathbb C^{m \times m}}r[\phi(V)]  = \min \left\{ r[\,G, \, B,\, F_B\,],
 \ \ r\!\left[\begin{array}{cc} G & B
\\B^{*} & 0
\end{array}\right]\!, \ \ r\!\left[\!\begin{array}{cc} G & F_B
\\ F_B & 0
\end{array}\!\right] \right\}\!,
\label{711}
\\
& \min_{V\in \mathbb C^{m \times m}}r[\phi(V)]  =
2r[\,G, \, B,\, F_B\,] + \max\{\, s_{+} + s_{-}, \ \
t_{+} + t_{-}, \ \ s_{+} + t_{-}, \ \ s_{-} + t_{+} \, \},
\label{712}
\\
& \max_{V\in \mathbb C^{m \times m}}i_{\pm}[\phi(V)] = \min\!\left\{i_{\pm}\!\left[\!\begin{array}{ccc} G
 & B  \\  B^{*}  & 0
   \end{array}\!\right], \ \ i_{\pm}\!\left[\!\begin{array}{ccc} G & F_B  \\  F_B  & 0
   \end{array}\!\right] \right\}\!,
\label{713}
\\
& \min_{V\in \mathbb C^{m \times m}}i_{\pm}[\phi(V)] = r[\,G, \, B,\, F_B\,] + \max\{\, s_{\pm}, \ \
t_{\pm} \, \},
\label{714}
\end{align}
where
\begin{align*}
s_{\pm}  = i_{\pm}\!\left[\!\!\begin{array}{cc} G & B \\ B^{*} & 0\end{array}\!\!\right]
 - r\!\left[\begin{array}{ccc} G & B  & F_B \\ B^{*} & 0 & 0
\end{array}\!\!\right]\!, \ \
 t_{\pm}  =i_{\pm}\!\left[\!\!\begin{array}{cc} G & F_B \\ F_B & 0\end{array}\!\!\right]
 - r\!\left[\begin{array}{ccc} G & B  & F_B \\ F_B & 0 & 0
\end{array}\!\!\right]\!.
\end{align*}
Applying (\ref{215}) and (\ref{216}), and
simplifying by elementary matrix operations and congruence matrix operations, we obtain
\begin{align}
r[\,G, \, B,\, F_B\,] & =r\!\left[\!\!\begin{array}{ccc}
A -  C(B^{\dag})^* - B^{\dag}C & B  & I_m \\ 0 & 0 & B \end{array}\!\!\right] - r(B) \nonumber
\\
& =r\!\left[\!\!\begin{array}{ccc} 0 & 0 & I_m \\ -BA + C  & -B^2 & 0 \end{array}\!\!\right] - r(B)
 = m +  r[\, B^2, \ C- BA \,] - r(B),
\label{715}
\\
 r\!\left[\!\!\begin{array}{ccc} G & B  & F_B \\ B^{*} & 0 & 0
\end{array}\!\!\right] & = r\!\left[\!\!\begin{array}{ccc} A & B  & I_m \\ B^{*} & 0 & 0 \\ 0 & 0 & B
\end{array}\!\!\right] - r(B) \nonumber
\\
& = r\!\left[\!\!\begin{array}{ccc} 0 & 0  & I_m \\ B^{*} & 0 & 0 \\ -BA  & -B^2 & 0
\end{array}\!\!\right] - r(B)
 = m + r\!\left[\begin{array}{cc} BA  & B^2 \\ B^* & 0 \end{array}\!\!\right] - r(B),
\label{716}
\\
r\!\left[\begin{array}{ccc} G  & B & F_B \\ F_B & 0 & 0
\end{array}\!\!\right] & = r\!\left[\begin{array}{cccc} A - B^{\dag}C  & B & I_m  & 0 \\
I_m & 0 & 0 & B^* \\ 0 & 0 & B & 0
\end{array}\!\!\right] - 2r(B)  \nonumber
\\
& = r\!\left[\!\!\begin{array}{cccc} 0 & 0 & I_m  & 0 \\
I_m & 0 & 0 & B^* \\  - BA + C  & -B^2 & 0 & 0
\end{array}\!\!\right] - 2r(B) \nonumber
\\
& = r\!\left[\!\!\begin{array}{cccc} 0 & 0 & I_m  & 0 \\
I_m & 0 & 0 & 0 \\  0 & -B^2 & 0 &  BAB^* - CB^*
\end{array}\!\!\right] - 2r(B) \nonumber
\\
& = 2m + r[\, B^2, \, CB^* -BAB^* \,] - 2r(B),
\label{717}
\\
i_{\pm}\!\left[\!\!\begin{array}{cc} G & B \\ B^{*} & 0\end{array}\!\!\right]
& = i_{\pm}\!\left[\!\!\begin{array}{ccc} A -  C(B^{\dag})^* - B^{\dag}C & B \\ B^{*} & 0\end{array}\!\!\right]
=i_{\pm}\!\left[\!\!\begin{array}{ccc} A & B \\ B^{*} & 0\end{array}\!\!\right],
\label{718}
\\
i_{\pm}\!\left[\!\!\begin{array}{cc} G
 & F_B  \\  F_B  & 0 \end{array}\!\!\right] & = r(F_B) + i_{\pm}(BGB^*)
  =  m - r(B) +  i_{\pm}(\,BAB^* - BC - CB^*\,).
\label{719}
\end{align}
Hence,
\begin{align}
& s_{\pm}   = i_{\pm}\!\left[\!\!\begin{array}{cc} A & B \\ B^{*} & 0\end{array}\!\!\right]
  - r\!\left[\begin{array}{cc} BA  & B^2 \\ B^* & 0 \end{array}\!\!\right] -  m + r(B),
\label{720}
\\
& t_{\pm}  = i_{\pm}(\,BAB^* - BC - CB^*\,) - r[\, B^2, \, BAB^* - CB^* \,] - m + r(B).
\label{721}
\end{align}
Substituting  (\ref{715})--(\ref{721}) into  (\ref{711})--(\ref{714})  yields (\ref{74})--(\ref{77}).
 \qquad $\Box$

\begin{corollary}\label{T72}
Let $A, \ C \in {\mathbb C}_{{\rm H}}^{m}$ and $B\in {\mathbb C}^{m \times m}$ be given, and
assume that $BXB^*=C$ and $BX + (BX)^*  = A$ are  consistent$,$ respectively$.$
Also let ${\cal S}(X)$ be as given in {\rm (\ref{73})}$.$ Then$,$ the following hold$.$
\begin{align}
\!\!\!\! \max_{X \in {\cal S}}\!r(\, A - BX - XB^*\,) &  = \min\{ \,  m +  r[\, B^2, \ C- BA \,] - r(B), \ 2r(B), \ 2m - 2r(B) + r(\,BAB^* - BC - CB^*\,) \, \},
\label{722}
\\
\!\!\!\!\min_{X \in {\cal S}}\!r(\, A - BX - XB^*\,) & = \max\{ \, s_1, \ \ s_2, \ \ s_3, \ \ s_4 \, \},
\label{723}
\\
\!\!\!\!\max_{X \in {\cal S}}i_{\pm}(\, A - BX - XB^*\,) & = \min\!\left\{ r(B), \ \
  m - r(B) + i_{\pm}(\,BAB^* - BC - CB^*\,) \right\},
\label{724}
\\
\!\!\!\! \min_{X \in {\cal S}}i_{\pm}(\, A - BX - XB^*\,) & = r[\,  B^2, \ C- BA\,]  +  \max\{ - r(B^2),  \ i_{\pm}(\,BAB^* - BC - CB^*\,) - r[\, B^2, \, CB^* - BAB^* \,] \},
\label{725}
\end{align}
where
\begin{align*}
& s_1  = 2r[\,  B^2, \ C- BA \,] - 2r(B^2),
\\
& s_2  = 2r[\,  B^2, \ C- BA \,] + r(\, BAB^* - BC - CB^*\,) - 2r[\, B^2, \, CB^* - BAB^*\,],
\\
& s_3 = 2r[\,  B^2, \ C- BA \,] - r(B^2) +  i_{-}(\,BAB^* - BC - CB^*\,) - r[\, B^2, \, CB^* - BAB^* \,],
\\
& s_4  = 2r[\,  B^2, \ C- BA \,]  - r(B_2)  + i_{+}(\,BAB^* - BC - CB^*\,) - r[\, B^2, \, CB^* - BAB^*\,].
\end{align*}
\end{corollary}

\begin{corollary}\label{T73}
The pair of matrix equations
\begin{align}
BX + XB^* = A, \ \  BXB^*=C
\label{726}
\end{align}
 have a common solution $X\in {\mathbb C}_{{\rm H}}^{m}$ if and only if
\begin{align}
{\mathscr R}(C) \subseteq {\mathscr R}(B),  \ \  r\!\left[ \!\!\begin{array}{cc} A & B  \\  B^* & 0
\end{array}\!\! \right] = 2r(B), \ \ BC + CB^* = BAB^*, \ \  {\mathscr R}( \,C - BA \, )
\subseteq {\mathscr R}( B^2).
\label{727}
\end{align}
In that case$,$ the general common Hermitian solution of {\rm (\ref{726})}
can be expressed as
\begin{align}
X = X_0 +  [\, F_A, \  0 \, ] F_G UE_H \left[\!\! \begin{array}{c} I_m  \\ 0
 \end{array}\!\! \right] + [\,  0, \  I_m \, ] F_G UE_H
 \left[\!\! \begin{array}{c} 0 \\ E_B \end{array} \!\!\right] + F_BSF_B,
\label{728}
\end{align}
where $ X_0$ is a special  common Hermitian solution of {\rm (\ref{726})}$,$  $G = [\, F_A, \ -A \, ],
\ H = \left[\! \begin{array}{c} B  \\ E_B \end{array} \!\right]\!,$  $U $ and $S$
are arbitrary$.$ In particular$,$ {\rm (\ref{726})} has a unique common  solution
if and only if $B$ is nonsingular and  $BC + CB^* = BAB^*$.
In the case$,$ the unique common solution is $ X = B^{-1}C(B^{-1})^*.$
\end{corollary}

\noindent {\bf Proof.} \  Suppose first that the two equations  in (\ref{726}) have a common solution.
This implies that $ BX + (BX)^* = A$ and $ BXB^* = C$ are consistent, respectively.
In this case, setting (\ref{722}) equal to zero leads to (\ref{727}).

We next show that under (\ref{727}), the two equations in (\ref{726})
have a common solution and their general common solution can be written as
(\ref{728}). Substituting (\ref{78}) into the first equation
 in (\ref{726}) yields
\begin{align}
F_BVB^*  + BV^*F_B  = A -  C(B^{\dag})^* - B^{\dag}C.
\label{729}
\end{align}
Solving for $V$ in (\ref{729}) by Lemma 6.1, we obtain the general solution
\begin{eqnarray*}
V  =  V_0 + [\, I_m, \, 0\, ]F_GWE_H\!\left[\!\! \begin{array}{c} I_m
 \\ 0 \end{array} \!\!\right] - [\,  0, \, I_m \, ]E_HW^{*}F_G\!\left[\!\! \begin{array}{c} 0 \\ I_m
\end{array} \!\!\right] + B^{\dag}BW_1 + W_2F_B,
\end{eqnarray*}
where $V_0$ is a special solution of (\ref{723})$,$ $G = [\, F_B, \, B \,],$
$H = \left[\!\!\begin{array}{c} B^* \\ F_B
\end{array} \!\!\right]\!,$ $W$, $W_1$ and $W_2$ are arbitrary.
Substituting this $V$ into (\ref{78}) yields
\begin{align*}
X & =  B^{\dag}C(B^{\dag})^* + F_BV + V^*F_B
\\
& = B^{\dag}C(B^{\dag})^* + F_BV_0 + V^*_0F_B  +
[\,  I_m, \  0 \, ]F_G UE_H \left[\! \begin{array}{c} I_n  \\ 0 \end{array}\! \right] +
F_AS_2E_B +  \  [\,  0,  \  I_m \, ] F_GUE_H \left[\! \begin{array}{c} 0 \\
 I_n \end{array}\! \right] + F_AT_1E_B,
\end{align*}
which can simply be written in the form of (\ref{728}).  \qquad $ \Box $

\begin{corollary}\label{T74}
Let $A, \ C \in {\mathbb C}_{{\rm H}}^{m}$ and $B\in {\mathbb C}^{m \times m}$ be given, and
assume that $BXB^*=C$ is consistent$.$ Also let $M_1,$  $M_2$ and $M_3$ be of the forms in {\rm (\ref{72})}.
Then$,$ the following hold$.$
\begin{enumerate}
\item[{\rm(a)}] There exists an $X\in {\mathbb C}_{{\rm H}}^{m}$ such that
\begin{align}
BX + XB^* \preccurlyeq A \ \  and \ \ \  BXB^*=C
\label{730}
\end{align}
if and only if
\begin{align}
r[\, B^2, \ C- BA \,] + i_{-}(M_1)- r(M_2) = 0, \ \
r[\, B^2, \ C- BA \,]  + i_{-}(M_3) - r[\, B^2, \, CB^* - BAB^*] =0.
\label{731}
\end{align}

\item[{\rm(b)}] There exists an $X\in {\mathbb C}_{{\rm H}}^{m}$ such that
\begin{align}
BX + XB^* \prec  A \ \  and \ \ \  BXB^*=C
\label{732}
\end{align}
if and only if
\begin{align}
 i_{+} (M_1) \geqslant m, \ \ and  \ \ i_{+}(\,BAB^* - BC - CB^*\,) = r(B).
\label{733}
\end{align}

\item[{\rm(c)}] There exists an $X\in {\mathbb C}_{{\rm H}}^{m}$ such that
\begin{align}
BX + XB^* \succcurlyeq A \ \  and \ \ \  BXB^*=C
\label{734}
\end{align}
if and only if
\begin{align}
r[\, B^2, \ C- BA \,] + i_{+}(M_1)- r(M_2) = 0, \ \
r[\, B^2, \ C- BA \,]  + i_{+}(M_3) - r[\, B^2, \, CB^* - BAB^*] =0.
\label{735}
\end{align}

\item[{\rm(d)}] There exists an $X\in {\mathbb C}_{{\rm H}}^{m}$ such that
\begin{align}
BX + XB^* \succ  A \ \  and \ \ \  BXB^*=C
\label{736}
\end{align}
if and only if
\begin{align}
i_{-} (M_1) \succcurlyeq  m \ \ and  \ \ i_{-}(\,BAB^* - BC - CB^*\,) = r(B).
\label{737}
\end{align}
\end{enumerate}
\end{corollary}

\section{Concluding remarks}
\renewcommand{\theequation}{\thesection.\arabic{equation}}
\setcounter{section}{8}
\setcounter{equation}{0}

In the previous sections, we showed that the three LMIs  of fundamental types in  (\ref{10b})--(\ref{10d})
can equivalently be converted to some quadratic matrix equations. Through the quadratic matrix equations
and a variety of known results on linear and quadratic matrix equations, we established necessary
and sufficient conditions for these LMIs to be feasible and obtained general solutions of these LMIs.
Since the results obtained in the previous sections
are represented in closed form by using the ranks, inertias and ordinary operations of the given matrices
and their generalized inverses, they can be easily used to approach various problems related to these basic
LMIs in matrix theory and applications. In particular, they can be used to solve mathematical programming and
optimization problems subject to LMIs in the (\ref{10b})--(\ref{10d}).

Based on the results in the previous sections, it is not hard to establish analytical solutions of the following constrained LMIs:
\begin{enumerate}
\item[{\rm(a)}] $AXB \succcurlyeq \,(\succ, \, \preccurlyeq, \, \prec ) \, C$ subject to $PX=Q$ and/or $XR = S;$

\item[{\rm(b)}] $AXA^{*} \succcurlyeq \,(\succ, \, \preccurlyeq, \, \prec ) \,B$ subject to $PX=Q$ and $X = X^*$, or $PXP^* = Q$ and $X  = X^*;$

\item[{\rm(c)}] $AX + (AX)^{*} \succcurlyeq \,(\succ, \, \preccurlyeq, \, \prec ) \,  C$ subject to $PX=Q$.
\end{enumerate}
The results obtained will sufficiently meet people's curiosity about analytical solutions of LMIs.
 In addition, the work in this paper will also motivate finding possible analytical solutions of some general LMIs, such as,
\begin{enumerate}
\item[{\rm(d)}] $AX + YB \succcurlyeq \,(\succ, \, \preccurlyeq, \, \prec ) \, C$;

\item[{\rm(e)}] $AXA^{*} + BYB^{*} \succcurlyeq \,(\succ, \, \preccurlyeq, \, \prec ) \,C$;

\item[{\rm(f)}] $AXA^{*} \succcurlyeq \,(\succ, \, \preccurlyeq, \, \prec ) \, B$ and $CXC^{*} \succcurlyeq \,(\succ, \, \preccurlyeq, \, \prec ) \, D$;

\item[{\rm(g)}] $AXB + (AXB)^{*}  \succcurlyeq \,(\succ, \, \preccurlyeq, \, \prec ) \,  C$,
\end{enumerate}
which are equivalent to the following linear-quadratic matrix equations:
\begin{enumerate}
\item[{\rm(d1)}] $AX + YB = C \pm UU^*$;

\item[{\rm(e1)}] $AXA^{*} + BYB^{*}  = C \pm UU^*$;

\item[{\rm(f1)}] $AXA^{*} =  B  \pm UU^*$ and $CXC^{*} = D \pm VV^*$;

\item[{\rm(g1)}] $AXB + (AXB)^{*}  = C \pm UU^*.$
\end{enumerate}
A special case of (g) for  $C \succcurlyeq 0$ was solved  in \cite{TV}.

In system and control theory, minimizing or maximizing the rank of a variable matrix over a
set defined by matrix inequalities in the L\"owner partial ordering is referred to as
a rank minimization or maximization problem, and is denoted collectively by RMPs.
 The RMP now is known to be NP-hard in general case, and a satisfactory
characterization of the solution set of a general RMP is currently not available.
Notice from the results in this paper that for some types of matrix inequality
in the L\"owner partial ordering, their general solutions can be written in closed form
by using the given matrices and their generalized inverses in the inequalities.
Hence, it is expected that the results in this paper can used to solve certain RMPs.
These further developments are beyond the scope of the present
paper and will be the subjects of separate studies.

After a half century's development of the theory of generalized inverses of matrices,
people now are widely using generalized inverses of matrices to solve a huge amount of problems in matrix theory
and applications. In particular, one can utilize them to represent solutions of matrix equations and
inequalities.

Since linear algebra is a successful theory with essential applications in most scientific fields,
the methods and results in matrix theory are prototypes of many concepts and content in other advanced branches of mathematics.
 In particular,  matrix equations and matrix inequalities in the
 L\"owner partial ordering, as well as generalized inverses of matrices were sufficiently extended to
 their counterparts for operators in a Hilbert space, or elements in a ring with involution,
 and their algebraic properties were extensively studied in the literature.
   In most cases, the conclusions on the complex matrices and their counterparts in general algebraic settings
    are analogous. Also, note that the results in this paper are derived from
 ordinary algebraic operations of the given matrices and their  generalized inverses. Hence,
 it is no doubt that most of the conclusions in this paper can trivially be extended to
  the corresponding  equations and inequalities for linear operators on a Hilbert space or elements
  in a ring with involution.

\end{document}